\documentclass[openbib,openany,a4paper,twoside,12pt]{book}
\usepackage{latexsym,graphics,makeidx,amsmath,amssymb}
\pdfoutput=1
\makeindex
\newtheorem{thm}{Theorem}[section]
\newtheorem{lem}{Lemma}[section]
\newtheorem{cor}{Corollary}[section]
\newtheorem{prop}{Proposition}[section]
\newtheorem{defi}{Definition}[section]
\newtheorem{rem}{Remark}[section]
\newtheorem{exam}{Example}[section]
\title{ON DILATIONS AND TRANSFERENCE FOR CONTINUOUS ONE-PARAMETER 
SEMIGROUPS OF POSITIVE CONTRACTIONS ON $\mathcal{L}^p$-SPACES.}
\author{GERO FENDLER}
\date{}
\begin{document}
\maketitle 
\vspace{1cm}
%
%
\thanks{\bf Acknowledgements : \rm 
These are notes of some lectures 
held at the Institute of Mathematics of the University of Wroc{\l}aw 
in November and December 1996. There the lectures were presented
under the title ``Functional calculus of operators on Hilbert spaces
and on other classes of Banach spaces''.
\par
The author is grateful to Professor M.~Bo\.{z}ejko for having had 
the possibility to present the material at his seminar. He thanks the audience
for their interest.
\par
We thank the referee for his comments. They were quite useful
for  improving the 
English as well as the mathematical presentation of the contents of this note.
\par 
Last not least we take this chance to thank Ms Georg from the mathematical institute of the
University of the Saarland for putting to \LaTeX parts of the
manuscript.\\[10ex]
{\bf \large{Pre-Introduction}}\\[5ex]
During the course held at the Wroc{\l}aw University things developed
a bit different from what was originally planned. Following the appearing interests
we partially changed the subject and now the title too, according to
what has been presented.
\par
We still found it worthwhile to include some remarks on functional
calculus on Hilbert spaces which now are presented in the introduction.
Readers interested only in one-parameter semigroups on
$\mathcal{L}^p$-spaces may thus leave out the introduction. Those 
who are only interested in Hilbert spaces may read the introduction and skip
the rest.\\}
\def\empty{}
\newcommand{\RR}{{\mathbb R}}
\newcommand{\NN}{{\mathbb N}}
\newcommand{\CC}{{\mathbb C}}
\newcommand{\Lp}[2]{L^p(#1,#2)}
\newcommand{\Lq}[2]{L^q(#1,#2)}
\newcommand{\Ulpr}[2]{\prod_{\mathfrak{#1}} #2  }
\newcommand{\Ulp}[2]{\prod_{{#1}} #2  }
\newcommand{\lp}[2]{l^{p}_{#1}(#2)}
\newcommand{\abs}[1][.]{\left| {#1} \right|}
\newcommand{\norm}[2][]{\left\| \, {#2} \, \right\|_{#1}}
\newcommand{\matr}[3]{%
\def\test{#2} \ifx \empty \test \left(#1_{i,j}\right)_{i,j=1}^{#3}%
\else \left(#2(#1_{i,j})\right)_{i,j=1}^{#3} \fi}%
\newcommand{\matrpol}[3][p]{%
\def\test{#2} \ifx \empty \test \left(#1_{i,j}\right)_{i,j=1}^{#3}
\else \left(#1_{i,j}(#2)\right)_{i,j=1}^{#3} \fi}%
\newcommand{\matrun }[5][k]{%
\def\test{#2} \ifx \empty \test \left(#2#1_{i,j}#3\right)_{i,j=1}^{#5}
\else \left(#3#1_{i,j}(#2)#4\right)_{i,j=1}^{#5} \fi}%
\newcommand{\lpnorm}[3][\infty]{\left(\sum_{#2 = 1}^{#1}{|#3|^{p}}\right)^{1/p}}
\newcommand{\cbnorm}[2][]{\left| \! \left| \! \left| \, #2 \, 
\right| \! \right| \! \right|_{#1}}
\newcommand{\qed}{\hspace*{\fill} \vspace{5mm}%
\ensuremath{\underline{~~~~\Box}}\\}
\newcommand{\qedm}{\ensuremath{\underline{~~~~\Box}}}
\newcommand{\BLIM}{\mbox{\tiny \rm LIM\normalsize}}
\newcommand{\id}{\mbox{\rm id}}
%
\tableofcontents
\chapter[Introduction]{Introduction}
\index{polynomial!calculus}%
\section[Norm Estimates for Polynomial Calculus]{Norm Estimates for
Polynomial Calculus}
\label{sec:normest}
Let $ H $ be a Hilbert space,
$ T : H \longrightarrow H $ 
\index{linear contraction}%
\index{contraction}%
a linear contraction, i.e. $ T $ is a linear map of norm at most one.
In 1951 
\index{inequality}%
\index{inequality!von~Neumann}%
\index{von~Neumann}%
J.~von~Neumann~\cite{Neumann_51} proved that for any 
\index{polynomial}%
polynomial
\begin{displaymath}
p(z) = \sum_{n=0}^N a_n z^n ,  \quad z \in \mathbb{C}
\end{displaymath}
the inequality
\begin{displaymath}
\norm{p(T)} \quad \leq \quad \sup_{\abs[z] \leq 1} \abs[p(z)]
\end{displaymath}
holds true.
Here, of course, we denote
$ p(T) \: = \: \sum_{n=1}^{N} \: a_n \: T^n $
and consider it as an element of the bounded linear operators 
$ B(H) $ of $ H. $
\par
Let me mention a short proof of this. First, if $ U \in B(H) $
is a unitary operator, thus a very special contraction, then
\begin{eqnarray*}
\norm{p(U)} ^{2} & = & \norm{p(U)^{*}p(U)} \\
              &  = & r(p(U)^{*}p(U)) \\
              & = & \sup \{\: \abs[ \lambda ]: \: \lambda \in \sigma (p(U)^*p(U)) \},
\end{eqnarray*}
where $ \sigma (A) \: = \: \{\: \lambda \: : \: 
\lambda - A \: \notin \: \mbox{Inv}(B(H)) \: \} $ denotes the 
\index{spectrum}%
\index{spectrum!operator}
spectrum of 
the linear operator $ A $ and $ r(A) $ its 
\index{spectral radius}%
spectral radius.
\par
On the other hand, 
\begin{eqnarray*}
p(U)^* p(U) & = & \sum_{n \in \mathbb{N}}^N \overline{a}_n (U^*)^n \cdot
                 \sum_{n \in \mathbb{N}}^N a_n U^n \\
           & = & \sum_{l=-N}^N \sum_{\{n,l \,: n+k=l \} } \: \overline{a}_n a_k \: U^l \\
           & = & F(U),
\end{eqnarray*}
where 
$ F(z)=\sum_{l=-N}^N \sum_{\{n,l \, : n+k=l \} } \: \overline{a} _n a_k \: z^l , \:
z \in \mathbb{C} \backslash \{0 \} $ is a rational function which is holomorphic
in a neighbourhood of 
$ \sigma (U) \subset \{\: z \in \mathbb{C} : \: \abs[z] = 1 \}$.
\par
By the spectral 
\index{theorem!spectral mapping}%
mapping theorem
$ \sigma (F(U)) \subset F(\sigma (U)) $,
and thus 
\begin{eqnarray*}
 r(p(U)^{*}p(U)) & = & \sup \{\: \abs[ \lambda ] \, : \: \lambda \in \sigma (F(U)) \} \\
              &  = & \sup \{\: \abs[ F(\lambda ) ] \, : \: \abs[ \lambda ] = 1 \} \\
              & = & \sup \{\: \abs[ \overline{(p(z)}p(z)] \, : \: \abs[z] = 1 \} \\
              & = & \sup \{\abs[ \, p(z) \, ] \, :\, \abs[z] = 1 \} ^2,
\end{eqnarray*}
which establishes von~Neumann's inequality in this case.
\qed
\par
Now, Foias and Sz.~Nagy \cite{NagyFoia_70}  constructed a unitary dilation for an arbitrary contraction:
\begin{thm}[Sz.~Nagy-Foias]
\label{thm:unitdil}
\index{theorem!Sz.~Nagy-Foias}%
\index{dilation}%
\index{dilation!unitary}%
Given a contraction
$T \in B(H)$
there exists a Hilbert space 
$ K $
together with a unitary
$U \in B(K)$,
an isometric embedding 
$ I\,: \, H \rightarrow K $
and an orthogonal projection
$ P \, :\, K \rightarrow H $
such that for 
$ k = 0,1,\dots $
\begin{eqnarray*}
T^k & = & P \circ U^k \circ I.
\end{eqnarray*}
\end{thm}
Using this theorem  we obtain, for any contraction 
$ T, $
\begin{eqnarray*}
p(T) & = & \sum_{n=0}^N \: a_n \, T^n \\
     & = & P \circ p(U) \circ I,
\end{eqnarray*}
where
$P, \, U,\, \mbox{and }I $
are as in the theorem. Thus,
\begin{eqnarray*}
\norm{p(T)} & \leq & \norm{p(U)} \\
            & \leq & \sup\{\, \abs[p(z)] \, : \, \abs[z] = 1 \},
\end{eqnarray*}
which is von~Neumann's inequality in the general case.
\qed
\par
It is clear, that
$T' \in B(H)$
is 
\index{polynomially bounded}%
\index{bounded!polynomially}%
polynomially bounded if it is
\index{similar!contraction}%
similar to a contraction
$T$
by a bounded operator
$ S \in \mbox{Inv}(B(H)), $
i.e.\ if
$T' = S^{-1} T S $.
Since in this case for any polynomial 
$p$:
$p(T')=S^{-1}p(T)S$, 
and hence
\begin{eqnarray*}
\norm{p(T')} & \leq & \norm{S^{-1}p(T)S} \\
               & \leq & \norm{S^{-1}} \, \norm{S} \, \norm[\infty]{p(T)}.
\end{eqnarray*}
\par
A natural question is whether the converse holds true too. 
This is known as the 
\index{Halmos}%
Halmos
\index{problem!Halmos}%
problem \cite{Halmos_70}. But it has some history:
\index{Sz.~Nagy}%
Sz.~Nagy~\cite{Nagy_47}~\cite{Nagy_59} in 1947 and 1959 respectively proved
\begin{itemize}
\item {if $ T \in B(H)  $ is invertible %
and $ T $ and its inverse are
\index{power-bounded}%
power-bounded, that is, if
\[ \sup \{\, \norm{T^n} \, : \, n \in \mathbb{Z} \, \} < \infty, \]
then $ T $ is similar to a unitary operator.}
\item {if $T \in B(H) $ is compact and power-bounded, for 
\index{power-bounded!positive potencies}%
positive 
potencies only, i.e.\ if
\[ \sup \{\, \norm{T^n} \, : \, n \in \mathbb{N} \, \} < \infty, \]
then $ T $ is similar to a contraction.}
\end{itemize}
His question then was:\\
If a general operator $ T \in B(H) $ is power-bounded, is then $ T $
similar to a contraction?
\par
\begin{itemize}
\item {Foguel \cite{Foguel_64} gave a counterexample in 1964.}
\item {Lebow \cite{Lebow_68} proved in 1968 that the operator, constructed by 
Foguel, is not even polynomially bounded.}
\item {Bo\.{z}ejko \cite{Bozejko_87a} 1987 produces a whole class of examples of 
power-bounded not polynomially bounded operators.}
\item {In between, Peller \cite{Peller_82} 1982 studied the space of functions
which act on power-bounded Hilbert space operators.}
\item {Halmos' question finitely has been answered in the negative by Pisier
    \cite{Pisier_96} 
1996, only recently.}
\item {An extension of Pisier's example has been obtained by 
Davidson and Paulsen \cite{DaviPaul_97} 1997.}
\end{itemize}

\section[Complete Boundedness]{Complete Boundedness}
\label{sec:2-cpbd}
For 
$n \in {\mathbb{N}}$ let $M_n = M_n({\mathbb{C}})=B(l^2_n)$ 
denote the space of complex
$n \times n$ 
\index{matrices}%
\index{matrices!$n \times n$}%
\index{matrices!$M_n$}%
\index{matrices!$B(l^2_n)$}%
matrices normed as acting
on the complex
$n$ 
dimensional Hilbert space $l^2_n$.\\

For 
$T \in B(H)$ 
a requirement which is stronger than its polynomial
boundedness is that of its 
\index{boundedness!complete polynomial}%
complete polynomial boundedness.
\begin{defi}
\label{defi:2cpbd}
An operator
$T \in B(H)$ 
is called completely poly\-no\-mi\-ally boun\-ded, if
there exists $C > 0$, such that for all 
$n\in{\mathbb{N}}$ and all 
$n \times n$ 
\index{matrix!of polynomials}%
matrices
of polynomials
$\matrpol{}{n}$ 
\begin{eqnarray*}
  \norm[M_n \otimes B(H)]{\matrpol{T}{n}} \le C\ \sup \left\{\, \norm[M_n]{
 \matrpol{z}{n} } : \, |z| \le 1  \, \right\} \,.
\end{eqnarray*}
\end{defi}
Here a 
\index{matrix!of operators}%
matrix 
$a=\matr{a}{}{n}$ 
of operators 
$a_{ij} \in B(H)$ 
acts on
\begin{eqnarray*}
l^2_n(H) \ = \ l^2_n \oplus_2 H \! &  =  & \! \left\{v= (v_1,\ldots ,v_n) \, : \,
 v_i \in H , \, \norm{v } =  \left(\sum_{i=1}^{n} 
\norm{v_i }^{2}\right)^{\frac{1}{2}} \right\} 
\end{eqnarray*}
by
\[
 a(v)_i = \sum_{j=1}^{n}\ a_{ij}v_j.
\]
\par
Based on the Sz.~Nagy-Foias dilation theorem it is not hard to show that any
contraction 
$T \in B(H)$, 
and hence any operator 
$T' \in B(H)$ 
which
is similar to a contraction, is completely polynomially bounded. Building on
work of Wittstock \cite{Wittstock_81},  Arveson \cite{Arveson_69}, \cite{Arveson_72} 
and Haagerup \cite{Haagerup_83}  Paulsen \cite{Paulsen_84}
showed the converse: \\
\index{theorem!Paulsen}%
\begin{thm}[Paulsen] 
\label{thm:similar}
An operator
$T' \in B(H)$ 
is similar to a contraction if and only if 
$T'$
is completely polynomially bounded.\\
\end{thm}
One more Definition:
\begin{defi}
\label{defi:cpbd}
If 
$K$ 
is another Hilbert space, 
$A \subset B(K)$ 
a linear subspace and
$\varphi : A \rightarrow B(H)$ 
a linear map, then 
$\varphi$ 
is called 
\index{bounded!completely}%
completely bounded, if there exists 
$C > 0$ 
such that for all 
$n\in {\mathbb{N}}$
and all 
$n\times n$ 
matrices
\linebreak[4] 
$a=\matr{a}{}{n} \in M_n(A)$
\begin{eqnarray}
\label{eq:cbnorm} 
  \norm[B(l^2_n \otimes_2 H)]{\matr{a}{\varphi}{n}} 
& \, \le \, & C
  \norm[B(l^2_n \otimes_2 K)]{\matr{a}{}{n}}.
\end{eqnarray}
We denote $ \cbnorm{.} : \varphi \mapsto \cbnorm{\varphi} $ the associated norm, i.e.\
\begin{displaymath}
  \cbnorm{\varphi } = \inf \{\, C \, : \, \mbox{\rm (\ref{eq:cbnorm})
holds true for all } a \in M_n(A) \mbox{ and all }n \in \mathbb{N}\, \}.
\end{displaymath}
\end{defi}
\par
In our case, with 
$S= \{z \in {\mathbb{C}} : |z|=1 \}$ 
and 
$\sigma$
 the surface measure
on 
$S$, 
let $ K=L^2 (S,\sigma)$ and consider the algebra of all 
\index{algebra!polynomials}%
polynomials
\begin{displaymath}
A \, = \, \{\, p \, :  p-\mbox{\rm polynomial} \} \, \subset \,  B(K)
\end{displaymath}
as embedded by means of its action of pointwise multiplication on 
\linebreak[4]
$L^2(S,\sigma)$-functions:
\begin{eqnarray*}
  f & \mapsto & p \cdot f ,\\
\mbox{where  for all} & p \in A, & f\in L^2(S,\sigma) \\
  p \cdot f(z) & = & p(z)f(z), \quad z \in S 
\end{eqnarray*}
denotes the pointwise product of functions.
\par
Then 
$T \in B(H)$
is completely polynomially bounded, exactly if the corresponding homomorphism
\[
p \mapsto p(T), \quad p \in A
\]
from $A$ into $B(H)$
is completely bounded.
\section[Analytic Semigroups on Hilbert Space]{Analytic Semigroups on
  Hilbert Space}
\label{sec:AnaSem}
It is obvious that the above mentioned problems have analogues for one-parameter
semigroups instead of the discrete semigroup 
${\mathbb{Z}}_+$
only.
\par
A 
\index{semigroup!$C_0$}%
$C_0$-semigroup
$(T_t)_{t\ge 0}$ 
acting on $H$ has an 
\index{generator!infinitesimal}%
infinitesimal generator
\begin{displaymath}
  Ax = \lim_{t\searrow 0}\ \frac{T_t - 1}{t}\ x,\quad x \in D(A) := \{
\,  \tilde{x} \,:
\,  \mbox{\rm the limit }  \lim_{t\searrow 0}\ \frac{T_t - 1}{t}\ \tilde{x} \;
\mbox{\rm exists } \}
\end{displaymath}
which is a closed densely defined operator.
\par
We shall consider here only semigroups 
which admit  bounded analytic extensions to some nontrivial 
\index{cone}%
cone
\index{$\Gamma_{\theta}$}%
\index{cone!$\Gamma_{\theta}$}%
$ \Gamma_\theta = \{z \not= 0 : | {\mbox{\rm arg}}(z) | < \theta \} \subset
{\mathbb{C}}$. (Analytic refers to weak analyticity, i.e.\ a map
$\Phi : \Gamma \rightarrow B(H) $ is called analytic, may be sometimes
holomorphic, if for all $x,y \in H$ $ z\mapsto (\Phi(z)x,y)$ 
is analytic on $\Gamma$.)\\
\begin{defi}
\label{defi:analytic}
\index{semigroup!analytic}%
\index{semigroup!$C_0$}%
A $C_0$-semigroup 
$(T_t)_{t\ge 0}$ 
acting on $H$
admits a bounded analytic 
\index{extension!bounded analytic}%
\index{analytic!extension}%
extension to some cone
$\Gamma_\theta$
if there exists a map
\begin{eqnarray*}
  T : \Gamma_\theta & \rightarrow & B(H) \\
              z & \mapsto & T_z,
\end{eqnarray*}
analytic on 
$\Gamma_\theta$
and extending
$T : t \mapsto T_t, $
such that
\begin{eqnarray*}
  T_{z+z'} & =  & T_zT_{z'} \quad \quad \forall\ z,z' \in \Gamma_\theta \,, \\
  \lim_{z \rightarrow 0, \, z \in \Gamma_{\theta}} T_z\ x & = & x \qquad
  \forall x \in H 
\end{eqnarray*}
and
\begin{displaymath}
  \sup_{z \in \Gamma_\theta}\ \norm{T_z } < \infty \,.
\end{displaymath}
\end{defi}
Assume that $A$ is injective then,
for 
$z \not= 0, 0 \leq {\rm Re}\ z<1$, 
complex powers 
$(-A)^z$ of the negative of $A$ 
can be defined as closable densely defined operators,
fulfilling 
$(-A)^z \cdot (-A)^{z'} = (-A)^{z+z'}$.
\begin{defi}
\label{defi:anapow}
We say that $(-A)$
\index{operator!imaginary powers}%
admits bounded imaginary powers if there exist 
$c > 0, C > 0$
such that 
$(-A)^{is} \in B(H)$ 
and
\begin{displaymath}
  \norm{(-A)^{is}}  \le C\ e^{c|s|} \quad \forall\ s\in \mathbb{R}.
\end{displaymath}
\end{defi}
Then in fact, 
$s\mapsto (-A)^{is}$ is a 
$C_0$-group of operators on $H$.
\par
\index{Le~Merdy}%
Le~Merdy~\cite{LeMerdy_95} connected the similarity problem for analytic
semigroups to the functional calculus of the negative of the generator:
\begin{thm}[Le~Merdy]
\label{thm:lemerdy} 
\index{theorem!Le~Merdy}%
Let 
 $(T_t)_{t\ge 0}$ 
 be a $C_0$-semigroup which admits an analytic extension. Let $A$ 
 be its generator and assume that $A$ is injective.\\
 Then $(-A)$ admits bounded imaginary powers if and only if the semigroup is 
 similar to a contraction semigroup.
\end{thm}
The implication 
``$ \Leftarrow $'' is entirely due to Pr\"uss and Sohr \cite{PrueSohr_90}.
For the other implication Le~Merdy adapts arguments of McIntosh
\cite{McInto_86} and Yagi
\cite{Yagi_84} to establish the complete boundedness of 
the functional calculus 
for rational functions of degree less than 
$-1$, with poles in the resolvent
set of $-A$ and uses a theorem of Paulsen analogous to the one cited above.   
\chapter[Dilation Theorems]{Dilation Theorems}
\section[${\mathcal{L}}^p$-Spaces as BanachLattices]{${\mathcal{L}}^p$-Spaces as Banach Lattices}
\label{sec:Lplat}
Since I shall treat an 
\index{$\mathcal{L}^{p}$-space}%
$\mathcal{L}^p$-space
$E \, = \, \Lp{\Omega}{\mu}$ 
as a Banach 
\index{lattice}%
\index{Banach lattice}%
lattice some remarks on
scalars are in order. Decomposing elements of $E$, that is, complex valued
functions, in their real and imaginary parts we obtain a decomposition, as
real linear spaces:
\index{$\mathcal{L}^{p}$-space!$L^p(\Omega ,\mu ;{\mathbb{C}})$}%
\begin{displaymath}
  L^p(\Omega ,\mu ;{\mathbb{C}}) = L^p(\Omega ,\mu ;\mathbb{R} ) \oplus i L^p(\Omega ,\mu ;\mathbb{R} )\,.
\end{displaymath}
Here, on 
$L^p(\Omega ,\mu ;\mathbb{R} )$, in addition to the linear operations the
pointwise operations of
\index{maximum}%
maximum and
\index{minimum}%
minimum of two elements and the
\index{value!absolute}%
absolute value are defined as
\begin{eqnarray*}
  f \vee g(\omega ) & = & \max\{ f(\omega ) ,g(\omega ) \},\quad \omega \in \Omega  ,\\
  f \wedge g (\omega ) & = & \min \{ f(\omega ), g(\omega ) \}, \quad
     \omega \in \Omega ,\ f,g \in L^p(\Omega ,\mu ; {\mathbb{R}}) ,\\
  \abs[f](\omega ) & = & f \vee (-f)(\omega )\ , \quad
     \omega \in \Omega , \quad f \in L^p(\Omega , \mu; {\mathbb{R}}).
\end{eqnarray*}%
\index{$f\vee g$}%
\index{$f\wedge g$}%
\index{$\abs[f]$}
Moreover, the norm on $L^p(\Omega ,\mu ;{\mathbb{C}})$ is related to the norm on 
$L^p(\Omega , \mu; {\mathbb{R}})$
by
\begin{displaymath}
  \norm{f} = \left( \int_{\Omega}\ \abs[f(\omega )]^p d\mu (\omega )\right)^{\frac{1}{p}} = 
  \left( \int_{\Omega} \left( \abs[{\rm Re}\ f]^2(\omega ) + \abs[{\rm Im}\ f]^2(\omega )\right)^{\frac{p}{2}}
  d\mu (\omega ) \right)^{\frac{1}{p}},
\end{displaymath}
where 
\index{${\rm Re}\ f$}%
\index{${\rm Im}\ f$}%
${\rm Re}\ f(\omega ) = \frac{1}{2}\left( f(\omega ) + {\overline{f (\omega )}}
\right),\; {\rm Im}\ f(\omega ) = \frac{1}{2} \left( f(\omega ) + {\overline{f(\omega )}}
\right) $ are elements of $L^p(\Omega ,\mu ;\mathbb{R} )$.
\par
One of the first facts to note is:
\begin{prop}
\label{prop:complexific}
If $T:L^p(\Omega ,\mu ;\mathbb{R} ) \rightarrow L^p(\Omega ,\mu ;\mathbb{R} )$ is a (real)
linear boun\-ded operation, then its 
\index{extension!complex}%
\index{extension!canonical}%
canonical extension
\index{$T_{\mathbb{C}}$}%
\begin{displaymath}
  T_{\mathbb{C}} (f + ig) = Tf + i\ Tg \qquad  f,g \in L^p(\Omega ,\mu ;\mathbb{R} )
\end{displaymath}
has the same norm bound as 
$T : \norm{ T_{\mathbb{C}} } = \norm{ T }$.
\end{prop}
{\bf Proof:}
I don't know, but for me this assertion does not appear completely trivial and thus
requires a proof.
\par
For
$\alpha,\beta \in \mathbb{R} $ 
and two independent Gaussian distributed 
\index{random variables!Gaussian}%
random variables $ X,Y $, with mean zero and variance 
one,
the variables  $Z = \alpha X + \beta Y$ 
and 
\linebreak[4]
$\tilde{Z}=(\alpha^2 + \beta^2)^{\frac{1}{2}} \cdot X$ 
are equidistributed.\
To prove this one computes the 
\index{Fourier~transform}%
Fourier~transforms
\begin{eqnarray*}
 {\mathcal E} (\exp it(Z)) & = & \int_{{\mathbb{R}}}\int_{{\mathbb{R}}}\ e^{-it(\alpha x + \beta y)}
 \frac{1}{\sqrt{2\pi}}\ e^{-\frac{1}{2}\abs[x]^2}\ \frac{1}{\sqrt{2\pi}}\ 
  e^{-\frac{1}{2}\abs[y]^2} \, dxdy\\
& = & \int_{\mathbb{R}}\ e^{-it\alpha x}\frac{1}{\sqrt{2\pi}}\
 e^{-\frac{1}{2}\abs[x]^2} \, dx\ 
  \int_{\mathbb{R}}\ e^{-it\beta y}\ \frac{1}{\sqrt{2\pi}}\ e^{-\frac{1}{2}\abs[y]^2}dy\\
& = & e^{-\frac{1}{2}(t\alpha)^2}e^{-\frac{1}{2}(t \beta)^2} = e^{-\frac{1}{2}[t^2(\alpha^2 + \beta^2)]}\\
& = & {\mathcal E} (\exp (it(\alpha^2 + \beta^2)^{\frac{1}{2}}X))\,.
\end{eqnarray*}
Thus for 
$f, \, g \in L^p(\Omega,\mu;{\mathbb{R}})$:
\begin{displaymath}
  {\mathcal E}(\abs[f(\omega)X + g(\omega)Y]^p) = {\mathcal E}\left(\abs[(f^2(\omega) +
  g^2 (\omega))^{\frac{1}{2}} \cdot X]^p \right),
\end{displaymath}
and with 
$\lambda \, = \, {\mathcal E}\abs[X]^p$:
\begin{displaymath}
  \int_{\Omega}(\abs[f]^2(\omega) + \abs[g]^2(\omega))^{\frac{p}{2}}d\mu (\omega) 
  = \frac{1}{\lambda}\ \int_{\Omega}\ {\mathcal E} \abs[f(\omega)X + g(\omega)Y]^p d\mu (\omega)\,.
\end{displaymath}
Let 
$h = f + ig$, 
then
\begin{eqnarray*}
  \norm{T_{\mathbb{C}} h }^p & = & \frac{1}{\lambda}\ 
\int_{\Omega}\ {\mathcal E} \abs[Tf(\omega)X +    Tg(\omega)Y]^p \, d\mu(\omega) \\
  & = & \frac{1}{\lambda}\ 
{\mathcal E}\ \int_{\Omega}\ \abs[T(f(\omega)X + g(\omega)Y)]^p\, d \mu(\omega) \\
 & \le & \norm{T   }^p\ \frac{1}{\lambda}\ 
{\mathcal E}\ \int_{\Omega}\ \abs[f(\omega)X +g(\omega)Y]^p \, d\mu (\omega)\\
  & = & \norm{T}^p \norm{h}^p\,.
\end{eqnarray*}
From this the inequality 
$ \norm{T_{\mathbb{C}}} \le \norm{T}$ 
is obvious.
\qed
\begin{rem} 
\label{rem:MarZyg}
\rm
We used here that a two-dimensional Hilbert space is isometrically isomorphic 
to a closed subspace of an 
$ \mathcal{L}^p$-space. This gave the equality of the norms of 
$T \mbox{ and its extension } T_{\mathbb{C}}$.
Our proof is a special case of the argument which Marcinkiewicz and
Zygmund gave for their famous extension theorem~\cite{MarcZygm_39}.
There, among others, the situation is considered of extending an operator
$ T \, : \, \Lp{\Omega}{\mu} \rightarrow \Lp{\Omega}{\mu}$
to an operator $T_H$ defined on the Hilbert space valued 
function space
$\Lp{\Omega}{\mu;H}$
by extending linearly the definition given on simple tensors
\begin{eqnarray*}
 \label{eq:MarcZygm}   
  T_H \, : \, f \otimes \xi & \mapsto & Tf \otimes \xi, 
\end{eqnarray*} 
where for $f \in \Lp{\Omega}{\mu}, \, \xi \in H$
the 
\index{tensor}%
\index{tensor!simple}%
simple tensor is the
\index{function!$H$-valued}%
$H$-valued function
$f \otimes \xi \ (\omega ) \ = \ f(\omega ) \xi, \quad \omega \in \Omega$.
\par
For a general  Banach lattice $ E $ and a  Hilbert space $H$
there is sense in an extension of the lattice and of an operator $T \in B(E)$. 
To prove the boundedness of the extended operator
one has to invoke Grothendieck's theorem to the result that
\[ 
\norm{T_H} \, \leq \, \mbox{K}_{\mbox{G}}\norm{T}, 
\] 
where 
$\mbox{K}_{\mbox{G}}$
denotes the  
\index{Grothendieck}%
\index{constant!Grothendieck}%
Grothendieck constant. The reader may find a discussion of
these facts in the book~\cite{LindTzaf_79}  of Lindenstrauss and Tzafriri.
\end{rem}

\begin{defi}
\label{defi:posop}
Now a 
linear\footnotemark %
\addtocounter{footnote}{-1}%
\footnotetext[1]{Whenever %
\index{scalar multiplication}%
\index{multiplication!scalar}%
scalar multiplication for complex numbers is defined, then linear actually means complex linear.}%
\addtocounter{footnote}{1}%
operator
\[   T : L^p (\Omega,\mu;\mathbb{C} ) \rightarrow L^p (\Omega,\mu;\mathbb{C} )
\]
is called 
\index{positive}%
\index{operator!positive}%
positive,
$T \ge 0$, 
if for all 
$f \in L^p (\Omega,\mu;\mathbb{C} )$
\[
f \ge 0\mbox{ implies } Tf \ge 0.
\]
\end{defi}
\begin{rem}
\label{rem:posop}
\rm
Clearly, a positive linear operator 
$T$ 
on a complex 
$\mathcal{L}^{p}$-space
leaves the real subspace
$L^p (\Omega,\mu;\mathbb{R})$ 
invariant, and it is the 
complexification of its restriction 
$T_{\mathbb{R}}$~\index{$T_{\mathbb{R}}$}.
\end{rem}
\begin{exam}
\rm
If 
$E \, = \, \Lp{\Omega}{\mu}$ 
is 
\index{finite dimensional}%
\index{$\mathcal{L}^{p}$-space!finite dimensional}%
finite dimensional, then, for some
$n$
and some 
$\omega_1 , \ldots , \omega_n > 0$,
\[
 E=l_n^p(\omega) = \left\{ \alpha = (\alpha_1 , \ldots , \alpha_n) \quad : 
\quad \alpha_i \in \mathbb{C} \quad \left( \sum^n_{i=1}\,\abs[\alpha_i ]^p \omega_i\right)^\frac{1}{p} =
  \norm{ \alpha } \right\}\,.
\]
For 
$1\, \leq \, i \, \leq \, n$
let
\index{$\delta_i$}%
$\delta_i = (0,0,1,0{\ldots} 0)$,
where the nonzero entry is at the i-th position, and let
$ \left\{ \delta_1 , \ldots , \delta_n \right\} $ 
be the 
\index{basis}%
\index{basis!standard}%
standard basis in
$\mathbb{C} ^n$ which we usually take as a basis of $l^{p}_{n}$.
\par
When a linear operator 
$T$ 
on 
$E$ is represented by its 
\index{matrix}%
matrix with respect to this basis, i.e.\
\[
(T\alpha)_i \, = \,  \sum_{j=1}^n \,T_{ij}\alpha_j \qquad i=1,\ldots , n\,,  
\]
then
$T$ is positive
if and only if 
$T_{ij} \ge 0 \mbox{ for all } i,j = 1, \ldots ,n$.\\
\begin{rem} 
\label{rem:poskern}
\rm
If 
$\mu$ is 
$\sigma$-finite and 
$T$ 
is given by a measurable kernel 
$k:\Omega \times \Omega \rightarrow \mathbb{C} $, i.e.\ if for $\mu$ almost all $ \omega \in
  \Omega\ $
\[
  Tf(\omega) = \int_{\Omega}\,k(\omega,\omega ')f(\omega ')d\mu(\omega '),
\] 
then 
$T \ge 0$ 
if 
$k \ge  0\quad \mu \times \mu$ 
almost everywhere.
\end{rem}
Examples of this type:
\begin{quote}
  Convolution on 
$L^p(\mathbb{R},\lambda)$ 
with Gauss, Poisson or other kernels.
\end{quote}
\end{exam}
For later use we note:

\begin{prop}
\label{prop:posmap}
If 
$ \,
T \, : \, L^p (\Omega,\mu;\mathbb{C} ) \rightarrow L^p (\Omega,\mu;\mathbb{C} ) $
is 
\index{operator!positive}%
positive, then
\[
  \norm{T} = \sup \left\{\, \norm{Tf } \, : \, f \ge 0,\,\norm{f }=1 \, \right\}\,.
\]
\end{prop}
{\bf Proof:}
It is obvious that the above
right hand side is dominated by
$\norm{ T } $. \\ 
Since 
$T =(T_{\mathbb{R}} )_{\mathbb{C}} $, 
we only need to prove
\[
  \norm{T_{\mathbb{R}}} \le  \sup \left\{\, \norm{Tf } \, : \, f \ge 0,\,\norm{f
  }=1 \, \right\}\ =: \, \lambda.
\]
But if 
$f \in L^p (\Omega,\mu;\mathbb{R}),\,$
then, denoting 
\index{$f^+$}%
$f^+=f{\vee}0$,
\index{$f^-$}%
$f^-= (-f)\vee 0 $,
\begin{eqnarray*}
  Tf & \, = \, & Tf^+ - Tf^- 
\end{eqnarray*}
and
\begin{eqnarray*}
 \abs[Tf] & \le & \abs[Tf^+] + \abs[Tf^- ]\\
    & = & Tf^+ + Tf^- \\
    & = & T(f^+ + f^- ).
\end{eqnarray*}
We obtain
\begin{eqnarray*}
  \norm{Tf } & \le & \norm{T(f^+ + f^- ) }\\
     &\le & \lambda \norm{f^+ + f^- } = \lambda \norm{f },
\end{eqnarray*}
since 
$\mu(\mbox{\rm supp}\,f^+ \cap \mbox{\rm supp}\,f^- ) = 0$.
\qed
\section[The finite dimensional Dilation Theorem of Ak\c{c}oglu and
Sucheston]{The Dilation Theorem of Ak\c{c}oglu and \\ Sucheston and its
Proof in the finite dimensional Case}
\label{sec:fdAkSu}
The results and proofs presented in this chapter are due to 
\index{Ak\c{c}oglu}%
Ak\c{c}oglu and
\index{Sucheston}%
Sucheston. We adapted them from the publications \cite{Akcoglu_75} and \cite{AkcoSuch_77}.
\begin{thm}[Ak\c{c}oglu and Sucheston\cite{AkcoSuch_77}]
\label{thm:AkcoSuchdilthm}
Assume
\index{dilation}%
\index{dilation!finite dimensional}%
$1 \le p < \infty$ 
and let $E = L^p(\Omega,\mu)$ 
be an $\mathcal{L}^p$-space,
$T:E \rightarrow E$ 
a positive contraction.
Then there exist another 
$\mathcal{L}^p$-space 
$\tilde{E} = L^p(\Omega ',\mu ')$
together with a positive invertible isometry 
$ S : \tilde{E} \rightarrow \tilde{E}$,
such that
\[ 
    DT^k = PS^k D\quad \mbox{\rm for}\quad k = 0,1,2,\ldots\, \; , 
\]
for some positive isometric embedding
$D : E \rightarrow \tilde{E}$ and a norm non-in\-crea\-sing positive projection
$P : \tilde{E} \rightarrow \tilde{E}$.
\end{thm}
We need prove this only for real scalars, and
we first assume that 
$E$ is finite dimensional. Then
for some 
${\omega}_1 > 0,\ldots ,{\omega}_n > 0$
\begin{eqnarray*} 
  E = l^p_n(\omega) = \left\{ \,(\alpha_1 , \ldots , \alpha_n ) \, : \, 
\alpha_i \in {\mathbb{R}} \right\} &
\mbox{ and } & 
  \norm[E]{ \alpha } = \left( \sum_{i=1}^{n}\, 
\abs[\alpha_i]^p {\omega}_i \right)^{\frac{1}{p}}.
\end{eqnarray*}
Moreover,
$E$ is isometric to 
$l^p_n$, 
by multiplication with 
${\omega}_i^{\frac{1}{p}}$:
\begin{eqnarray*}
 m_\omega \, : \, (\alpha_1 , \ldots , \alpha_n) & \rightarrow & 
(\alpha_1 \omega_1^{\frac{1}{p}} , \ldots , \alpha_n \omega_n^{\frac{1}{p}})\\
m_\omega \, : \,l_n^p(\omega) & \rightarrow &  l^p_n.
\end{eqnarray*}
Thus we may assume that
$E=l^{p}_{n}$,
for otherwise it would be sufficient to argue for 
$T' =m_\omega\, Tm_\omega^{-1}$,  $ \quad T' : l^p_n \rightarrow l^p_n\ $.
\par
Let $q$ be such that 
$\frac{1}{p} + \frac{1}{q} = 1$,  
and let
$T^\ast : l^{q}_n \rightarrow l^q_n$ 
denote the 
\index{adjoint}%
adjoint to
$T$ 
so that for 
$\alpha \in l_{n}^p = E , \,\beta \in l_n^q = E^\ast$ 
\[
 (T \alpha,\beta) =
(\alpha, T^\ast \beta), \]
where the 
\index{pairing}%
\index{pairing!bilinear}%
bilinear pairing
$ (.,.) \, : \, E \times E^{\ast} \rightarrow {\mathbb{C}}$
is given by
\[ 
(\alpha \ , \ \beta ) \, = \, \sum_{i=1}^{n} \ \alpha_i \beta_i.
\]
For 
$1<p<\infty$
we define a mapping 
\index{$\alpha^{\ast}$}%
\begin{eqnarray*} 
  ^{\ast} & \, : \,&  E^{+} \rightarrow E^{\ast \ +}\\
(\alpha_1 , \ldots , \alpha_n ) & \, \mapsto & (\alpha_1^{p-1}, \ldots , \alpha_n^{p-1}),
\end{eqnarray*}
such that  
\begin{eqnarray*}
 \norm[{q}]{ \alpha^\ast }^q & = & \sum_{i=1}^{n}\ \abs[\alpha^{p-1}_{i}]^q = 
                      \sum_{i=1}^n \ \alpha_{i}^{pq-q} =  \\ 
 & =   & \sum_{i=1}^n \ \alpha_{i}^p = \norm[p]{ \alpha }^p\,.\end{eqnarray*}
(If 
$p=1$,
then let
$ \alpha^\ast_{i} = 1 \mbox{ \rm if }  \alpha _{i} \not= 0, \,
\alpha^\ast_{i} = 0 \mbox{ \rm if } \alpha _{i} = 0$.) \\

We further define 
\[\begin{array}{rcl}
\index{M}%
M : E^+ &\rightarrow & (E^\ast)^+ \\
\mbox{by }\qquad M\alpha & = & T^\ast(T \alpha)^\ast. 
\end{array}\]
\begin{lem}
\label{lem:Mextr}
Assume
$1<p< \infty$.
For 
$ \alpha \in E^+$ 
with 
$\norm[p]{T\alpha } = \norm{ T} \norm[p]{\alpha} $ 
there holds true:
\[
M\alpha = \norm{ T }^p \alpha^{p-1}.
\] 
\end{lem}
{\bf Proof:}
We have
\[ 
   \norm[q]{M\alpha } \le \norm{T^\ast } \norm[q]{(T\alpha)^\ast} = \norm{T}
\norm[p]{T\alpha}^{p-1} =
   \norm{T}^p \norm[p]{\alpha}^{p-1},
\]
and if $q$ is defined by $ \frac{1}{p} + \frac{1}{q}=1 \ :$
\[
 \norm[p]{ \alpha } \norm[q]{ M \alpha} \geq (\alpha, M \alpha) = (T \alpha, (T \alpha)^{p-1}) = \norm{ T \alpha}^p = \norm{ T }^p
  \norm[p]{ \alpha }^p.
\]
Hence there is equality in H\"{o}lders 
\index{inequality!H\"{o}lder}%
inequality. It follows that for some
$\lambda \ge 0$
$M \alpha = \lambda\ \alpha^{p-1}$.
Then
$  (\alpha, \lambda \alpha^{p-1}) = \norm{T}^p \norm[p]{\alpha}^p
 \, \mbox{\rm shows } \, 
    \lambda = \norm{ T }^p\,. $ 
\qed
\index{$E(X)$}%
Let $E(X) = \left\{ \alpha \in E : \mbox{\rm supp } \alpha 
:= \{i \in \{1,\ldots ,n \}: 
\alpha_i \not= 0 \} \subset X \right \}$ 
and define
\index{$T_X$}%
$T_X\ \alpha =
   T(\chi_X \cdot \alpha)$.
We denote
$\lambda_X = \sup\{\norm{T\alpha}: \alpha \in E(X)^+, \, \norm{ \alpha }=1 \}$ 
the norm of $T_X$.
\begin{cor}
\label{cor:Mset}
Assume
$1<p< \infty$. 
 If $\alpha \in E(X)^+$ 
with 
$\norm{T\alpha} = \lambda_X \norm{\alpha} $, 
then
\[
  \chi_X \cdot M\alpha = \lambda^{p}_{X} \alpha^{p-1} \,.
\]
\end{cor}
{\bf Proof:}
This follows immediately from the above lemma, applied to the operator 
$T_X$
acting on the space
$E(X)$:\\
In fact, from
\[
   \chi_X \cdot T^\ast = (T_X)^\ast
\]
we infer
\begin{eqnarray*}
\chi_X \cdot M \alpha 
& \, = \, & \chi_X \cdot T^{\ast} \left( (T \alpha)^{p-1}\right) 
\, = \, (T_X)^{\ast} \left( (T \alpha)^{p-1}\right) \\
& \, = \, & (T_X)^{\ast} \left( (T_X \alpha)^{p-1}\right) 
\, = \, \lambda^{p}_{X} \alpha^{p-1} \,.
\end{eqnarray*}
\qed
\begin{lem}
\label{lem:Mextradd}
Assume
$1<p< \infty$.
Let 
$\alpha,\beta \in E^+ $ fulfil $ \alpha \cdot \beta = 0$ 
and 
$M\alpha \le \alpha^{p-1}$.
Then 
$\alpha M \beta = 0$ 
and 
$M(\alpha + \beta) = M\alpha + M \beta$.
\end{lem}
{\bf Proof:} 
We have 
$  0 \le (T\beta,(T\alpha)^{p-1}) = (\beta, M\alpha) \le 
(\beta,\alpha^{p-1}) = 0$, 
since 
$\alpha \cdot \beta = 0$. 
Since $T$ is a positive operator, 
$T \beta $ 
and 
$T \alpha $ 
must have disjoint supports.
Hence \qquad 
$0 = \left(T\alpha, (T\beta)^{p-1}\right)= (\alpha, M\beta)$,
from which we infer 
$ \alpha \cdot M\beta = 0$.
\par
Again the disjointness of the supports of
$ T\alpha \mbox{ and } T\beta$
implies  
$ (T(\alpha + \beta))^{p-1} \! = (T\alpha + T\beta)^{p-1} = (T\alpha)^{p-1} +
  (T\beta)^{p-1}$.
Applying
$T^{\ast}$
to this equality finally shows 
$ M(\alpha + \beta) = M\alpha + M\beta$.
\qed
\begin{lem}
\label{lem:Mexterw}
Assume
$1<p< \infty$.
Let 
$\alpha \in E^+$ 
satisfy 
$M\alpha \le \alpha^{p-1}$,
and assume that some
coordinates of 
$\alpha$ vanish. 
Then there exists 
$\tilde{\alpha}$ 
with strictly larger support than 
$\alpha$, 
such that
\[    
    M \tilde{\alpha} \le \tilde{\alpha}^{p-1}\,.
\]
\end{lem}
{\bf Proof:} 
Let 
$\quad X = \{ i : \alpha _{i} = 0 \}$. 
Since 
$E$ 
is finite dimensional there exists, by compactness of 
$\{ r \in E (X) : \norm{ r } = 1 \}$,
some
$\beta \in E(X) \mbox{ with } \norm{ \beta } = 1$ 
such that
\[
 \norm{ T \beta } = \norm{ T_X \beta } = \norm{T_X } \norm{\beta }.
\]
It follows that
\[
   \chi_X M\beta \le \beta^{p-1}.
\]
Then, since 
$\alpha \cdot \beta = 0$,
we may apply the last lemma to the result that 
$  \alpha M \beta = 0 $, 
and hence
$\mbox{supp }M \beta \subset X$.\\
Now we obtain
\begin{eqnarray*}
  M(\alpha + \beta ) & = & M \alpha + M \beta = M \alpha + \chi_X M \beta \\
  & \le & \alpha^{p-1} + \beta^{p-1} = (\alpha + \beta )^{p-1},
\end{eqnarray*}
and thus
$\tilde{\alpha} = \alpha + \beta$ 
can be chosen.
\qed
\begin{thm}
\label{thm:Mextrfkt}
There exists 
$u \in E^+$ 
with strictly positive coordinates such that
\[
  Mu \le u^{p-1}\,.
\]
\end{thm}
{\bf Proof:}
If 
$ p = 1$ 
take 
$u \equiv 1\;\mbox{, i.e.\ } u_i =1\; \mbox{ for }i=1,\ldots , n$.\\
Otherwise this follows immediately by Lemma~\ref{lem:Mexterw}.
\qed
\begin{rem}
\label{rem:isom}
\rm
 Assume 
$ p > 1 \mbox{ and } \norm{ T } = 1$.~\\[-2em] 
\begin{description}
\item[(i)~]
For
$u \in E^+ \mbox{ with } \norm[p]{u} = 1$
the assertions
$Mu = u^{\ast}$
and
$\norm[p]{ Tu } = 1$
are equivalent. In fact, denoting
$ v = Tu$,
\begin{eqnarray*}
\norm[p]{ Tu } & = & (v,v^\ast) \, = \, (Tu,v^\ast) \\
  & = & (u,T^\ast v^\ast )  \, = \, (u,Mu).
\end{eqnarray*}
Hence 
$\norm[p]{ Tu } = 1$
implies, by the converse to H\"older's inequality, that
$Mu = u^{\ast} = u^{p-1}$.
The converse is evident.
\par
From Proposition~\ref{prop:posmap} we know that 
$u$ can be chosen in $E_+$. In the case that 
$\norm{T} = 1 $ the proof of Theorem~\ref{thm:Mextrfkt} is thus much easier.
\item[(ii)]
Assume further all entries of the matrix of $T$ to be strictly positive.
Then, if $u \in E^+$ and if
$\norm[p]{ Tu } = 1 = \norm[p]{u}$,
the vectors
\[
   v = Tu \quad \mbox{\rm and} \quad u^\ast = T^\ast v
\]
have strictly positive coordinates. This can be seen from the equalities
\[
  v_j = \sum_{i=1}^{n}\ T_{ji}u_i\ , \qquad 
u^{p-1}_i = \sum_{j=1}^{n}\ T_{ji}v_j^{p-1}\,.
\]
\end{description}
\end{rem}
For 
$E= l^p_n$ 
the measure space 
$\Omega '$ appearing in the statement of the dilation
theorem, will be a subset of 
${\mathbb{R}}^2$, $\mu '$ 
will be the restriction of the
$2$-dimensional Lebesgue measure 
to $\Omega '$ 
and 
$S$ will be constructed
from a point transformation 
$\tau : \Omega ' \rightarrow \Omega '$ 
taking
into account the Radon-Nikodym derivative of 
$\mu ' \circ \tau^{-1}$ with
respect to 
$\mu '$.
\par
For 
$i=1,\ldots , n$ 
let 
$I_i$ 
be pairwise disjoint intervals on the
$x$-axes of 
${\mathbb{R}}^2$ 
of length 
$l(I_i )=1$ 
each. Let 
$J_i,\ i=1, \ldots ,n$ 
be mutually disjoint intervals too, again each of length one  
but on the
$y$-axes.\\
Set
\[
  X_i = I_i \times J_i\ , \quad i=1, \ldots ,n\ ; 
\quad Z_0 = \bigcup^n_{i=1}\ X_i\,.
\]
For 
$ k \not= 0,\ k \in {\mathbb{Z}}$ 
let 
$Z_k$ 
be mutually disjoint rectangles, each disjoint
 from 
$Z_0$ 
too, of positive finite 2-dimensional Lebesgue measure.
\par
By Theorem~\ref{thm:Mextrfkt} there is 
$u = (u_1, \ldots , u_n) \in E^+,\, u_i > 0$ 
for all 
$i \in \{ 1, \ldots , n\}$, 
with
\[
   Mu \le u^{p-1}.
\]
Let
\[   v = (v_1 , \ldots , v_n ) = Tu
\]
be its image under 
$T$. 
(It well might happen that $v_j = 0\mbox{ for some }j$.)
 Let 
$I = \{ 1, \ldots , n\}$ 
and 
$J = \{ j \in I \, : \, v_j \not= 0 \} $.\\
Define 
$   P =   I \times J $ 
as an index set,
\begin{eqnarray*}
    \xi_{ij} & := & T_{ji}\ \frac{u_i}{v_j},  \qquad (i,j) \in P\ \\
    \eta_{ij}& := & T_{ji}\ \left( \frac{v_j}{u_i} \right)^{p-1},  \qquad (i,j)
    \in I \times I.
\end{eqnarray*}
Since
\begin{eqnarray*} 
     v_j \; = \; (Tu)_j &  = & \sum_{i=1}^{n} \ T_{ji}u_i \qquad \mbox{\rm holds true, we have}\\
     \sum_{i=1}^{n} \ \xi_{ij} & = & 1 \qquad \mbox{\rm for all}\; j \in J\,.
\end{eqnarray*}
Similarly, from
\begin{eqnarray*}
      T^\ast (Tu)^{p-1} & = & Mu \le u^{p-1}, \qquad \mbox{\rm we obtain}\\
      \sum_{i=j}^{n} \ \eta_{ij} & \le & 1 \qquad \qquad \mbox{\rm for all} \ i \in I\,.
\end{eqnarray*}
Divide each 
$I_j,\ j \in J$, 
in 
$n$ 
subintervals 
$I_{ij}$ 
with length 
$\xi_{ij}$, 
and for each 
$i \in I$ 
choose
$n$
subintervals 
$J_{ij}$ 
in 
$J_i$
 with length 
$\eta_{ij}$.
 It well might happen, that some of those intervals degenerate, e.g.\ if
$(i,j) \notin P $
then $ \eta_{ij} =0$.\\
Let for
$ (i,j) \in P $
\begin{eqnarray*}
       S_{ij} & = & I_{ij} \times J_j,  \\
       R_{ij} & = & I_i \times J_{ij} ,
\end{eqnarray*}
and define
\[ 
   S = \bigcup_{(i,j)\in P}\ S_{ij}\ , \quad R = \bigcup_{(i,j)\in P}\ R_{ij}\,.
\]
Now there are affine transformations
\begin{eqnarray*}
\tau_{ij} : R_{ij} &\rightarrow &S_{ij}, \\
\tau_{ij}(x,y) &=& (a_{ij}x + b_{ij}, c_{ij}y + d_{ij}) \qquad (x,y) \in R_{ij}
\end{eqnarray*}
for some 
$a_{ij} , b_{ij}, c_{ij} , d_{ij} \in {\mathbb{R}}$,
which are
surjective up to sets of Lebesgue measure zero.
\par
\begin{figure}
\includegraphics{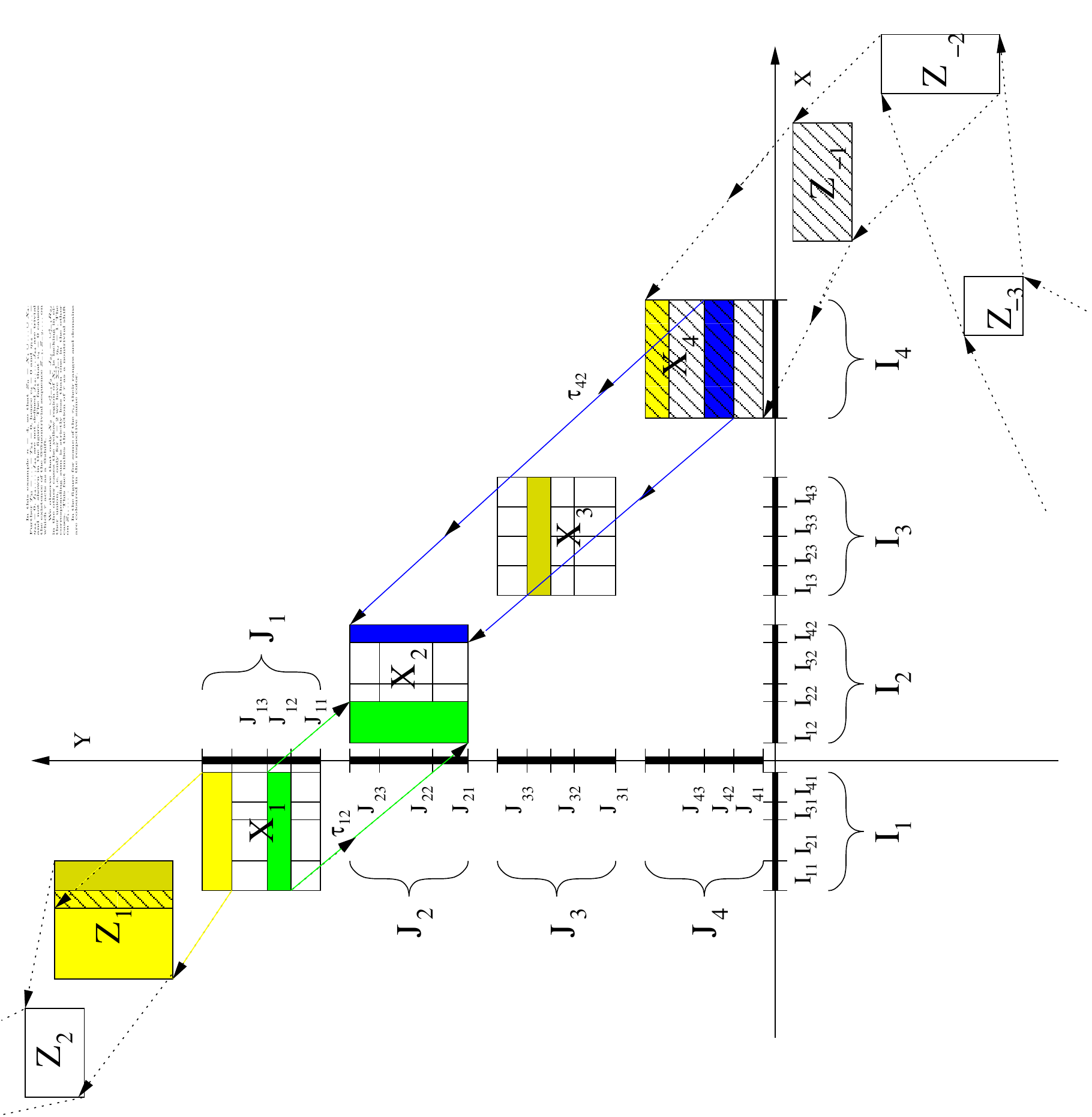}
\caption{The Transformation $\tau$.}
\label{fig:tau}
\end{figure}

We are going to define a point transformation $ \tau$ of 
$\bigcup^{-\infty}_{k=-1} Z_k \ \cup Z_0 \cup \ \bigcup^{+\infty}_{k=1}
Z_k$
onto itself:\\
\begin{enumerate}
\item
If 
$R = Z_0$ 
then define 
$\tau$ 
as the identity on 
$\bigcup^{+\infty}_{k=1}\ Z_k$
otherwise (piecewise affine) to transport 
$Z_0 \setminus R$ 
onto 
$Z_1, \quad Z_k$ 
onto
$Z_{k+1}\ \quad k \ge 1$.\\
\item
If 
$S = Z_0$ 
let 
$\tau$ 
be the identity on
 $\bigcup^{-\infty}_{k=-1}\ Z_k$
otherwise let 
$\tau$ 
map 
$Z_k$ 
onto 
$Z_{k+1}\; k \le -2$ 
and 
$Z_{-1}$ 
onto
$Z_0 \setminus S$.\\
\item
Now it only remains to define
\[
\tau \, : \, R \rightarrow S
\]
by
\[
    \tau_{| R_{ij}} = \tau_{ij} \,.
\]
\end{enumerate}
The Figure~\ref{fig:tau} 
indicates the properties of the point transformation 
$\tau$ in an example for the dimension $n = 4$.
\par
Let 
$ \Omega ' := \bigcup_{k \in {\mathbb{Z}}}\ Z_k $
and let
$\mu '$ 
be the restriction of the two dimensional Lebesgue measure 
to
$ \Omega '. $
Define 
$\nu = \mu '  \circ \tau^{-1}$ 
and denote
$\rho $
the Radon Nikodym derivative of $\nu $ with respect to $\mu ' $. \\
For $(i,j)\in P, \quad (x,y) \in S_{ij}$
\begin{eqnarray*}
  \rho(x,y) & = &  \frac{d\nu}{d\mu ' }(x,y) \; = \;\frac{\mu '  (\tau^{-1}(S_{ij})}{\mu ' (S_{ij})} \; = \; 
     \frac{\mu ' (R_{ij})}{\mu ' (S_{ij})} \; = \; \frac{\eta_{ij}}{\xi_{ij}}\\
  & = & \frac{T_{ji}(\frac{v_j}{u_j})^{p-1}}{T_{ji}(\frac{u_i}{v_j})} \; = \; 
     \left( \frac{v_j}{u_i} \right)^p (=: \rho_{ij})\,.
\end{eqnarray*}
\begin{enumerate}
\item
Then 
$S : f \mapsto Sf$
defined by
\[ 
    Sf(x,y):= \rho(x,y)^{\frac{1}{p}}\ f(\tau^{-1}(x,y)), \qquad (x,y) \in
    \Omega ',
\]
is an invertible isometry of 
$\Lp{\Omega '}{ \mu '  } = \tilde{E}$. \\
\item
$ P : \tilde{E}
\rightarrow \tilde{E}$ 
is just averaging over the sets 
$X_i$:
\[
  Pf(x,y) = 
\left\{
\begin{array}{ll}
        {\displaystyle  \frac{1}{\mu ' (X_i)}\ \int_{X_i} f d\mu '  }& \mbox{
        whenever } (x,y) \in X_i \mbox{ for some } i \in I\\
         0 & \mbox{ if } (x,y) \in Z_k \mbox{ for } k \not= 0.       
\end{array} 
\right.
\]

\item
$D : E \rightarrow \tilde{E}$
is given by
\[
    Df (\alpha_1 , \ldots , \alpha_n) = \sum_{i \in I}\ \alpha_i \chi_{X_i} (x,y)\ , 
\qquad (x,y) \in \Omega ' \,.
\]
\end{enumerate}
\begin{lem}
\label{lem:x-coord}
If 
$f,g \in E'$ 
are two functions supported in 
$ Z_0 \cup \bigcup_{k \geq 1} Z_k $ 
which depend only on the 
$x$-coordinate,
then 
$Pf = Pg $ implies $ PSf = PSg$.
Furthermore, $ PSf\mbox{ and } PSg$ have their supports in
$ Z_0 \cup \bigcup_{k \geq 1} Z_k $ again.
\end{lem}
{\bf Proof:}
If
$ f \, : \, \Omega ' \rightarrow {\mathbb{R}} $ depends only on the 
$x$-coordinate, then, for some function
$F \, : \, {\mathbb{R}} \rightarrow {\mathbb{R}}$,
we have
$f(x,y) = F(x), \qquad (x,y) \in \Omega '$.\\
For $ j \in J$ we compute
\begin{eqnarray*}
 \int_{X_j}\ Sf(x,y)\ d\mu ' (x,y) & = & \sum_{i \in I}\ \int_{S_{ij}}\; \frac{v_j}{u_i}\ 
     f(\tau^{-1}_{ij}(x,y)) \ dxdy\\
 & = & \sum_{i \in I} \ \frac{v_j}{u_i}\ \frac{\mu ' (S_{ij})}{\mu ' (R_{ij})}\ 
\int_{R_{ij}}\ f(x_{ij}) \ dxdy\\
 & = & \sum_{i \in I} \left( \frac{v_j}{u_i}\right)^{1-p}\ \int_{R_{ij}}\ F(x)
     \ dxdy\\
 & = & \sum_{i \in I} \left( \frac{v_j}{u_i} \right)^{1-p} \eta_{ij}\ \int_{I_i}\ 
     F(x) \ dx\\
 & = & \sum_{i \in I}\ T_{ji}\ \int_{X_i}\ f(x,y) \ d\mu ' (x,y)\,.
\end{eqnarray*}
If $ j \in I \setminus J $ then 
$X_j \subset Z_0 \setminus S$ and $ \tau^{-1}(X_j) \subset Z_{-1} $.
It follows that
$ PSf = 0 \mbox{ on } X_j $
whenever $ f=0 \mbox{ on } Z_{-1}$.
Hence,
\begin{eqnarray}
 PSf & = & \sum_{j\in J}\ \sum_{i \in I}\ T_{ji}\ \int_{X_i}\ f(x,y) \ d\mu ' (x,y)\ \chi_{X_j}
\label{PS}\\ 
&  = & \sum_{j\in I}\  \sum_{i \in I}\ T_{ji}\ \int_{X_i}\ f(x,y) \ d\mu ' (x,y)\ \chi_{X_j}.
\end{eqnarray}
Here we could extend the summation with respect to
 $ j $ 
from $J$
to all of 
\linebreak[4]
$I= \{1 , \ldots , n \}$,
since for 
$j_0 \in I \setminus J$
we have $T_{j_0 i} = 0 \quad \forall i \in I$.
In fact, because of
$  u_i > 0 \quad \forall \ i \in I $,
it follows that
$ j_0 \in I \setminus J , \mbox{ i.e.\ } 
0=v_{j_0} = \sum_{i \in I} \ T_{j_0i} u_i$,
implies $ T_{j_0 i} = 0 \quad \forall \ i \in I $.

If $Pf = Pg$ and $f,g$ are both supported in 
$ Z_0 \cup \bigcup_{k \geq 1} Z_k $, 
then for all $ i \in I$ 
$\int_{X_i}\ f(x,y) \ d\mu ' (x,y) = \int_{X_i}\ g(x,y) \ d\mu ' (x,y)$,
and hence
\[ 
    PSf = PSg\,.
\]
\qed
Further, for 
$\alpha \in E$ 
we may restate (\ref{PS}) as: 
\begin{lem}
\label{lem:indvor}
$PS(\sum_{i \in I}\ \alpha_i \chi_{X_i}) = \sum_{j\in J}\ (T\alpha)_j \ \chi_{X_j} \ $.
\end{lem}
\begin{thm}
\label{thm:findilthm}
If 
$k \in {\mathbb{Z}}_+, \quad \alpha \in E$, 
then
\[
  PS^k\ \sum_{i \in I}\ \alpha_i \chi_{X_i} = \sum_{j \in J}\ (T^k \alpha)_j\chi_{X_j}\,.
\]
\end{thm}
{\bf Proof:}
We prove the theorem by induction on k.
The assertion is clear for 
$k = 0$
and just proved and stated in the last lemma  
for 
$k = 1$. If $f \in \tilde{E}$ 
depends
only on the 
$x$-coordinate, 
then the same is true for 
$Sf$, since $\tau$ 
 is 
piecewise affine and also the Radon Nikodym derivative depends
only on the 
$x$-coordinate.
Hence 
$S^k f$ 
depends only on the 
$x$-coordinate. Then from Lemma~\ref{lem:x-coord}
\begin{eqnarray*}
  PS^{k+1} \sum\alpha_i \chi_{X_i} & = & PS\ P\ S^k\ \sum\alpha_i \chi_{X_i}\\
    & = & PS \left( \sum_{i \in I}\ (T^k \alpha)_i \chi_{X_i} \right)\\
    & = & \sum_{j \in J}\ (T^{k+1} \alpha)_j \chi_{X_j}\,.
\end{eqnarray*} 
This proves the Ak\c{c}oglu-Sucheston theorem in the finite dimensional case.
\qed
\chapter[Ultraproducts of Banach Spaces]{Ultraproducts of Banach Spaces}
\section[The general Banach space Ultraproduct Construction]{The general Banach space Ultraproduct Construction}
\label{sec:genUlpr}
This construction might be based either on the notion of an 
\index{ultrafilter}%
ultrafilter or
equivalently (at least in the case we are interested in) on points of the
\index{Gelfand space}%
\index{space!Gelfand}%
Gelfand space of the algebras $l^{\infty}$.
\par
Let 
$(A,\le)$ 
be a partially ordered, directed set and let 
\index{$\overline{A}$}%
$\overline{A}$
denote the Gelfand space of 
\index{$l^{\infty} (A;{\mathbb{C}})$}%
$l^\infty (A;{\mathbb{C}})$,
i.e.\ the set of 
\index{multiplicative functional}%
\index{functional!multiplicative}%
multiplicative
functionals
$\alpha : l^\infty (A;{\mathbb{C}}) \rightarrow {\mathbb{C}}$ 
with the weak--$\ast$--topology. Then 
$A \hookrightarrow {\overline A}$ 
is canonically embedded
in the compact Hausdorff topological space 
\index{$\overline{A}$}%
$\overline{A}$
and is a dense subset 
in there.\\
\par
For 
$f \in l^\infty (A;{\mathbb{R}})$ 
define
\[
   \lim_{\alpha \in A} \inf f(\alpha) = \sup_{\beta \in A} \inf_{\alpha \ge \beta}
     f(\alpha)
\]
and
\[  \lim_{\alpha \in A} \sup f(\alpha) = \inf_{\beta \in A} \sup_{\alpha \ge \beta}
      f(\beta)\,.
\]
\begin{prop}
\label{prop:Blimex}
There exists a point
\index{$\mbox{LIM}$}%
\index{LIM}%
$\mbox{\rm LIM}\ \in \overline{A}$ 
such that
\[
   \lim_{\alpha \in A} \inf f(\alpha) \le \mbox{\rm LIM} (f) \le 
    \lim_{\alpha \in A} \sup f(\alpha)\quad \forall\, f \in l^\infty (A;\RR)\,.
\]
\end{prop}
{\bf Proof:}
For 
$\alpha \in A$ 
let 
$A_{\alpha} = \{ \beta : \beta \ge \alpha \}$ 
and let
$\overline{A}_\alpha$ 
be its closure in 
$\overline{A}$.
Since 
$A$ 
is a 
\index{directed}%
directed set, there exists, for
$\alpha, \beta \in A$, 
some 
$\gamma \in A$ 
such that
\[
   \emptyset \not= A_\gamma \subset A_\alpha \cap A_\beta\,.
\]
Hence the sets 
$(\overline{A}_\alpha)_{\alpha \in A}$ 
have the finite intersection
property, i.e.\ for any finitely many 
$\alpha_1 , \ldots , \alpha_n:$
$\bigcap_{i=1}^n \ \overline{A_{\alpha_i}} \not= \emptyset$. 
In the compact space
$\overline{A}$ 
there then exists a point
\[
   \mbox{\rm LIM}\ \in \bigcap_{\alpha \in A}\ \overline{A_\alpha}\,.
\]
But this implies that for any 
$f \in l^\infty (A;{\mathbb{R}})$
\[
   \mbox{\rm LIM}\ (f) \in \bigcap_{\alpha \in A}\ f(\overline{A}_\alpha )
   \subseteq \bigcap_{\alpha \in A}\ \overline{f(A_\alpha)}\ ,
\] 
and the last set is contained in the interval
\[
  \left[ \, \lim_{\alpha \in A}\inf f(\alpha) \, , 
\,\lim_{\alpha \in A} \sup f(\alpha) \, \right]\ ,
\]
which proves the proposition.
\qed
\begin{rem} \rm \hspace*{\fill}~\\[-2em]
\label{rem:Blimit}
\begin{description}
\item[(i)~~] 
If 
$f \in l^\infty (A;{\mathbb{C}})$ 
and 
$u : \overline{f(A)}
  \rightarrow C$ 
is continuous, then 
\[
\quad \mbox{\rm LIM}\ u \circ f =
  u(\mbox{\rm LIM}\ f).
\]
In fact if 
$\hat{f}$ 
denotes the 
\index{Gelfand transform}%
\index{transform!Gelfand}%
Gelfand transform of
$f$, 
  then 
$u \circ f = \left(u \circ \hat{f}\right)|_A$ 
and 
\begin{eqnarray*}
\mbox{\rm LIM}\ u \circ f & = &
 (u \circ f \hat{)} (\mbox{\rm LIM}) \, = \, u \circ \hat{f}(\mbox{\rm
  LIM})\\
& = & u \left( \hat{f} (\mbox{\rm LIM}) \right) \, = \, u(\mbox{\rm LIM}\ f).
\end{eqnarray*}
\item[(ii)~] 
If 
$f,g \in l^\infty (A;{\mathbb{R}})$
 are such that 
$f \le g$, 
then
\[
    \mbox{\rm LIM}\ f \le \mbox{\rm LIM}\ g\,.
\]
Clearly 
$f \le g$ 
on 
$A$ 
implies 
$\hat{f} \le \hat{g}$ 
on 
$\overline{A}$. 
\item[(iii)]
\mbox{\rm LIM} is usually called a 
\index{limit!generalised}%
generalised limit, or Banach
\index{limit!Banach}%
\index{Banach limit}%
limit.
Its value at a function
$ f \in l^\infty (A;{\mathbb{C}})$
we sometimes denote
$\mbox{\rm LIM}_{\alpha \in A}\, f(\alpha)$.
\end{description}
\end{rem} \rm
Now let 
($A,\le $, 
LIM) be as above. For a net 
$\left(E_\alpha , \norm[\alpha]{.}\right)_{ \alpha \in A}$
of Banach spaces let 
\index{$\Lambda (A,E)$}%
${\Lambda}(A,E)$
denote the space of those functions 
$f : A \rightarrow \bigcup_{\alpha \in A}\ 
E_\alpha$ 
such that
\[
  f(\alpha) \in E_\alpha \quad \forall\, \alpha \in A \mbox{ and }
 \norm[\infty]{f} 
= \sup_\alpha
  \norm[\alpha]{ f (\alpha) } < \infty\,.
\]
Then clearly 
$\Lambda (A,E), \norm[\infty]{.}$ 
is a normed vector space, it
is even complete and a Banach lattice if all 
$E_\alpha$ 
are 
Banach lattices.
The operations are defined pointwise:
\[   \begin{array}{lcll}
     (f + g)(\alpha) & = & f(\alpha) + g(\alpha), & \alpha \in A\\
     f \wedge g (\alpha & = & f(\alpha) \wedge g(\alpha), & \alpha \in A\\
     \mbox{\rm etc.\ } & & &
     \end{array}
\]
Let us define a semi-norm 
\index{$\norm[{\mbox{\rm \BLIM}}]{.}$}%
$\norm[{\mbox{\rm \BLIM}}]{.}$
on 
$\Lambda (A,E)$ 
by
\[
  \norm[{\mbox{\rm \BLIM}}]{ f } = \mbox{\rm LIM}_{\alpha \in A} \norm[\alpha]{ f(\alpha)}  \,.
\]
We let 
$N = \{ f \in \Lambda (A,E) \, : \, \norm[{\mbox{\rm \BLIM}}]{ f } = 0 \}$ 
denote
its kernel and denote
\index{$\prod_{\mbox{\rm \BLIM}} E_\alpha$}
\[
   \prod_{\mbox{\rm \BLIM}} E_\alpha = \Lambda (A,E)/N \,. 
\]
Further, for 
$f,g \in \Lambda (A,E)$, 
we write
\index{$f \sim g$}%
$f \sim g$ 
if
\index{$N$}%
$f-g \in N$
 and 
\index{$[f]$}%
$[f] = f + N$.
\par
The following useful observation seems to be due to Ak\c{c}oglu and Sucheston~\cite{AkcoSuch_77}. 
\begin{prop}
\label{prop:complete}
$\Lambda (A,E)/N$ 
is a Banach space.
\end{prop}
This proposition will be proved by the following two lemmata.
\begin{lem}
\label{lem:equiclass}
For 
$f\in \Lambda (A,E)$ 
there exists 
$g \in \Lambda (A,E)$ 
such that 
\[
f \sim g  \mbox{ and }
\norm[\alpha]{g(\alpha)} \le \norm[{\mbox{\rm \BLIM}}]{ f } \; 
\forall \alpha \in A.
\]
\end{lem}
{\bf Proof:}
For 
$\alpha \in A $
denote 
$ \lambda_\alpha =\norm[{\mbox{\rm \BLIM}}]{ f} / \max\{ \norm[\alpha]{ f
  (\alpha )},
\norm[{\mbox{\rm \BLIM}}]{ f } \} $
and define 
$g \in \Lambda (A,E)$ 
by
\[
  g(\alpha) = \lambda_\alpha f(\alpha),\qquad \alpha \in A.
\]
Then,
\begin{eqnarray*}
  \norm[\alpha]{ f(\alpha) - g(\alpha) } & = & \norm[\alpha]{ (1 - \lambda_\alpha )f(\alpha) }\\
     & = & (1-\lambda_\alpha) \norm[\alpha]{ f (\alpha )} \,.
\end{eqnarray*}
Hence,
\begin{eqnarray*}
  \norm[ {\mbox{\rm \BLIM}}]{ f - g} & \le & \mbox{\rm LIM}_{\alpha \in A}\ (1-\lambda_\alpha)
  \cdot  \mbox{\rm LIM}_{\alpha \in A}\ \norm[\alpha]{ f(\alpha)}\\
      & = & 0 \cdot \norm[{\mbox{\rm \BLIM}}]{ f }\,.
\end{eqnarray*}
It is clear that for all $\alpha \in A $:
$\quad \norm[\alpha]{ g(\alpha) } = \lambda_\alpha \norm[\alpha]{ f(\alpha) } 
\le \norm[{\mbox{\rm \BLIM}}]{ f }$. 
\qed 
Now the completeness of 
$\Lambda (A,E)/N$ 
follows from:
\begin{lem}
\label{lem:compl}
Let 
$\left(f_{n}\right)_{n=1}^{\infty} \in \Lambda (A,E)$ 
be a sequence such that \\
$\sum_{n=1}^{\infty}\ \norm[ {\mbox{\rm \BLIM}}]{ f_{n} } < \infty$.
Then there exists
\[
  f \in \Lambda (A,E)\quad \mbox{\rm such that}\quad \lim_{N\rightarrow \infty}\ 
  \norm[{\mbox{\rm \BLIM}}]{ f - \sum_{n=1}^N\ f_{n} } = 0\,.
\]
\end{lem}
{\bf Proof:}
Let 
$g_{n} \sim f_{n}$ 
be such that 
$\norm[\alpha]{ g_{n} (\alpha) } \le \norm[{\mbox{\rm \BLIM}} ]{ f_{n} }\ $
for all 
$\alpha \in A$. 
Then 
$g(\alpha) = \sum_{n=1}^{\infty}\ g_{n} (\alpha) \in E_\alpha$
exists, since 
$ \sum_{n=1}^{\infty}\ \norm[\alpha]{g_{n} (\alpha) } \le \sum_{n=1}^{\infty}\ \norm[ {\mbox{\rm \BLIM}}]{ f_{n} }
< \infty$. 
This estimate further shows 
$g \in \Lambda (A,E)$. 
It is clear that for all $ N \in {\mathbb{N}}:$
$\quad \sum_{n=1}^N\ g_{n}\ \sim \sum_{n=1}^{N}\ f_{n}$.\\
Now,
\begin{eqnarray*}
  \norm[{\mbox{\rm \BLIM}} ]{ g - \sum_{n=1}^N\ f_{n} } & = & 
\norm[{\mbox{\rm \BLIM}}]{ g - \sum_{n=1}^N\ g_{n} } =\\
  & = & {\mbox{\rm LIM}_{\alpha \in A}} \norm[\alpha]{ \sum^\infty_{n=N+1} g_{n}(\alpha )} \le \
      {\mbox{\rm LIM}}_{\alpha \in A} \sum^\infty_{n=N+1} \norm[\alpha]{ g_{n}(\alpha) }\\
  & \le & \sum^\infty_{n=N+1} \norm[{\mbox{\rm \BLIM}}]{f_{n} } \rightarrow 0 \;
    \mbox{\rm if}\; N \rightarrow \infty\,.\; \,  \quad \qquad  \qquad \qquad \qedm
\end{eqnarray*}
\begin{rem} \rm \hspace*{\fill}~\\[-2em]
\label{rem:opulpr}
\begin{description}
\item[(i)~~]  
If 
$E =(E_\alpha )_{\alpha \in A},\, F= (F_\alpha )_{\alpha \in A}$ 
are as
above and 
$u_\alpha : E_\alpha \rightarrow F_\alpha$ 
is a bounded linear
operator such that 
$\quad \sup_{\alpha \in A}\norm[\alpha]{ u_\alpha } < \infty$, 
then there exists a linear operator
\[
   u : \Lambda (A , E) \rightarrow \Lambda (A,F)
\]
defined
by
\[ 
 u(f_\alpha ) = u_\alpha (f_\alpha), \quad \alpha \in A\,.
\]
Furthermore,
\[
  \norm[{\mbox{\rm \BLIM}}]{ u(f_.)} \le \sup_{\alpha \in A} \norm[\alpha]{u_\alpha }\; 
  \norm[{\mbox{\rm \BLIM}} ]{f_.}.
\]
Hence 
$u$ 
passes to the quotient spaces and defines a bounded linear 
\index{ultraproduct!maps}%
map, the
ultraproduct
of the net $\left(u_{\alpha}\right)_{\alpha \in A}$, denoted 
\[
\prod_{\BLIM} u_{\alpha} \,: \, 
\prod_{\mbox{\rm \BLIM}}\ E \rightarrow \prod_{\mbox{\rm \BLIM}}\ F,
\]
of norm at 
most 
$\sup_{\alpha \in A} \norm{u_\alpha }  $.
\item[(ii)~] 
If 
$E$ 
is a 
Banach space and 
$\left(E_\alpha\right)_{\alpha \in A} $ 
is a
net of its subspaces 
$(E_\alpha \subseteq E \quad \forall \alpha \in A )$,
ordered by inclusion
$(\mbox{i.e.\ }\alpha' \le \alpha \Leftrightarrow E_{\alpha '} \subseteq 
E_\alpha )$, 
such that 
$\bigcup_{\alpha \in A}\ E_\alpha$ 
is dense in 
$E$,
then 
$E$ 
is canonically embedded in 
$\prod_{\mbox{\rm \BLIM}}\ E_{\alpha}$ 
by just taking the ultraproduct of the 
\index{ultraproduct!inclusions}%
inclusions.\\
For, if 
$f \in E$, 
then there exists a net 
$f_\alpha \in E_\alpha$ 
such that
$\norm[E]{ f - f_\alpha } \rightarrow 0$.\\
Define 
$\tilde{f} \in \prod_{\mbox{\rm \BLIM}}\ E_{\alpha}$ 
as 
$\tilde{f} = [f_\alpha ]$,
then
\[
   \norm{ \tilde{f} } = {\mbox{\rm LIM}}_{\alpha \in A}\ \norm[{\alpha}] {f_\alpha} =
   {\mbox{\rm LIM}}_{\alpha \in A}\ \norm[E]{f_\alpha } = \lim_{\alpha \rightarrow \infty} 
   \norm{f_\alpha } = \norm{ f } \,.
\]
It is straightforward to check, that the definition of 
$\tilde{f}$ 
does not depend on the special net 
$(f_\alpha )_{\alpha \in A}$ 
and that 
$f \mapsto \tilde{f}$ 
is linear.
\item[(iii)] 
The (real) 
\index{$\Lp{\Omega}{\mu;{\mathbb{R}}}$}%
$L^p(\Omega,\mu;{\mathbb{R}})$
spaces, 
$1 \le p < \infty$, 
are 
characterized by the Kakutani and the 
\index{theorem!Kakutani}%
\index{theorem!Bohnenblust-Nakano}%
Bohnenblust-Nakano theorem (see e.g.\ 
\cite{Lacey_74} Chap.5 \S 15 Theorem 3) as those Banach lattices 
$E$ 
such that
\[
  f,g \in E^+, \; f \wedge g = 0 \Rightarrow \norm{ f + g }^p = \norm{f }^p +
  \norm{ g }^p \,.
\]
\end{description}
\end{rem}
As a corollary to this fact we have: 
\begin{cor}
\label{cor:Lpulpr}  
\index{$\prod_{\mbox{\rm \BLIM}} E_\alpha$}%
Let 
$(E_\alpha )_{\alpha \in A}$ 
be a net of 
$\mathcal{L}^p$-spaces, then there exists
a measure space 
$(\Omega^{\circ},\mathfrak{A}^{\circ} \mu^{\circ})$ 
such that
\[
   \prod_{\mbox{\rm \BLIM}}\ E_\alpha = \Lp {\Omega^{\circ}}{ \mu^{\circ}}\,.
\]
\end{cor}
{\bf Proof:}
Let 
$f,g \in \prod_{\mbox{\rm \BLIM}}\ E_\alpha^+$ 
with 
$f \wedge g = 0$ 
be
represented by 
$(f_\alpha)_{\alpha \in A}$ 
respectively $(g_\alpha )_{\alpha \in A}$.
Then 
$h = (f_\alpha \wedge g_\alpha )_{\alpha \in A} \in N$, 
since it is
a representation for 
$f \wedge g = 0$, 
and we may further assume that 
$f_\alpha, g_\alpha \in E^{+}_{\alpha} $
for all 
$ \alpha \in A$. 
\par
Now, for all $\alpha \in A$,
the functions
$f_\alpha - (f_\alpha \wedge g_\alpha)$ 
and 
$g_\alpha - (f_\alpha \wedge g_\alpha)$ 
have disjoint supports, hence, using 
$E_{\alpha} \in \mathcal{L}^p$,
\begin{eqnarray*}
 \lefteqn{  \norm[E_{\alpha}]{f_\alpha - (f_\alpha \wedge g_\alpha) + g_\alpha     - 
(f_\alpha \wedge g_\alpha )}^p } & & \\
& = &   \norm[E_\alpha] {f_\alpha -(f_\alpha \wedge g_\alpha)}^p + 
\norm[E_\alpha]{g_\alpha - (f_\alpha \wedge g_\alpha)}^p\,.
\end{eqnarray*}

Because
\[
    [(f_\alpha - g_\alpha \wedge f_\alpha)_\alpha] = f, \quad [(g_\alpha - 
    (f_\alpha \wedge g_\alpha))_{\alpha}] = g\ ,
\]
we have, by passing to the limit: 
\begin{eqnarray*}
  \norm{ f + g }^p & = & {\mbox{\rm LIM}}_{\alpha \in A} \norm[\alpha]{f_\alpha + g_\alpha - 2(f_\alpha 
     \wedge g_\alpha) }^p = \\
  & = & {\mbox{\rm LIM}}_{\alpha \in A} \norm{ f_\alpha -f_\alpha \wedge g_\alpha }^p + 
     {\mbox{\rm LIM}}_{\alpha \in A} \norm{ g_\alpha -(f_\alpha \wedge g_\alpha) }^p\\
  & = & \norm{ f }^p + \norm{g }^p\,.
\end{eqnarray*}   
\begin{rem} \rm
A more direct approach to the above Corollary~\ref{cor:Lpulpr},
using directly ultraproducts and not relying on the Bohnenblust-Nakano or on
Kakutani's theorem, can be found in the publication~\cite{DaCaKriv_72} of
Dacunha Castelle and Krivine.
\end{rem}
\section[The Ak\c{c}oglu-Sucheston Dilation %
Theorem]{\!\!The \!Ak\c{c}oglu-Sucheston \!Dilation \!Theorem}
\label{sec:AkSuProof}
In this section we shall use the foregoing constructions to complete the
proof of the Ak\c{c}oglu-Sucheston dilation theorem in the general case.
We shall rather closely follow their arguments.
\begin{defi}
\label{defi:semipart}
A 
\index{semi-partition}%
semi-partition
$\alpha$ of $\Omega$ 
is a finite collection of pairwise
disjoint subsets (measurable of course) 
\[
\alpha = \{ X_1 , \ldots , X_{n_\alpha} : X_i \cap X_j =
\emptyset \mbox{ if } i \not= j \}
\]
each one of finite measure.
\end{defi}
For a semi-partition 
$\alpha$ 
we let
$\mathcal{E}_\alpha : E \rightarrow E $ denote the corresponding conditional
expectation operator, defined by:
\[
   \mathcal{E}_\alpha(g)(\omega) =
  \left\{ \begin{array}{ll}
    {\displaystyle \frac{1}{\mu (X)}\, \int_X\ g(\omega')\ d\mu(\omega')} & \mbox{\rm if } \omega \in X \in \alpha\\
    0 & \mbox{\rm otherwise}. 
          \end{array} \right. 
\]
\begin{rem} \rm \hspace*{\fill}~\\[-2em]
\label{rem:semipart}
\begin{description}
\item[(i)~~] 
The set 
$A = \{ \alpha : \alpha \; \mbox{\rm  a semi-partition} \}$ 
of all
semi-partitions is 
\index{ordered!directed by refinement}%
\index{ordered!partially}%
partially ordered,
\index{directed}%
directed by
\index{refinement}%
refinement, and for
$g \in E$
we have
\begin{eqnarray*}
   \lim_{\alpha \in A}\ \norm[p]{ \mathcal{E}_\alpha\ g - g} \rightarrow 0.
\end{eqnarray*}
\item[(ii)~]
Further, if 
$T : E \rightarrow E$ 
is a positive contraction, and if
\[
   T_\alpha : E \rightarrow E \mbox{ is defined by } T_\alpha = \mathcal{E}_\alpha T \mathcal{E}_\alpha,
\]
then it may be restricted to 
$E_\alpha := \mathcal{E}_\alpha E$ 
and there be viewed as a
positive contraction on 
$E_\alpha \simeq l^p_{n_\alpha}(\omega) \simeq l^p_{n_\alpha}$.
(Here $\simeq $ is a positive isometrical isomorphism and $\omega$ is the
\index{weight}%
\index{sequence!weight}%
weight sequence $ (\mu(X_1), \ldots ,\mu(X_n)) $.)
\item[(iii)]
Since 
$A$ 
is directed by refinement 
and since 
$\norm{T}<\infty$
we obtain by induction from (i):
\[
  \lim_{\alpha \in A} \norm{ T^n_\alpha\ f - f } = 0\ , \quad n = 0,1, \ldots \quad
    \forall\, f \in E\,.
\]
(For a proof just note:
$\norm{ \mathcal{E}_\alpha T \mathcal{E}_\alpha f - Tf }  \le  \| { \mathcal{E}_\alpha
  T \mathcal{E}_\alpha f - \overbrace{T {\mathcal{E}_\alpha} f}^{=: g} } \| + \\
  \norm{T \mathcal{E}_\alpha f - Tf} \leq 
  \norm{ \mathcal{E}_\alpha g - g } + \norm{ T } \norm{ \mathcal{E}_\alpha f - f }
  \stackrel{\alpha \rightarrow \infty}{\longrightarrow} 0
$.)
\end{description}
\end{rem}
For a semi-partition $\alpha$ we display in Figure~\ref{fig:dil1} a 
\index{dilation}%
dilation of
$T_\alpha$ 
according to what already  has been proved (Theorem~\ref{thm:AkcoSuchdilthm}).
\begin{figure}[h]
\unitlength1cm
\begin{picture}(13,5)
\put(1.5,4){$E_\alpha $ }
\put(10.5,4){$ E_\alpha$ }
\put(1.5,1){$E'_\alpha$}
\put(10.5,1){$E'_\alpha$}
\put(8,1){$E'_\alpha$}
\put(2.5,4.1){\vector(1,0){7.5}}
\put(2.5,1.1){\vector(1,0){5}}
\put(9,1.1){\vector(1,0){1}}
\put(1.8,3.6){\vector(0,-1){2}}
\put(10.7,3.6){\vector(0,-1){2}}
\put(6,4.4){$ T^k_\alpha $}
\put(4.2,1.4){$ S^k_\alpha $}
\put(1.8,2.3){$ D_\alpha $}
\put(10.8,2.3){$ D_\alpha $}
\put(9,1.4){$ P_\alpha $}
\end{picture}
\caption{Approximating a Dilation}
\label{fig:dil1}
\end{figure}
\par
Take a Banach-limit LIM on $l^\infty (A, {\mathbb{C}})$ such that
\[
  \liminf f \le \mbox{\rm LIM} f \le \limsup f \quad \forall\, f \in l^\infty 
  (A,{\mathbb{R}})\,.
\]
Denote 
$E' = \Lambda (A,E_\alpha )/N$ 
and define:
\[     \begin{array}{lclcl}
   D:E \rightarrow E'& \mbox{\rm by} & f \mapsto [D_\alpha \mathcal{E}_\alpha
   f] & & \\
   S:E' \rightarrow E' & \mbox{\rm by} & Sf = [S_\alpha f_\alpha] & \mbox{\rm if} &
       f = [f_\alpha ]\\
   P:E' \rightarrow E' & \mbox{\rm by} & Pf = [P_\alpha f_\alpha ] &
      \mbox{\rm if} & f = [f_\alpha].
   \end{array}
\]
By what has been proved,
\[ 
P_\alpha S^n_\alpha D_\alpha \mathcal{E}_\alpha f = D_\alpha T^n_\alpha \mathcal{E}_\alpha f =
\left\{    \begin{array}{lll}
     D_\alpha T^n_\alpha f & n = 1,2, \ldots & \\
     D_\alpha \mathcal{E}_\alpha f & n = 0 & \forall \; \alpha \in A\,  ,
            \end{array} \right.
\]
and, since 
$D$ 
is an isometry and 
$ \mathcal{E}_\alpha$ a norm non-increasing projection:
\begin{eqnarray*}
\norm{ D_\alpha T^n_\alpha f - D_\alpha \mathcal{E}_\alpha T^nf } & = & \norm{ T^n_\alpha f
   - \mathcal{E}_\alpha T^n f } \\
   & \le & \norm{ T^n_\alpha f - T^n f } \rightarrow 0 \; \mbox{\rm with}\;
   \alpha \in A\,.
\end{eqnarray*}
But this shows
\[ 
   PS^nDf = DT^n f\ , \qquad n = 1,2,\ldots \quad f \in E,
\]
which concludes the proof.
\qed
\begin{rem} \rm
\label{rem:subpos}
\index{Coifman}%
Coifman,
\index{Rochberg}%
Rochberg and
\index{Weiss}%
Weiss~\cite{CoRoWe_77} extended the result of \\
Ak\c{c}oglu and Sucheston to cover the case
of a 
\index{contraction!sub-positive}%
sub-positive contraction $T$
(see the pages 58 and 59 of \cite{CoRoWe_77}).
\end{rem}
For the readers convenience we recall the definition. 
\par
For
$U:L^p(\Omega, \mu ; {\mathbb{C}}) \rightarrow L^p(\Omega, \mu ; {\mathbb{C}})$, 
bounded and linear, define\\
\index{$\overline{U}$}%
$\overline{U} \in B(L^p(\Omega, \mu ; {\mathbb{C}}))$ 
by
\[
  \overline{U}f = \overline{(U{\bar f})} \quad f \in L^p(\Omega, \mu ; {\mathbb{C}})
\]
($^-$ 
denoting complex conjugation),
and define further
$ {\rm Re}\ U \in B(\Lp{\Omega}{\mu;{\mathbb{C}}} $ by
\[
    {\rm Re}\ U = \frac{1}{2}(U + \overline{U})\,.
\]%
~\index{${\rm Re}\ U $}
\begin{defi}
\label{defi:subposcontr}
$T \in {B}(L^p(\Omega, \mu ; {\mathbb{C}}))$ 
is called a 
\index{sub-positive}%
\index{contraction!sub-positive}%
sub-positive
contraction, if there exists a positive contraction 
$R$ 
such that 
$R + {\rm Re}(e^{i\theta}T) $
is positive for all $\theta \in [0,2\pi)$.
\end{defi}
\begin{prop}
\label{prop:subposmatr}
If 
$T \in B(l^p_n)$ 
is represented by its matrix 
$(T_{ij})^n_{i,j=1}$ 
with
respect to the usual basis of 
$l^p_n$, 
then 
$T$ 
is a
\index{contraction!sub-positive}%
sub-positive contraction
if and only if 
$\left( \abs[T_{ij}] \right)^n_{i,j=1}$ 
represents a contraction.
\end{prop}
{\bf Proof:}
Let $R$ be a contraction such that
\[ 
    R + {\rm Re}\ (e^{i \theta}T) \ge 0 \quad \forall\; \theta\,.
\]
Then
\[   
 R_{ij} +({\rm Re}\ (e^{i \theta}T))_{ij} =  R_{ij} + {\rm Re}\ (e^{i \theta }T)_{ij} \ge 0 \quad \forall\, i,j,\theta.
\]
Hence 
$R_{ij} - \abs[T_{ij}] \ge 0 \quad \forall\, i,j$, 
and by Proposition~\ref{prop:posmap} 
it follows that the matrix 
$\left(\abs[T_{ij}]\right)_{i,j=1}^{n}$ 
represents a contraction which we  denote 
$\abs[T]$.
\par
The other implication is clear, since we can choose 
$\abs[T]$
as the positive contraction, named 
$R$,
in the definition of sub-positivity.
\qed
\par
for a sub-positive contraction $T$
a modification of the proof of the dilation theorem has to be done in the
finite dimensional case.
For the construction then, Coifman and Weiss  use 
$\abs[T]$ instead of 
$T$ to obtain the
measure space 
$L^p(\Omega', \mu')$, 
 the affine mappings 
$\tau_{ij}$ 
and the
Radon-Nikodym derivative 
$\rho$.
\par
Now the operator 
$S$ 
is defined as
\[
  S\ f(x,y) = \sum_{i,j}\ \sigma_{ij}\ \rho (x,y)^{\frac{1}{p}} 
  f(\tau^{-1}_{ij}(x,y))\chi_{S_{ij}}x,y), \qquad (x,y) \in \Omega',
\]
where 
$\sigma_{ij} = \mbox{\rm sign}\; T_{ij}$.\\
\par
\begin{defi}
We shall say that an operator 
$R : \Lp{\Omega}{\mu} \rightarrow \Lp{\Omega'}{\mu'}$ 
\index{operator!separation preserving}%
\index{separation!preserving}%
preserves separation if
$f \cdot g = 0 $ implies $ Rf \cdot Rg = 0 \quad \forall\; f,g \in \Lp{\Omega}{\mu}$, 
i.e., if the supports of two elements are disjoint, then so are
the supports of their images under 
$R$.
\end{defi}
\begin{rem}~\\[-2em]
\label{rem:isometry}
\begin{description}
\item[(i)~]For %
$1 \leq p < \infty, \; p \neq 2$
it is easy to see that an isometry 
$S : \Lp{\Omega}{\mu} \rightarrow  \Lp{\Omega'}{\mu'}$ 
is 
\index{isometry!separation preserving}%
separation preserving.\\
\item[(ii)]In the case that $p=2$ the appropriate requirement is that $S$ is 
a positive isometry.
\end{description}
\end{rem}
For the readers convenience we prove the assertions of this  remark.
\par
First assume $1 \leq p < \infty, \; p \neq 2$.
If 
$f \cdot g = 0$ 
then 
$\norm{f + g }^p + \norm{f-g}^p = 2(\norm{f}^p + \norm{g}^p)$
and conversely. 
Now, if 
$S$ 
is an isometry, then this equality holds true for 
$f$ and $g$
if and only if it holds true for the images 
$Sf$ and $Sg$.
\par
Now the first part is proved and to establish the second one we proceed
similarly:
\par
For non-negative $f,g \in L^2(\Omega,\mu)$ 
(resp. $f,g \in L^2(\Omega',\mu')$) we
have that 
$\norm[2]{f\, + \, g}^2 = \norm[2]{f}^2 \, + \, \norm[2]{g}^2 $
if and only if 
$f \cdot g \, = \, 0$.
Hence  for such functions the same reasoning as above applies.
\par
When $f$ is real valued, then $Sf^+ = S(f^+) $ and  $Sf^{-} = S(f^{-})$, since
\linebreak[4] 
$S(f^+) \cdot S(f^{-}) = 0 $ and $ S(f^+), \, S(f^{-}) \, \geq 0 $.
Now it is easy to see that ${\rm Re}\ (Sf) = S ({\rm Re}\ f) $ and ${\rm Im}\ (Sf) = S ({\rm Im}\ f)$ for
a possibly complex valued element $f \in L^2(\Omega, \mu)$.
\par
Since $ f \cdot g \, = \, 0$ exactly if 
all products of one element from
\linebreak[4]
$ \{ ({\rm Re}\ f )^+,\ ({\rm Re}\ f )^{-},\ ({\rm Im}\ f )^+ ,\ ({\rm Im}\ f )^{-} \}$ 
with an element from
\linebreak[4]
$ \{ ({\rm Re}\ g )^+,\ ({\rm Re}\ g )^{-},\ ({\rm Im}\ g )^+ ,\ ({\rm Im}\ g )^{-} \}$ 
vanish one easily concludes using the properties of $S$ established in the
last paragraph. \qed
\begin{lem}
\label{lem:disjoint}
A 
\index{contraction!separation preserving}%
separation preserving contraction
$R$ 
is a sub-positive contraction.
\end{lem}
{\bf Proof:}
For $A \subset \Omega, \;  0 < \mu (A) < \infty$ 
define 
$h_A(\omega) = R(\chi_A)(\omega)$.
\par 
We claim that for
$\omega \in \Omega$ 
and two sets 
$B, C$
of finite nonzero measure 
one of the following cases occurs:
\begin{enumerate}
\item
$ h_C(\omega) =0 \mbox{ and } h_B(\omega) = 0 $
\item
$ h_C(\omega) =0 \mbox{ and } h_B(\omega) \not= 0 $
\item
$ h_C(\omega) \not=0 \mbox{ and } h_B(\omega) = 0 $
\item
$ h_C(\omega) = h_B(\omega) \not= 0 $
\end{enumerate}
To establish the claim, we may assume that
$ R\chi_B(\omega) \cdot R\chi_C (\omega) \not= 0$.\\ 
Since $R$ preserves separation we have
$R\chi_{B-B\cap C}(\omega) \cdot R\chi_{C-B\cap C}(\omega) = 0$,
and from
\begin{eqnarray*}
R \chi_{B \cap C}(\omega) & = & R\chi_B (\omega) - R\chi_{B-B\cap C}(\omega)\\
  & = & R\chi_C(\omega) - R\chi_{C-B\cap C}(\omega)\ 
\end{eqnarray*}
we obtain
\begin{eqnarray*}
R\chi_B (\omega)=R\chi_{B\cap C}(\omega) 
& \mbox{ or } & R\chi_C(\omega)=R\chi_{B\cap C}(\omega).
\end{eqnarray*}
In the first case
$R\chi_B (\omega) - R\chi_{B-B\cap C}(\omega) =  R\chi_C (\omega) \not= 0$.
Again since $R$ preserves separation, we may use
$R\chi_{B-B\cap C}(\omega) \cdot R\chi_{C}(\omega) = 0$
to obtain
\linebreak[4]
$R\chi_{B-B\cap C}(\omega)=0$,
and hence
\[
R\chi_B (\omega)=R\chi_{C}(\omega). 
\]
The second case yields by symmetry the same result, which establishes the
claim.
\par
For 
$\omega \in \Omega$
we may thus define unambiguously
\[
h(\omega) \, = \,
\left\{ 
\begin{array}{ll}
h_A(\omega) & \mbox{ if } A \mbox{ is such that } 0< \mu(A) < \infty \mbox{
  and } h_A(\omega) \not= 0 \\
0 & \mbox{ if for all sets of finite measure } A: \, h_A(\omega) =0.
\end{array}
\right.
\]
This function might not be measurable, but for any set 
$B$
of finite nonzero measure (denote 
$\tilde{B}=\mbox{\rm supp}R{\chi_{B}}$ ) the function
\( h \cdot \chi_{\tilde{B}} = R{\chi_B}
\) 
is measurable.
\par
If 
$(A_j)^k_{j=1}$ 
are pairwise disjoint and 
$ f = \Sigma\ \lambda_j \chi_{A_j}$
is a 
\index{simple function}%
\index{function!simple}%
simple function, then
\begin{eqnarray}%
\label{eq:seppreserv}
Rf & = & h \sum^n_{j=1}\ \lambda_j\ \chi_{\tilde{A_j}}\,.
\end{eqnarray}
If we define 
\begin{eqnarray*} 
P & : & \Lp{\Omega}{\mu} \rightarrow \Lp{\Omega}{\mu}
\end{eqnarray*}
by
\begin{eqnarray}%
\label{eq:seppreserv1}
  Pf(\omega) = |h(\omega)|\ \sum^n_{j=1}\ \lambda_j 
  \chi_{\tilde{A_j}} (\omega)\, ,
\end{eqnarray}
then we see that for all simple functions $f$:
\begin{eqnarray*}
 \norm{Pf }  & = &\norm{ Rf }, \\   
   Pf + {\rm Re}\ e^{i \theta} Rf & \ge & 0 \mbox{ \rm if }  f \ge 0\,.
\end{eqnarray*}
Using the boundedness of $R$ and $P$ this extends to all of
$\Lp{\Omega}{\mu}$.
\qed
\begin{rem}
\label{rem:seppreserv}
\rm
For later use it is worthwhile to note that for a separation preserving
contraction
$R$ there exists a positive contraction $P$, such that for all 
$f \in \Lp{\Omega}{\mu}$:
\[
\abs[Rf](\omega) = Pf(\omega),\qquad \mu \mbox{-almost everywhere}.
\]
This follows from (\ref{eq:seppreserv}) and (\ref{eq:seppreserv1}),
whenever $ f = \Sigma\ \lambda_j \chi_{A_j}$
is a simple function as above,
since the sets $\tilde{A_j}, \, j=1, \ldots , k $ are pairwise disjoint.
By continuity it continues to hold for arbitrary $f \in \Lp{\Omega}{\mu}$.
\end{rem}
\section[A further Ultraproduct Construction]{A further Ultraproduct Construction}
\label{sec:furthUlpr}
Only later, in section~\ref{sec:condil}, we will proof a dilation
theorem, see Theorem~\ref{thm:dilthm}, for certain continuous semigroups on 
$ \mathcal{L}^p$-spaces. For its proof we will use once more  a 
Banach space ultraproduct construction, applied to dilations, 
which are given by the 
Ak\c{c}oglu-Sucheston dilation theorem, Theorem~\ref{thm:AkcoSuchdilthm}.
We shall need some considerations on representations of
continuous groups on reflexive spaces too, which will be presented in 
section~\ref{sec:GenRem}, but we found it best to present the ultraproduct part of the
proof of Theorem~\ref{thm:dilthm} here.
\par
To fix ideas, in this section we shall consider 
a continuous one-parameter semigroup 
$ \{T_t : t \geq 0 \}$,
with
$ T_0  = \id $,
of positive contractions acting on
$ \Lp{\Omega}{\mu} $, where
$ 1\leq  p < \infty $.
We shall prove a preliminary lemma on the existence of a
dilation of the sub-semigroup
$ \{T_t : t \geq 0 , \ t \in {\mathbb{Q}}\}$
to a group of isometries of some 
$ \mathcal{L}^p$-space.
(Here
$ \mathbb{Q} $ 
denotes the rational numbers.)
\par
Furthermore we remark that, as almost 
always in this paper, the letters
$ T, S, $
$P, D $ 
stand for 
\index{positive}%
positive operators between
$ \mathcal{L}^p $-spaces,
\index{$D$}%
$D$
will denote an isometric embedding,  
\index{$P$}%
$P$
a contractive projection, 
\index{$T$}%
$T$
will be a contraction and 
\index{$S$}%
$S$
an invertible isometry.
\begin{lem}
\label{lem:ulprlem}
If 
$ \{ T_t : t \geq 0 \} $
is a semigroup in $ B(\Lp{\Omega}{\mu}) $, 
which fulfils the above requirements, then there exists a measure space 
$ (\Omega^{\circ}, \mathfrak{A} ^{\circ}, \mu^{\circ}) $ 
\index{${(\Omega^{\circ}, \mathfrak{A} ^{\circ}, \mu^{\circ})} $}%
and  a group
\index{$S_{s}^{\circ}$}%
$ \{ S_{s}^{\circ} : s \in \mathbb{Q} \} $ 
of positive isometries of 
\index{$\Lp{\Omega^{\circ}}{\mu^{\circ}} $}%
$ \Lp{\Omega^{\circ}}{\mu^{\circ}} $ 
such that 
\begin{displaymath} 
D^{\circ} \circ T_{s} = P^{\circ} \circ S_{s}^{\circ} 
\circ D^{\circ},  \hspace{10mm}  s \in \mathbb{Q}^+,
\end{displaymath}
for a suitable positive isometric embedding
$ D^{\circ} : \Lp{\Omega}{\mu} \rightarrow \Lp{\Omega^{\circ}}{\mu^{\circ}} $
and a positive contractive projection
$ P^{\circ} $ acting on $ \Lp{\Omega^{\circ}}{\mu^{\circ}} $.  
\end{lem}
In our paper~\cite{Fendler_97a} the proof of this lemma depends on certain
filters on the rational numbers. Here instead, we shall use the partial 
ordering of the natural numbers given by 
\index{order!divisibility}%
divisibility.
\begin{defi}
\label{defi:preceeds}
For 
$n,m \in {\mathbb{N}}$
define
\index{$n\preceq m$}%
\[
n \preceq m \mbox{ if } n \mbox{ divides } m. 
\]
\end{defi}
\begin{rem} \rm \hspace*{\fill}~\\[-2em]
\label{rem:directed}
\begin{description}
\item[(i)~~]
If we take for 
$k,l \in {\mathbb{N}}$
$m=\mbox{\rm lcm}(k,l)$
 to be the least 
\index{multiple!least common}%
common multiple of $k$ and $l$, then
\[ k \preceq m \mbox{ and } l \preceq m. \]
Thus 
$({\mathbb{N}}, \ \preceq )$
is directed.
\item[(ii)~]
For a finite set 
$ B \subset \mathbb{Q} $ 
let 
$ U_B := \{ n \in \mathbb{N} : ns \in \mathbb{Z} \hspace{3mm} \forall s \in B \} $.
The latter set just consists of the common multiples of the denominators 
of the reduced fractions representing the rationals from
$B$.
\par
The fact that the set of all sets
$ \{ U_B : B \subset \mathbb{Q} \, , \hspace{3mm} B \mbox{ finite } \} $ 
is  closed under 
\index{intersection!finite}%
finite intersections just corresponds to
$({\mathbb{N}}, \ \preceq )$
being directed.
\item[(iii)]
We note that the Banach limits given by Proposition~\ref{prop:Blimex}~for 
$ ({\mathbb{N}}, \ \preceq )$
are in one to one correspondence to the 
\index{Banach limit}%
\index{filter!maximal}%
maximal filters
containing the above 
\index{filterbasis}%
filter-basis.
\end{description}
\end{rem} 
{\bf Proof of the Lemma~\ref{lem:ulprlem}:} \rm  
 According to Theorem~\ref{thm:AkcoSuchdilthm}, for any 
$ n \in \mathbb{N} $, 
there exists a dilation of 
$ \{ T^{k}_{1/n} : k \in \mathbb{Z}^+ \} $: 
\begin{displaymath}
D_{1/n} \circ T^{k}_{1/n} = P_{1/n} \circ S^{k}_{1/n} 
\circ D_{1/n} \, , \hspace{10mm}  k \in \mathbb{Z}^+,
\end{displaymath}
where 
$ D_{1/n},  P_{1/n}, S_{1/n} $ 
are as above, 
$ S_{1/n} $ 
acting on some space
$ \Lp{\Omega_{1/n}'}{\mu_{1/n}' } $. \par
If we define, for $ s \in \mathbb{Q} $\,,\,
$ S_{n,s} \in B(\Lp{\Omega_{1/n}'}{\mu_{1/n}'}) $ 
by 
\begin{displaymath} 
S_{n,s} = \left\{ \begin{array}{r@{\quad \quad}l}
 S_{1/n}^{ns}  &  \mbox{if} \quad ns \in \mathbb{Z} \\ 
\id & \mbox{if} \quad ns \notin \mathbb{Z},
\end{array} \right.  
\end{displaymath}
and if 
$ B = \{s_1,\dots ,s_k \} \subset \mathbb{Q^+} $ 
is a finite subset, then 
for 
$ s \in B $ 
and 
$ n \in U_B $ 
the diagram in Figure~\ref{fig:dil2} 
commutes.
\begin{figure}
\unitlength1cm
\begin{picture}(13,5)
\put(0.5,4){$ \Lp{\Omega}{\mu} $ }
\put(10.5,4){$ \Lp{\Omega}{\mu} $ }
\newsavebox{\raum}
\savebox{\raum}{$ \Lp{\Omega_{1/n}'}{\mu_{1/n}'} $ }
\put(0.5,1){\usebox{\raum}}
\put(10.5,1){\usebox{\raum}.}
\put(6,1){\usebox{\raum}}
\put(2.5,4.1){\vector(1,0){7.5}}
\put(3.5,1.1){\vector(1,0){2}}
\put(9,1.1){\vector(1,0){1}}
\put(1.5,3.6){\vector(0,-1){2}}
\put(11.5,3.6){\vector(0,-1){2}}
\put(6,4.4){ $ T_{1/n}^{ns} $ }
\put(4.2,1.4){ $ S_{n,s} $ }
\put(1.8,2.3){ $ D_{1/n} $}
\put(11.8,2.3){ $ D_{1/n} $}
\put(9,1.4){ $ P_{1/n} $ }
\end{picture}
\caption{Dilations for Rationals}
\label{fig:dil2}
\end{figure}
\par
Let 
\index{LIM}%
LIM
be a Banach limit for 
$({\mathbb{N}}, \ \preceq)$. Then
we may form 
\index{ultraproduct}%
ultraproducts of the spaces and operators involved.
From Corollary~\ref{cor:Lpulpr} we know
that there exist measure spaces 
$ (\Omega^{\circ} , \mathfrak{A}^{\circ}, \mu^{\circ}) $ 
and 
$ (X^{\circ},\mathfrak{B}^{\circ},\nu^{\circ}) $ 
such that
\begin{eqnarray*}
\Ulp{\BLIM}{\Lp{\Omega}{\mu}} & = &  \Lp{X^{\circ}}{\nu^{\circ}} \\
\Ulp{\BLIM}{\Lp{\Omega_{1/n}'}{\mu_{1/n}'}} & = & \Lp{\Omega^{\circ}}{\mu^{\circ}}
\end{eqnarray*} 
can be identified as Banach lattices in each of the two cases. \par
Let I denote the canonical inclusion 
$ I : \Lp{\Omega}{\mu} \rightarrow \Ulp{\BLIM}{\Lp{\Omega}{\mu}} = 
\Lp{X^{\circ}}{\nu^{\circ}} $ 
and denote
 \begin{eqnarray*}
D^{\circ} & = & \Ulp{\BLIM}{ D_{1/n}} \, \circ I \; , \\
P^{\circ} & = &\Ulp{\BLIM}{ P_{1/n}} \; , \\
S^{\circ}_{s} & = & \Ulp{\BLIM}{ S_{n,s}} \; , \quad s \in \mathbb{Q} \,. 
\end{eqnarray*}  
Then all the asserted properties of the operators and spaces involved
are rather immediate. 
\par
To give examples,  let us check that  
\begin{eqnarray*}
 S : \mathbb{Q} &  \rightarrow & B \left(\Ulp{\BLIM}{\Lp{\Omega_{1/n}'}{\mu_{1/n}'}} \right), 
\quad s \mapsto S_s \,, 
\end{eqnarray*} 
is a group homomorphism, and that we obtained a dilation of the semigroup
$ \{ T_s : s \in \mathbb{Q}^+ \} $.
 \par
If 
\mbox{ $ s, s' \in \mathbb{Q} $} 
are given, and if \mbox{ $ f \in \Ulp{\BLIM}{\Lp{\Omega_{1/n}'}{\mu_{1/n}'}} $}
is represented by a sequence 
\mbox{ $ (f_n)_{n \in \mathbb{N}} $},
which is possible by Proposition~\ref{prop:complete},
then for $n$ sufficiently large with respect to our partial ordering
$\preceq, $ i.e.\ if
\mbox{ $ n \in U_{ \{ s,s' \} } $,} 
there holds true:
\begin{displaymath}
\begin{array}{r@{\quad = \quad}l}
S_{n, s+s'}(f_n) & S_{1/n}^{n(s+s')}(f_n ) \\
   &  S_{1/n}^{ns}(S_{1/n}^{ns'}(f_n )) \\
   &  S_{n, s}(S_{n,s'}(f_n)) \,.
\end{array}
\end{displaymath}
By the definition of the partial ordering and by Proposition~\ref{prop:Blimex},
\begin{eqnarray*}
\lefteqn{\norm[\BLIM]{(S_{n,s+s'}(f_n ))_ {n \in \mathbb{N}} - (S_{n,s} (S_{n,s'} (
  f_{n} )))_{ n \in \mathbb{N}}} \, \leq}\\
 & \le & \limsup_{n\in {\mathbb{N}}, \preceq}
\norm{S_{n,s+s'}(f_n ) - S_{n,s} (S_{n,s'} ( f_{n} )) } \\
&= & 0.
\end{eqnarray*} 
From this we infer that they represent the same elements in \\
$ \Ulp{\BLIM}{\Lp{\Omega_{1/n}'}{\mu_{1/n}'}} $, 
and thus 
\begin{eqnarray*}
 S^{\circ}_{s} \circ S^{\circ}_{s'} (f ) & = & S^{\circ}_{s+s'} (f). 
\end{eqnarray*} \par 
Furthermore, for $ s \in \mathbb{Q}^+ $, 
the commutativity of the above
diagram for large enough $n$, i.e.\ for 
$ n \in U_{\{s\}}, $ 
implies, as one can see using a reasoning analogous to the one just given, 
\begin{displaymath} \qquad \qquad \qquad \qquad \qquad \quad
D^{\circ} \circ T_{s} = P^{\circ} \circ S_{s}^{\circ} 
\circ D^{\circ}. \qquad \qquad \qquad \qquad \quad \qedm
\end{displaymath}
\chapter[Commutative Groups on 
reflexive Spaces]{Representations of Commutative topological Groups on 
reflexive Banach Spaces}
\section[General Remarks]{General Remarks}
\label{sec:GenRem}
This topic has been discussed in an broader context by 
\index{Glicksberg}%
I.~Glicksberg and 
\index{de~Leeuw}%
K.~de\-~Leeuw\-~\cite{GlicdeLe_65}. Here we shall develop only those parts
of their theory which are necessary for our purpose.
\par
Let $E$ be a reflexive Banach space, 
\index{$G~$}%
$G$ a commutative topological group
and 
\index{$\pi~:~G~\rightarrow~B(E)$}%
$ \pi \, : \, G \rightarrow  B(E)$
a 
\index{representation}%
\index{representation!uniformly bounded}%
uniformly bounded representation of $ G$ on $E$.
Thus we are given an algebraic group homomorphism
$\pi$
from
$G$
to the invertible elements of
$ B(E) $
such that 
$ C:= \sup_{\pi \in G} \norm{\pi(x)} < \infty$. 
\par 
Our topic here is to study how the subspace 
\begin{displaymath}
 E_c = \{ \xi \in E : \; x \mapsto \pi(x) \xi \mbox{ is continuous from } 
G \mbox{ to } E \}
\end{displaymath}
of continuously translating elements of
$ E $ is situated in $E$.
\par
First, since $\pi$ is uniformly bounded, 
\index{$E_{c}$}%
$E_{c}$ is surely a closed subspace of 
$E$ and it is clearly invariant for $\pi(G)$.
\par
If $E$ were a Hilbert space, then one would take an orthogonal projection 
$P$
onto
$E_{c}$. Since $G$ is amenable as a discrete group with an invariant
mean on $l^{\infty}(G)$, $m$ say,
we then could define a new projection $P_{G}$ onto $E_{c}$ in the commutant of
$\pi(G)$ by 
\begin{displaymath}
<\,P_{G} \xi \, , \, \eta \, > \, = \, 
m(x \mapsto <\pi(x^{-1}) P \pi(x) \xi, \eta > )  
\quad \quad \forall \xi, \eta \in E.
\end{displaymath}
Clearly its norm is bounded by $C^{2}$.
\par
In a reflexive space there is in general no bounded projection 
onto a closed subspace and we have to use a refinement of the above 
construction.
\par
Let 
\index{$\mathcal{U}(e)$}%
$\mathcal{U}(e)  $ denote the system of open 
\index{neighbourhood}%
\index{neighbourhood!system}%
neighbourhoods of the
identity 
\index{$e~\in~G$}%
$e \in G$
and let, for $U \in \mathcal{U}(e)$,
\begin{displaymath}  
U_{\pi} = \{ \, \pi(x) \, : \, x \in U \, \}^{-wot}
\end{displaymath}
denote the closure of $\pi(U) $ in the 
\index{topology!weak}%
weak operator topology in $ B(E) $.
Then, since $G$ is abelian, 
\begin{displaymath}
 \Gamma \, = \bigcap_{U \in {\mathcal{U}}(e)} U_{\pi}
\end{displaymath}
is a  set of commuting operators.
When endowed with the weak operator topology, it
is a compact topological space.
Furthermore, if $E$ is given its weak topology,
then the action of 
\index{$\Gamma$}%
$\Gamma$ on $E$, i.~e.\ the map
\begin{eqnarray*}
(S,\xi) & \mapsto & S \xi, \\
\Gamma \times E & \rightarrow & E,
\end{eqnarray*}
is separately continuous. 
\par
\begin{lem}
\label{lem:fixpt} 
$E_c$ coincides with the set of 
\index{$\Gamma$-fixed}%
\index{fixed points!$\Gamma~$}%
$\Gamma$-fixed points in $E$.
\end{lem}
{\bf Proof:}
If $ \xi \in E_c $, then for $\epsilon > 0 $ there exists an 
$ U \in {\mathcal{U}}(e) $ such that 
$\norm{\pi(x)\xi - \xi}  <  \epsilon \mbox{ for all }  x \in U.  
\mbox{ But then}
\norm{S \xi - \xi }  <  \epsilon \mbox{ for all } S \in U_{\pi} $,
and it follows that $ S \xi = \xi \mbox{ for all } s \in \Gamma$. 
\par
On the other hand, if $ \{ x_{\alpha} \}_{\alpha \in I} $ is a net
in $G$ converging to the identity, then the accumulation points of the net of operators
$\{\pi(x_{\alpha}) \}_{\alpha \in I} $, there
exists at least one since $\Gamma$ is compact, 
are in $\Gamma$.
\par 
We claim that for  any $ \xi \in E $ which is fixed by all elements of $\Gamma$ we  have
\begin{displaymath} 
\xi = \lim_{\alpha} \pi(x_{\alpha})\xi.
\end{displaymath} 
Establishing the claim will prove the Lemma, and to do so we are going 
to prove that all weak topology 
accumulation points of the net
$\{\pi(x_{\alpha}) \xi \}_{\alpha \in I} $
are contained in
\begin{displaymath}
\bigcap_{U \in {\mathcal{U}}(e)} U_{\pi} \xi \, = \, \Gamma \xi \, = \, \{ \xi \}.
\end{displaymath}
For this let 
$ \eta \in E $
be an accumulation point of 
$\{\pi(x_{\alpha}) \xi \}_{\alpha \in I} $
in the weak topology of $ E $.
Given
$U \in {\mathcal{U}}(e)$, we find, for any weak topology neighbourhood
$W$ of $\eta$, some $\alpha$, such that
\[
x_{\alpha} \in U \mbox{ and } \pi(x_{\alpha})\xi \in W.
\]
Hence,
\[
U_{\pi}\xi \cap W \neq \emptyset
\] for any such $W$ and we infer, since $U_{\pi}\xi$ is closed as the
continuous image of the compact set $U_{\pi}$, that $\eta \in  U_{\pi}\xi$.
\qed
For $ \xi \in E $ let 
\index{$\mbox{conv}\{\pi(U)\xi\}^-$}%
$ \mbox{conv}\{ \pi(U) \xi \}^- $ be the norm closure 
of the convex hull of 
$ \pi(U) \xi $.
This  set is bounded, weakly closed because of its convexity, and hence 
weakly compact in the reflexive space $ E $. 
\par
Let
\[ C_{\xi} := \bigcap_{U \in {\mathcal{U}}(e)} \mbox{conv}\{ \pi(U) \xi \}^-. \]
This then is a non-void, convex, weakly compact and 
$ \Gamma $-invariant subset of
$ E $ reducing to 
\begin{eqnarray*}
 C_{\xi} & \, = \, & \{ \xi \} \mbox{ if } \xi \in E_c.
\end{eqnarray*}
Furthermore
\begin{eqnarray}
\label{eq:homogen}%
C_{\alpha \xi } & = & \alpha C_{\xi},~~~~~~~~~~~~~~\alpha \in \mathbb{C},
 \xi \in E,\\
\label{eq:additive}%
C_{\xi + \eta} & \subset & C_{\xi} + C_{\eta},~~~~~~~~~~~~\xi , \eta
\in E.
\end{eqnarray}
For 
$\xi \in E $
the 
\index{theorem!fixed-point}%
\index{theorem!Markov-Kakutani}%
Markov-Kakutani fixed-point theorem, 
see e.g. Chap. IV, Appendix, Theorem~1 of \cite{Bourbaki_87}, may be applied 
to the action of
$ \Gamma \mbox{ on } C_{\xi} $.
Hence there exists at least one point in 
$ C_{\xi} $
fixed by all
$ S \in \Gamma$.
This, by Lemma~\ref{lem:fixpt},  shows 
\begin{eqnarray*}
 C_{\xi} \cap E_c & \, \not= \, & \emptyset.
\end{eqnarray*}
\par
We claim that 
$ C_{\xi} \cap E_c $ 
contains exactly one element.\\
If there are
$\, \eta,\zeta \in C_{\xi} \cap E_c\,$,
then for any 
$ \epsilon > 0$
there exists 
$ U_0 \in {\mathcal{U}}(e)$ 
such that for all $ x \in U_0 $
\[
\begin{array}{rcll}
\norm{\pi(x)\eta - \eta}{}  <  \epsilon & \mbox{and} & 
\norm{\pi(x)\zeta - \zeta}{} <  \epsilon  .
\end{array}
\]
By the definition of $C_{\xi}$  there are approximations
$ \tilde{\eta}, \tilde{\zeta} \in \mbox{conv}\{ \pi(U_0) \xi \} $
such that
\[
\begin{array}{rcl}
\norm{\tilde{\eta} - \eta}{}  <  \epsilon & \mbox{and} & 
\norm{\tilde{\zeta} - \zeta}{} <  \epsilon.  
\end{array}
\]
We may write 
$ \tilde{\eta} = \sum_{i=1}^n \lambda_i \pi(y_i) \xi $
and 
$\tilde{\zeta} = \sum_{j=1}^m \mu_j \pi(z_j) \xi  $, where
\linebreak[4] 
$ y_1, \dots, y_n,z_1, \dots, z_m \in U_0, $
and
$\lambda_1, \dots, \lambda_n > 0,\, \mu_1, \dots, \mu_m > 0  $
fulfil 
$  \sum_{i=1}^n \lambda_i = \sum_{j=1}^m \mu_j =1$.
\par
Then
\begin{eqnarray*}
\lefteqn{\norm{ \eta - \sum_{j=1}^m \mu_j \pi(z_j) \tilde{\eta}}{}}\\
 &\leq & \norm{ \eta - \sum_{j=1}^m \mu_j \pi(z_j) \eta}{} + 
\norm{\sum_{j=1}^m \mu_j \pi(z_j)(\eta - \tilde{\eta})}{}\\
 & \leq & \sum_{j=1}^m \mu_j \sup_j \norm{\eta - \pi(z_j) \eta}{} +
\sum_{j=1}^m \mu_j \sup_j \norm{\pi(z_j)}{} \norm{\eta - \tilde{\eta}}{}\\
 & \leq & \epsilon + C \epsilon.
\end{eqnarray*}
Similarly we obtain
\[
 \norm{ \zeta - \sum_{i=1}^n \lambda_i \pi(y_i) \tilde{\zeta}}{} 
\leq (1+C)\epsilon.
\]
\par
Since $ G $ is commutative
\[
\sum_{j=1}^m \mu_j \pi(z_j) \tilde{\eta}  = 
\sum_{j=1}^m \sum_{i=1}^n \mu_j \lambda_i \pi(z_j y_i) \xi =
\sum_{i=1}^n \lambda_i \pi(y_i) \tilde{\zeta} 
\]
coincide and we infer 
\[ \norm{\eta - \zeta}{} \leq 2(1+C) \epsilon, \] 
which proves the claim since $ \epsilon > 0 $ is arbitrary. 
\vspace{5mm}
\par
Now we have almost proved the following proposition.
\begin{prop}
\label{prop:projection}
Let $ E $ be a reflexive Banach space, $G$ a commutative topological group and 
$ \pi : G \rightarrow B(E) $ a uniformly bounded representation of $ G $ as 
above. Then there exists a  projection 
\index{$E_{c}$}%
\index{$Q~$}
\index{projection!$Q~\in~B(E)$}%
$ Q \in B(E) $ with range $ E_c $ 
such that: \rm
\begin{description}
\item[(i)~~]{~~~~~ $\pi (x) Q  =  Q \pi (x),~~~~~~ x \in G$,} 
\item[(ii)~]{~~~~~ $\norm{Q}{}  \leq  \sup_{x \in G} \norm{\pi(x)}{}$,}
\item[(iii)]{~~~~~ $  Q \xi  \in  C_{\xi},~~~~~~~~~~\xi \in E$.}
\end{description}
\end{prop}
{\bf Proof:}
For all 
$  \xi \in E $ we know that $E_c \cap C_{\xi} $ is a one point set,
and hence we may define a map 
$ Q : E \rightarrow E $
by requiring for $  \xi \in E $
\[ \{Q \xi \} = E_c \cap C_{\xi}. \]
By its definition it is clear that 
$ Q\xi \in C_{\xi} \; \mbox{ for all } 
\xi \in E$ and from the equation (\ref{eq:homogen}) the homogeneity
of $Q$ is obvious.
\par
To see that $Q$ is additive, we claim that for $\xi$ and $\eta$ in $E$ we
have 
\begin{eqnarray*}
  Q\xi \, + \, Q\eta & \in & C_{\xi + \eta},
\end{eqnarray*}
since then we can infer from the uniqueness that
\begin{eqnarray*}
   Q\xi \, + \, Q\eta & = & Q(\xi +\eta ).
\end{eqnarray*}
To establish the claim, it suffices to show that
for arbitrary $U \in {\mathcal{U}}(e)$
\[
 Q\xi \, + \, Q\eta \in \mbox{conv}\{ \pi(U) (\xi +\eta )\}^-.
\]
To do this, denote $C=\sup_{x \in G} \norm{ \pi(x) }{}$, and choose
for $\epsilon >0$ 
a symmetric
$V \in {\mathcal{U}}(e)$ 
such that 
$V^{2} \subset U$
and
$\norm{Q\xi -\pi (x)Q\xi } < {\epsilon}$, 
$\norm{Q\eta  - \pi (x)Q\eta } < {\epsilon}  $
for all \hfill $ x \in V$.
\hfill
Now \hfill there \hfill exist \hfill approximations \hfill
$\,  \tilde{\xi} = \sum_{i=1}^n \lambda_i \pi(y_i) \xi \,$
\hfill and
\linebreak[4] 
$\tilde{\eta} = \sum_{j=1}^m \mu_j \pi(z_j) \eta \, $
such that 
\[
\norm{Q\xi \, - \, \tilde{\xi}}{} < \frac{\epsilon}{C} \quad \mbox{and} \quad
\norm{Q\eta \, - \, \tilde{\eta}}{} < \frac{\epsilon}{C}.
\] 
Here again
$ y_1, \dots, y_n,z_1, \dots, z_m \in V, $
{ and } 
$\lambda_1, \dots, \lambda_n > 0$
$ \mu_1, \dots, \mu_m > 0  $
fulfil 
$  \sum_{i=1}^n \lambda_i = 1
 \mbox{ and } \sum_{j=1}^m \mu_j =1$.
We note that 
\[
\sum_{j=1}^m \sum_{i=1}^n \mu_j \lambda_i \pi(z_j y_i) (\xi\, + \, \eta )
\, \in \mbox{conv}\{ \pi(V^2) (\xi\, + \, \eta ) \} \, \subset \,
\mbox{conv}\{ \pi(U) (\xi\, + \, \eta ) \},
\]
and estimate, using the commutativity of $G$ again:
\begin{eqnarray*}
\lefteqn{\norm{Q\xi \, + \, Q \eta \, - \, \sum_{j=1}^m \sum_{i=1}^n \mu_j \lambda_i
  \pi(z_j y_i) (\xi\, + \, \eta)}{} }& & \\
& \leq &  
\norm{ Q\xi \, - \,  \sum_{j=1}^m \mu_j \pi(z_i) Q\xi }{}
\, + \, \norm{ \sum_{j=1}^m \mu_j \pi(z_j) \left(Q\xi \, - \, \sum_{i=1}^n
  \lambda_i \pi(y_i) \xi\right)}{}\\
& & + \, \norm{  \sum_{i=1}^n \lambda_i \pi(y_i) \left(Q\eta \, - \, \sum_{j=1}^m
  \mu_j \pi(z_j) \eta\right)}{}
\, + \, \norm{ Q\eta \, - \, \sum_{i=1}^n \lambda_i \pi(y_i)Q\eta}{}\\
& \leq &
\epsilon \, + \, \sum_{j=1}^m \mu_j \, C \, \norm{Q\xi \, - \, \tilde{\xi}}{}
\, + \, \sum_{i=1}^n \lambda_i  \, C \, \norm{Q\eta \, - \, \tilde{\eta}}{}  
\, + \, \epsilon \, \leq \, 4 \epsilon.
\end{eqnarray*}
\par
Since $Q^2\,=\,Q$ is evident, it is now proved that 
$ Q $ is a linear projection onto $E_c$.
\par
That
$ Q $
is in the commutant of
$  \pi(G) $ 
is implied by 
\[ C_{\pi(x)\xi} = \pi(x) C_{\xi}, \quad  x \in G , \; \xi \in E \]
which again is a consequence of the  commutativity of $ G $. 
\par
Finally the norm estimate {\rm (ii) } is obvious from
\[ \sup\{ \, \norm{\eta} : \; \eta \in C_{\xi} \, \} \, \leq \,  
\sup_{x \in G} \norm{ \pi(x) }{} \norm{\xi}{}, \quad  \xi \in E. \]
\qed
\section[Application to Dilations of continuous one-parameter
  Semigroups]{Application to Dilations of continuous \\one-parameter
  Semigroups}
\label{sec:condil}
 In this section we shall collect our foregoing considerations to prove 
the following dilation theorem for continuous one-parameter semi\-groups:
\begin{thm}
\label{thm:dilthm}
\index{theorem!dilation!semigroup}%
Let 
$ \{ T_t : t \geq 0 \} $,
with
$ T_0  = \id $,
be a strongly continuous semigroup of positive contractions acting on
$ \Lp{\Omega}{\mu} $, where
$ 1< p < \infty $. 
Then there exists a measure space 
$ (\tilde{\Omega} , \tilde{\mathfrak{A}} , \tilde{\mu} ) $
and a strongly continuous group of 
\index{group!isometries}%
isometries
\index{$\{S_s:s\in\mathbb{R}\}$}%
$ \{ S_s : s \in \mathbb{R} \} $, with $ S_0= {\rm id}  $, 
acting on 
\index{$\Lp{\tilde{\Omega}}{\tilde{\mu}}$}%
$ \Lp{\tilde{\Omega}}{\tilde{\mu}} $
such that
\begin{displaymath}
D \circ T_t = P \circ S_t \circ D , \hspace{10mm}\mbox{$  t \geq 0 $},
\end{displaymath}
where $ D $ is an isometric embedding of
$ \Lp{\Omega}{\mu} $ into $ \Lp{\tilde{\Omega}}{\tilde{\mu}} $ 
and \\ 
$ P : \Lp{\tilde{\Omega}}{\tilde{\mu}} \rightarrow \Lp{\tilde{\Omega}}{\tilde{\mu}} $
is a contractive projection; further
\index{$D$}%
$ D $,
\index{$P$}%
$ P $ 
and 
the isometries
$ \{ S_s : s \in \mathbb{R} \} $
can be chosen to be positive.
\end{thm}
Given a semigroup as in the theorem
\begin{displaymath}  T_t \, : \, \Lp{\Omega}{\mu} \rightarrow \Lp{\Omega}{\mu},
\quad t \geq 0 ,
\end{displaymath}
we obtain for all 
$s \in \mathbb{Q}_{+}$ the commutative diagram displayed in Figure~\ref{fig:dil3}.
\begin{figure}
\unitlength1cm
\begin{picture}(13,7)
\put(0,6){$ \Lp{\Omega}{\mu} $ }
\put(11,6){$ \Lp{\Omega}{\mu} $ }
\newsavebox{\rauma}
\savebox{\rauma}{$ \Lp{\Omega^{\circ}}{\mu^{\circ} } $ }
\put(0,3){\usebox{\rauma}}
\put(11,3){\usebox{\rauma}}
\put(7,3){\usebox{\rauma}}
\newsavebox{\raumb}
\savebox{\raumb}{$ \Lp{\Omega^{\circ}}{\mu^{\circ} }_{c} $ }
\put(0,0){\usebox{\raumb}}
\put(7,0){\usebox{\raumb}}
\put(2,3.1){\vector(1,0){4.5}}
\put(1,5.6){\vector(0,-1){2}}
\put(12,5.6){\vector(0,-1){2}}
\put(2,6.1){\vector(1,0){8.5}}
\put(1,2.6){\vector(0,-1){2}}
\put(8,0.6){\vector(0,1){2}}
\put(2,0.1){\vector(1,0){4.5}}
\put(9,3.1){\vector(1,0){1.75}}
\put(5,6.4){ $ T_{s}$ }
\put(4,3.4){ $ S_{s}^{\circ} $ }
\put(1.3,4.8){ $ D^{\circ} $}
\put(12.3,4.8){ $ D^{\circ} $}
\put(4,0.4){ $ S_{s} $ }
\put(1.3,1.8){ $Q $}
\put(8.3,1.8){ ${I}$}
\put(9.75,3.4){$ P^{\circ} $}
\end{picture}
\caption{The Part of Continuity}
\label{fig:dil3}
\end{figure}
\hspace*{\fill}
\par
Here 
$ \{ S^{\circ}_{s} \, :\, s \in \mathbb{Q} \}$
is the group constructed in Lemma~\ref{lem:ulprlem}
$ D^{\circ} \mbox{ and } P^{\circ}$
equally come from there, whereas 
$\Lp{\Omega^{\circ}}{\mu^{\circ}}_{c} \mbox{ and the projection } Q $
are constructed in Proposition~\ref{prop:projection}.
 The map $I$ is just the inclusion of 
the subspace
$\Lp{\Omega^{\circ}}{\mu^{\circ}}_{c}$ and the operators 
$ S_{s} \, = \, S^{\circ}_{s}{|_{\Lp{\Omega^{\circ}}{\mu^{\circ}}_{c}}}$
are just the restrictions of the 
$ \{ S^{\circ}_{s} \, :\, s \in \mathbb{Q} \}$
to this subspace.
\par
The following lemma shows that the range of $ D^{\circ} $ is included
in the set of continuously translating elements of
$\Lp{\Omega^{\circ}}{\mu^{\circ}}. $
\begin{lem}
\label{lem:contlem}
For 
\mbox{$ f \in \Lp{\Omega}{\mu}$}, \quad 
\mbox{$ S^{\circ}_. \circ D^{\circ}(f) : s \mapsto 
S^{\circ}_s \circ D^{\circ}(f) $}
is continuous from 
$ \mathbb{Q} \subset \mathbb{R}$
to 
$ \Lp{\Omega^{\circ}}{\mu^{\circ}} $ 
with its norm topology.
\end{lem}
{\bf Proof:}
Since 
$ \{ S^{\circ} : s \in \mathbb{Q} \} $ 
is a group of isometries, it suffices
to show, for 
\mbox{$ f \in \Lp{\Omega}{\mu} $}, 
the continuity from the right of the map
$ s \mapsto S^{\circ}_s(D^{\circ}(f)) $
at
\mbox{$ s = 0 $}. 
\par
The space 
$ \Ulpr{ {\BLIM}}{ \Lp{\Omega_{1/n}'}{\mu_{1/n}'}}  =  
\Lp{\Omega^{\circ}}{\mu^{\circ}} $  
is 
\index{Banach space!uniformly convex}%
\index{convex!uniformly}%
uniformly convex. Thus to 
$ \epsilon>0 $ 
there exists 
$ \eta_p(\epsilon)>0 $ 
such that for any 
$ f,h \in \Lp{\Omega^{\circ}}{\mu^{\circ}} $ 
of norm one 
$ \norm[p]{ \frac{1}{2} (f + h)} \geq 1-\eta_p(\epsilon)  $ 
implies 
$ \norm[p]{ f - h } \leq \epsilon$. 
\par
Given a norm one element in the range of 
$ D^{\circ}$
this may be written 
$D^{\circ}(f)$
for an
$ f \in \Lp{\Omega}{\mu} $ 
with 
$ \norm[p]{f}=1 $.
Then for
$ \epsilon > 0 $ 
there exists $ \delta > 0 $
such that  
$ \norm[p]{ T_sf - f } \leq 2 \eta_p(\epsilon) $
whenever 
$ 0 \leq s \leq \delta $.
Hence, for $ s \in \mathbb{Q} \cap [0,\delta) $,
\begin{displaymath}
\begin{array}{rcl}
\norm[p]{S_s^{\circ}(D^{\circ}(f)) + D^{\circ}(f)}  & \geq & 
\norm[p]{ P^{\circ} \circ S^{\circ}_s \circ D^{\circ} (f) + P^{\circ} \circ D^{\circ}(f)} \\
& =&  \norm[p]{ D^{\circ} \circ T_s(f) + D^{\circ}(f)}  
\, = \,  \norm[p]{ T_s(f) + f } \\ 
& = & \norm[p]{ 2f - (f- T_s(f))} \\
& \geq & \norm[p]{2f} - \norm[p]{T_s(f)-f} \\
& \geq  & 2-2\eta_p(\epsilon).
\end{array}
\end{displaymath}
Since 
$ \norm[p]{S_s^{\circ}(D^{\circ}(f))}=\norm[p]{ D^{\circ}(f)}=1 $,
we infer 
$ \norm[p]{S_s^{\circ}(D^{\circ}(f)) - D^{\circ}(f)} \leq \epsilon $.
\qed
Now we may continue with the proof of Theorem~\ref{thm:dilthm}. 
\newline
{\bf Proof of Theorem~\ref{thm:dilthm}:}
The above lemma asserts that 
the range of $Q$ includes 
$ D^{\circ}(\Lp{\Omega}{\mu}) $. 
Therefore, denoting
$Y:=\Lp{\Omega^{\circ}}{\mu^{\circ}} $,
\begin{displaymath}
\begin{array}{rcl}
Q\circ D^{\circ} \circ T_s & = & Q \circ P^{\circ}_{|Y} \circ S^{\circ}_s 
\circ Q \circ D^{\circ},
\quad \mbox{ $ s \in \mathbb{Q^+} $}.
\end{array}
\end{displaymath}
The range of 
$ Q $ 
is a sub-lattice, as follows from   Lemma 6, Chap~6 \S17, of Lacey's book
\cite{Lacey_74},
and it is closed, since 
$ Q $ 
is a contractive projection. 
As an abstract 
$ L^p $-space
it is,
by the Bohnenblust-Nakano theorem, see e.g.\ Chap~5 \S15, Theorem 3 of
\cite{Lacey_74},  isometrically and order isomorphic to
$ \Lp{\tilde{\Omega}}{\tilde{\mu}} $,
for some measure space
$ (\tilde{\Omega}, \tilde{\mathfrak{A}}, \tilde{\mu}) $,
by a linear map
$ \Phi $,
say. 
Define
\begin{displaymath}
\begin{array}{rcl}
D & = & \Phi^{-1} \circ D^{\circ}, \\
P & = & \Phi^{-1} \circ P^{\circ} \, _{|Y} \circ \Phi, \\
S_s &=& \Phi^{-1} \circ Q \circ S^{\circ}_s \, _{|Y} \circ \Phi, \quad  s \in \mathbb{Q}.\\
\end{array}
\end{displaymath}
By continuity, the representation 
\index{$S_.~$}%
$ S_. $
can be extended to  a continuous representation of
$ \mathbb{R}  $,
still acting on 
$ \Lp{\tilde{\Omega}}{\tilde{\mu}} $.
By abuse of notation this extension will still be denoted 
$ S_. $. \par
For all 
$ f \in \Lp{\Omega}{\mu} $ 
we obtain
\begin{displaymath}
D \circ T_t (f) = P \circ S_t \circ D (f), \quad  t \in \mathbb{R} \, ,
\end{displaymath}
since both sides are continuous functions of 
$ t \in \mathbb{R} $ 
and the above equality is valid for the dense subgroup 
$ \mathbb{Q} $. 
\qed
\par
In the last part of this section we indicate the changes necessary,
for proving the analogue of Theorem~\ref{thm:dilthm} for a semigroup
$ \{ T_t : t \geq 0 \} $
of sub-positive contractions. 
In this case it can be shown, by the same reasoning as above, that 
$ Q $ 
is a contractive projection. From {\rm (iii)} of our proposition it 
can be seen to be a sub-positive contraction, even.
Anyway, in the case 
$ 1 < p < \infty $
and
$ \quad p \not= 2 $,
the structural description of the range of 
\index{projection!contractive}%
contractive projections on 
$ \mathcal{L}^{p} $-spaces, cf.\ e.g.\ chap.\ 6, \S 17, Theorem 3 \cite{Lacey_74},
guarantees that, 
for a direct sum 
$ U : \Lp{\Omega^{\circ}}{\mu^{\circ}} \rightarrow \Lp{\Omega^{\circ}}{\mu^{\circ}} $
of unitary multiplication operators,
$ UY $
is isometrically and order isomorphic to 
$ \Lp{\tilde{\Omega}}{\tilde{\mu}} $,
for some measure space
$ (\tilde{\Omega}, \tilde{\mathfrak{A}}, \tilde{\mu} ) $,
by an isomorphism which we  again call
$ \Phi $.
We note that 
$ U_{|Y} $ 
is invertible and
acts as the identity on the range of
$ D $,
since 
$ Q $ 
does, as follows from Lemma~\ref{lem:contlem}.
All we have to do is a further conjugation of
$  P^{\circ} \, _{|Y} $ 
and of
$ \{ Q \circ S^{\circ}_s \, _{|Y} : s \in \mathbb{Q} \} $
with
$ U_{|Y} $ 
when defining
$ P $ 
and
$ \{ S_s : s \in \mathbb{Q} \} $
before we conjugate the respective results with 
$ \Phi $. 
\par 
If 
$ p = 2 $,
then the closed subspace 
$ Y $ 
is isomorphic to
$ l^2(I) $
for some set 
$ I $.
In this case no assertion on (sub)positivity properties of 
the involved operators can be made 
and we simply transport the group 
$ \{ S_s : s \in \mathbb{Q} \} $
by means of this isomorphism. \par
The completion to a representation of
$\mathbb{R}$
and the last conclusion on the strong continuity of this representation
can be done exactly as before.
\chapter[Transference]{Transference}
\section[Transference for linear Maps]{Transference for linear Maps}
\label{sec:lintran}
Our theme is now to obtain reasonable norm estimates for operators
which can be defined by strongly convergent integrals
\begin{eqnarray*}
\int\limits_{0}^{\infty}{k(t)T_t} \; dt, 
\hspace{10mm} \mbox{$ k \in L^{ \rm 1} \it (\mathbb{R}_{+} , \lambda) $}.
\end{eqnarray*}
Here $ (T_{t})_{t \geq 0} $ is still a 
\index{semigroup!strongly continuous}%
\index{semigroup!one-parameter}%
one-parameter semigroup, strongly continuous,
of 
\index{semigroup!sub-positive contractions}%
\index{contraction!sub-positive}%
sub-positive contractions, acting on 
$ E \, = \, \Lp{\Omega}{\mu} $.
We shall be merely interested in the case that 
$1 < p < \infty$, 
though some results are valid for $p=1$ too.
Often in this case only slight modifications of the proofs are necessary, 
but for semigroups on $\mathcal{L}^1$-spaces we did not prove an analogue 
of Theorem~\ref{thm:dilthm}. Hence we can not take advantage of reducing
the problems given for one-parameter semigroups to problems
for one-parameter groups.
\par
We display the dilation given by that theorem, respectively its 
sub-positive version alluded to in the last part of section~\ref{sec:condil}
in a commutative diagram (see Figure~\ref{fig:dil4}).
Then it is clear that for all 
$ k \in L^{1}(\mathbb{R_{+}}, \lambda  )$
\begin{figure}
\unitlength1cm
\begin{picture}(13,5)
\put(0.5,4){$ E $}
\put(9.7,4){$ E  $ }
\put(0.5,1){$\tilde{E}$}
\put(9.7,1){$\tilde{E}$}
\put(6.7,1){$\tilde{E}$}
\put(1.5,4.1){\vector(1,0){7.5}}
\put(1.5,1.1){\vector(1,0){4.2}}
\put(7.3,1.1){\vector(1,0){2}}
\put(0.7,3.6){\vector(0,-1){2}}
\put(9.8,3.6){\vector(0,-1){2}}
\put(5,4.4){ $ T_t $ }
\put(3.7,1.4){ $ S_{t} $ }
\put(1,2.5){ $ D $}
\put(10,2.5){ $ D $}
\put(7.9,1.4){ $ P $ }
\end{picture}
\caption{Dilation for Reals}
\label{fig:dil4}
\end{figure}
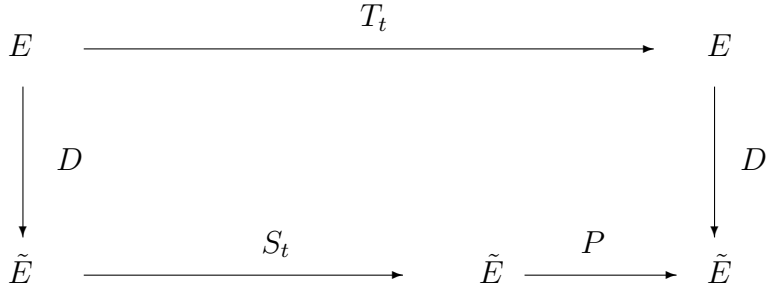
\begin{eqnarray*}
\norm{\int\limits_{0}^{\infty}{k(t)T_t} \; dt} & \, \leq \, &
\norm{\int\limits_{0}^{\infty}{k(t)S_t} \; dt},
\end{eqnarray*}
and we only need to estimate the latter.
\par
For this let for $k \in L^{1}(\mathbb{R}, \lambda)$
denote
$ \lambda_{p}(k) \, : \, \Lp{\mathbb{R}}{\lambda}\, \rightarrow  \,
\Lp{\mathbb{R}}{\lambda} $ 
the 
\index{convolution!on ${\mathbb{R}}$}%
convolution with $ k $:
\begin{eqnarray*}
\lambda_{p}(k)f\ (x) & \, = \, & \int_{\mathbb{R}} {k(y)f(x-y)} \, dy \, ,
 \quad f \in \Lp{\mathbb{R}}{\lambda}, \, x \in \mathbb{R}.
\end{eqnarray*}
Further we denote by
\begin{displaymath}
\norm[p,p]{k} \, = \, \sup \left \{ \norm[p]{\lambda_{p}(k)f} \, : 
\, \norm[p]{f} \,    = \, 1 \right \}
\end{displaymath}
its 
\index{$\norm[p,p]{k}$}%
\index{norm!operator}%
\index{norm!operator!$\norm[p,p]{k}$}%
operator norm on 
$\Lp{\mathbb{R}}{\lambda}$. 

\begin{thm}
\label{thm:mapnorm}%
\index{theorem!transference}%
\index{theorem!transference!$\mathbb{R}$}%
Assume $ 1<p<\infty$
and let
$E \, = \, \Lp{\Omega}{\mu}$ 
be  an $\mathcal{L}^{p}$-space. 
If 
$ \pi \, : \, \mathbb{R} \rightarrow B(E) $
is a uniformly bounded strongly continuous representation of $\mathbb{R}$
on E,
then, with 
$C=\sup \left \{ \norm{\pi(x)} \, : \, x \in \mathbb{R} \right \}$,
there holds for all $k \in L^{1}(\mathbb{R} , \lambda)$:
\begin{displaymath}
\norm{\int_{\mathbb{R}}k(x) \pi(x) \ dx } \, \leq \, C^2 \norm[p,p]{k}.
\end{displaymath}
\end{thm}
\begin{rem} \rm%
\label{rem:l1norm}
We do not prove the 
above theorem here since it is a special case of a result which
we shall establish  
(c.f.\ Theorem~\ref{thm:cbrepn})
soon.
But we would like to remark that for $p=1$ a stronger
estimate
\begin{eqnarray*}
\norm{\int_{\mathbb{R}}k(x) \pi(x) \ dx } \, \leq \, C \norm[1,1]{k}
\end{eqnarray*}
holds true for all $k \in L^1(\mathbb{R},\lambda)$,
simply because the norm $\norm[1,1]{.}$ coincides with the $L^{1}$-norm.
\end{rem}
\par 
The following corollary is due to Coifman and Weiss \cite{CoifWeiss_77}.
Note that the case $p=1$ holds true for the reasons mentioned in the last
remark.
\begin{cor}[Coifman-Weiss]
\label{cor:maptransfer}%
\index{theorem!transference!semigroup}%
\index{theorem!transference!Coifman~and~Weiss}%
Let $ (T_{t})_{{t \geq 0}} $ be a strongly continuous se\-mi\-group
of (sub)\-po\-si\-tive contractions acting on 
$ E \in \mathcal{L}^{p} $, where $1\leq p <\infty$.
Then for all $ k \in  L^{1}(\mathbb{R}_{+},\lambda)$:
\begin{displaymath}
\norm{\int\limits_{0}^{\infty}{k(t)T_t} \; dt} \, \leq
\, \norm[p,p]{k}.
\end{displaymath}
\end{cor}
\par
\section[Complete Boundedness of Transference]{Complete Boundedness of %
Transference}
\label{sec:transcpbd}
For the convenience of the reader we recall the notion of 
\index{bounded!p-complete}%
p-complete boundedness from \cite{Fendler_87} and \cite{Pisier_90}.
To this end we have to introduce some notation.
\par
Assume 
$1\leq p<\infty$ and let $E$ be a Banach space,
then for 
$ n \in \mathbb{N} $ 
and an  
$ n \times n $ 
\index{matrix!of operators}%
\index{operator!matrix}%
matrix 
\index{$\matr{m}{}{n}~\in~M_n~\otimes~B(E)$}%
$ \matr{m}{}{n} \in M_n \otimes  B(E) $ 
of operators $ m_{i,j} \in B(E) $
denote
\index{$\norm[(n)]{.}$}%
\begin{eqnarray}
\label{eq:p-cb}
\norm[(n)]{ \matr{m}{}{n}} &  = & %
\sup  \left( \sum_{i=1}^{n}\norm{ %
\sum_{j=1}^n m_{i,j} (g_j)}^{p} \right)^{\frac{1}{p}},
\end{eqnarray}
the supremum being taken over all %
$g_1, \ldots ,  g_n  \in  E, \; \mbox{with }$
\[   \left( \sum_{j=1}^{n}\norm{g_j}{}^{p} \right)^{\frac{1}{p}} 
\leq  1.
\]
\par
The reader may note that the above norm arises from considering the matrix of
operators
\begin{displaymath}
\matr{m}{}{n}
\end{displaymath}
as acing on the Banach space
$\lp{n}{E}$.
If 
$E \, = \, \Lp{\Omega}{\mu}$ 
is a space of class 
$\mathcal{L}^{p}$
then 
$\lp{n}{E}$
is clearly again of the same class. Therefore the notion we define next seems
to be most useful if the Banach spaces involved are $\mathcal{L}^{p}$-spaces or belong
to some related class of Banach spaces, e.g.\ 
subspaces or 
\index{$\mathcal{L}^{p}$-space!subspaces}%
subspaces of 
quotient spaces of 
\index{$\mathcal{L}^{p}$-space!quotient spaces}%
$\mathcal{L}^{p}$-spaces.
\par
\begin{defi}
\label{defi:pcbnorm}
Assume that $E,\; F$ are Banach spaces and that $1\leq p < \infty$.
If 
$ S \subset B(E) $ 
is a subspace, then we call a linear operator
$ u : S \rightarrow B(F) $
\index{bounded!p-complete}%
p-completely bounded if there exists a finite constant 
$ C>0 $ 
such that for all 
$ n \in \mathbb{N}$ 
and for all
$ \matr{m}{}{n} \in M_n \otimes  S $ 
\begin{displaymath}
\norm[(n)]{ \matr{m}{u}{n} } \leq C \norm[(n)]{ \matr{m}{}{n}}.
\end{displaymath}
We denote 
\index{$\cbnorm[p]{u}$}%
$ \cbnorm[p]{u} $
the least such constant.
\end{defi}
\par 
Since the results presented in this 
section, and their proofs, are almost verbatim the same in the more general 
situation when one is concerned with an 
\index{group!locally compact}%
\index{group!amenable}%
amenable locally compact group and a 
\index{Haar measure}%
\index{Haar measure!left}%
left Haar measure on it, instead of the 
locally compact abelian (hence amenable) additive group
$ \mathbb{R}$
with the (translation invariant) Lebesgue measure,
we chose this generality and we let denote
$ G $
a locally compact amenable group endowed with a Haar measure 
\index{Haar measure!$~\lambda~$}%
$ \lambda $ 
on it. 
\par
For
$ k \in L^{1}(G, \lambda) $ 
and
$ g \in \Lp{G}{\lambda} $
\index{convolution product}%
\index{$k\star g$}%
 the convolution product   
\[ k \star g (x) = \int_{G}{k(y)g(y^{-1}x)} \, d \lambda (y), \qquad
x \in G 
\] 
is defined. 
Further, 
\index{$\lambda_{p}:k\mapsto~(g\mapsto~k\star~g)$}%
$ \lambda_{p} : k \mapsto (g \mapsto k \star g ) $ 
is injective from 
\index{$L^{1}(G,\lambda)$}%
$ L^{1}(G,\lambda) $ 
into 
$ B(\Lp{G}{\lambda}) $ 
and we may consider its range
$ \lambda_{p}(L^1(G,\lambda)) $ 
as a normed subspace of 
$ B(\Lp{G}{\lambda})$.
\par
A 
\index{representation!continuous}%
continuous 
\index{representation}%
representation of
$ G $
on
$ \Lp{\Omega}{\mu} $
is, by definition, a group homomorphism 
$ \pi $
mapping 
$ G $ 
into the invertible operators on
$ \Lp{\Omega}{\mu} $
which is continuous when
those are endowed with the strong operator topology.
If, furthermore, 
$ \pi $
is uniformly bounded,
i.e.
$ \sup_{x \in G} \norm{ \pi (x) }{} < \infty $,
then we can consider its extension by integration,
$ \lambda_{\pi} : L^1(G,\lambda) \rightarrow B(\Lp{\Omega}{\mu}) $,
defined by
\begin{displaymath}
\lambda_{\pi}(k)f = \int_G k(x) \pi (x)f \; d\lambda (x), \quad
f \in \Lp{\Omega}{\mu}, \; k \in L^1(G,\lambda).
\end{displaymath}
We remark that
$ \lambda_{\pi}$
is an algebra homomorphism for the convolution structure on
$L^1(G,\lambda)$.
\begin{thm}
\label{thm:cbrepn}
Assume $1 \leq p <\infty$. 
Let 
$\pi :  G \rightarrow B(\Lp{\Omega}{\mu}) $
be a continuous uniformly bounded representation of
$ G  $.
Then
\index{$\lambda_{\pi}$}%
\index{representation!$\lambda_{\pi}$}%
\[
 \lambda_{\pi} : \lambda_p(L^1(G,\lambda)) \rightarrow B(\Lp{\Omega}{\mu}) 
\]
is a p-completely bounded algebra homomorphism with norm
\[
\cbnorm[p]{ \lambda_{\pi }} \leq \sup_{x \in G } \norm{ \pi (x)}{}^2 .
\]
\end{thm}
{\bf Proof:} 
We apply the amenability of the group 
$ G $
in a manner closely related to a F{\o}lner-Leptin condition.This seems to be
due to 
\index{Herz}%
C.~Herz \cite{Herz_73}, compare also \cite{CoifWeiss_77}. 
\par
For $ 1<p<\infty $ let $ q $ be defined by 
$ \frac{1}{p} + \frac{1}{q} = 1 $.
Then 
for 
\mbox{$ \alpha \in \Lp{G}{\lambda} $} 
and 
\mbox{$ \beta \in \Lq{G}{\lambda} $}
the convolution
\index{$\beta \star \alpha^{\vee} $}%
\mbox{$ \beta \star \alpha^{\vee} $}, 
where
\index{$\alpha^{\vee}$}%
\mbox{$ \alpha^{\vee}(x) := \alpha(x^{-1}), \quad x \in G $},
coincides 
$ \lambda $ 
almost everywhere with a continuous function
vanishing at infinity, and
$ \norm[{\infty}]{\beta \star \alpha^{\vee}} \leq 
\norm[p]{\alpha} \norm[q]{\beta} $. 
\par
This holds true for $p=1$, in which case $q=\infty$, given that 
$\beta$ has support of finite Haar measure.
\par
Since 
$ G $ 
is amenable, there exist nets 
$ (\alpha_{\tau})_{\tau \in \Delta} \subset \Lp{G}{\lambda} $ \\
and
\mbox{$ (\beta_{\tau})_{\tau \in \Delta} \subset \Lq{G}{\lambda} $}
with
\begin{displaymath}
\sup_{\tau \in \Delta} \norm[p]{ \alpha_{\tau}} \leq 1 ,\
\sup_{\tau \in \Delta} \norm[q]{ \beta_{\tau}} \leq 1
\end{displaymath}
such that 
\begin{displaymath}
\lim_{\tau \in \Delta} \beta_{\tau} \star \alpha^{\vee}_{\tau} = 1 
\end{displaymath}
uniformly on compact sets. 
\pagebreak[4] 
\par
Hence, whenever 
$ k \in L^1(G,\lambda), \; f \in \Lp{\Omega}{\mu} $
and 
$ g \in \Lq{\Omega}{\mu} $
are given, then there holds true:
\begin{eqnarray*}
\lefteqn{ \int_{\Omega} \lambda_{\pi}(k) f(\omega) g(\omega)
\; d\mu(\omega)  } \\ 
&  & = \lim_{\tau \in \Delta} \int_{\Omega} \int_G 
\beta_{\tau} \star \alpha_{\tau}^{\vee}(x) k(x) \pi (x) f(\omega) g(\omega)
\; d\lambda(x) d\mu(\omega)  \\
& & = \lim_{\tau \in \Delta} \int_{\Omega}
\lambda_{\pi}((\beta_{\tau} \star \alpha_{\tau}^{\vee}) \cdot k) 
f(\omega) g(\omega) \; d\mu(\omega) ,
\end{eqnarray*}
where $ \cdot $ denotes the pointwise product of functions.
This is an abuse of the dominated convergence theorem, since 
$ \Delta $ 
might be uncountable. But here we are concerned with one, later with
finitely many, integrable functions on 
$ G $.
They all vanish 
$ \lambda $ 
almost everywhere outside a 
$ \sigma $-compact 
subgroup for which one can arrange 
$ \Delta $
to be a sequence. Similar arguments justify the use of Fubini's theorem 
in the arguments given below. 
\par
If we denote by 
\index{$\pi^t~:G~\rightarrow~B(\Lq{\Omega}{\mu})$}%
$ \pi^t : G \rightarrow B(\Lq{\Omega}{\mu}) $
the representation 
\index{representation!adjoint}%
adjoint to $ \pi $, given by
\begin{eqnarray*}
\lefteqn{ \int_{\Omega}f(\omega)(\pi^t(x)g)(\omega) \; d\mu(\omega) }  \\
& & = \int_{\Omega}(\pi(x^{-1})f)(\omega)g(\omega) \; d\mu(\omega), \quad x \in G,
\end{eqnarray*}
then
\begin{eqnarray*}
\lefteqn{ \int_{\Omega} 
\lambda_{\pi}(\beta \star \alpha^{\vee} \cdot k) f(\omega) g(\omega)
\; d\mu(\omega)  } \\ 
& & = \int_{\Omega} \int_G \int_G k(yx) F^{\omega}(x^{-1})
G^{\omega}(y) \; d\lambda(x) d\lambda(y) d\mu(\omega),
\end{eqnarray*}
where 
\begin{eqnarray*}
F^{\omega}(x) & = & \alpha(x)\pi(x^{-1})f(\omega), 
\quad x \in G,\quad \omega \in \Omega \\
G^{\omega}(y) & = & \beta(y)\pi^t(y^{-1})g(\omega), 
\quad y \in G,\quad \omega \in \Omega. 
\end{eqnarray*}
\par
Finally for $ f_1,\dots,f_n \in \Lp{\Omega}{\mu} , \quad
g_1,\dots,g_n \in \Lq{\Omega}{\mu} $
and  \\
\mbox{$ \matr{k}{}{n} \in M_n \otimes \lambda _p(L^1(G,\lambda)) $}
we compute, with
$ F_{\tau ,j}^{\omega} $  
and
$ G_{\tau,i}^{\omega} $  
defined in analogy to the above functions
$ F^{\omega} $
and 
$ G^{\omega} $: 
\begin{eqnarray*}
\lefteqn{ | \sum_{i,j=1}^n \int_{\Omega} \lambda_{\pi}(k_{i,j}) (x) 
f_j(\omega) g_i(\omega) \; d\mu(\omega) | } \\
& & =  \lim_{\tau \in \Delta} | \sum_{i,j=1}^n \int_{\Omega}
\lambda_{\pi}((\beta_{\tau} \star \alpha_{\tau}^{\vee}) \cdot k_{i,j}) 
f_j(\omega) g_i(\omega) \; d\mu(\omega) | \\
& & =  \lim_{\tau \in \Delta} | \int_{\Omega} \sum_{i,j=1}^n
\int_G \int_G k_{i,j}(yx) F_{\tau ,j}^{\omega}(x^{-1}) G_{\tau,i}^{\omega}(y)
\; d\lambda(x) d\lambda(y) d\mu(\omega) |. 
\end{eqnarray*}
We dominate this by
\begin{eqnarray*}
\lefteqn{\lim_{\tau \in \Delta} \int_{\Omega} \left\{
\norm[(n)]{\matr{k}{\lambda_p}{n}}
\left( \sum_{j=1}^n \int_G |F_{\tau ,j}^{\omega}(x) |^p \; d\lambda(x)
\right) ^{\frac{1}{p}} \right.
\!\! \cdot} \\
& &  \quad \quad \left. \left( \sum_{i=1}^n \int_G 
|G_{\tau,i}^{\omega}(y) |^q
d\lambda(y) \right) ^{\frac{1}{q}} \right \} \; d\mu(\omega)  \\
& & \leq \norm[(n)]{\matr{k}{\lambda_p}{n}} \lim_{\tau \in \Delta} \left\{ \!\!
\left( \sum_{j=1}^n \int_G \!\! 
|\alpha_{\tau}(x)|^p \! \int_{\Omega} \!\! | \pi(x^{-1}) f_j(\omega) |^p 
\, d\mu(\omega) d\lambda(x) \right) ^{\!\! \frac{1}{p}} \right.
\!\!\!\! \cdot \\
& & \quad \quad \left. \left (\sum_{i=1}^n \int_G \!\! |\beta_{\tau}(y)|^q \!\!
\int_{\Omega} \!\! | \pi^t(y^{-1}) g_i(\omega) |^q 
\; d\mu(\omega) d\lambda(y) \right) ^{ \!\! \frac{1}{q}} \right\} \\
& & \leq \norm[(n)]{\matr{k}{\lambda_p}{n}}  \lim_{\tau \in \Delta}
\left(\sup_{x \in G } \norm{ \pi (x^{-1})}{} \norm[p]{\alpha_{\tau}}
(\sum_{j=1}^n \norm[p]{f_j}^p)^{\frac{1}{p}} \cdot \right. \\
& & \quad \quad \left. \sup_{y \in G } \norm{ \pi^{t} (y^{-1})}{}
\norm[q]{\beta_{\tau}}(\sum_{i=1}^n \norm[q]{g_i}^q)^{\frac{1}{q}} \right) \\
& & \leq \sup_{x \in G } \norm{ \pi (x)}{}^2 \norm[n]{\matr{k}{\lambda_p}{n}}
(\sum_{j=1}^n \norm[p]{f_j}^p)^{\frac{1}{p}}
(\sum_{i=1}^n \norm[q]{g_i}^q)^{\frac{1}{q}}. \\
\end{eqnarray*}
For the case $p=1$, i.e.\ $q=\infty$,
this estimate has to be modified accordingly. 
\qed
Similarly as Theorem~\ref{thm:mapnorm} combined with the Dilation Theorem~\ref{thm:dilthm}~had a nice conclusion which we stated as Corollary~\ref{cor:maptransfer}, we here
obtain a strengthening of that corollary for $p$ in the range $1<p<\infty$.
\par
\begin{cor}
\label{cor:estimcor}%
\index{transference!p-completely bounded}%
\index{theorem!transference!p-completely bounded}%
Let 
$ 1<p<\infty $
and 
$ \{ T_t : t \geq 0 \}  $ 
be a strongly continuous semigroup of positive or subpositive contractions
acting on
$\Lp{X}{\mu}$. 
If for some 
\mbox{ $ n \in \mathbb{N} $ }
and some matrix 
\index{$\matr{k}{}{n} \in M_n \otimes  L^{1}(\mathbb{R},\lambda) $}%
\mbox{$ \matr{k}{}{n} \in M_n \otimes  L^{1}(\mathbb{R},\lambda) $},
whose entries 
\mbox{$ k_{i,j} \in L^{1}(\mathbb{R}, \lambda) $}
have support in 
$ \mathbb{R_{+}} $, 
there exists a constant 
\mbox{$ C \geq 0 $} 
such that for any 
$ n $ elements
\mbox{$ g_1, \dots ,g_n \in \Lp{\mathbb{R}}{\lambda} $}
\begin{displaymath}
\left( \sum_{i=1}^{n}\int_{\mathbb{R}}\abs[ \sum_{j=1}^{n} \lambda_{p}(k_{i,j}) g_j (y) ]^p \;
dy\right)^{\frac{1}{p}}
\leq C\left( \sum_{j=1}^n\int_{\mathbb{R}}\abs[ g_j(y)]^p \; dy\right)^{\frac{1}{p}}, 
\end{displaymath}
then for all 
$  f_1, \dots , f_n \in \Lp{X}{\mu} $
\begin{displaymath}
\left( \sum_{i=1}^n \norm[p]{ \sum_{j=1}^n {\int_{0}^{\infty}  k_{i,j}(t) \
T_{t}f_j \, dt}} ^p \right)^{\frac{1}{p}}
 \leq C \left(\sum_{j=1}^n \norm[p]{f_j}^p \right)^{\frac{1}{p}}.
\end{displaymath} 
\end{cor}
{\bf Proof:} 
Let
$ \{ T_t : t \geq 0 \} $
be the semigroup under consideration.
We have to prove that for a matrix
$ \matr{k}{}{n} \in M_n \otimes L^1(\mathbb{R^+},\lambda) $
\begin{eqnarray*}
\norm[(n)]{ \left(\int_{0}^{\infty}{k_{i,j}(t) \ T_{t} \, 
                                         dt}\right)_{i,j=1}^{n}} 
& \, \leq \, & \norm[(n)]{ \matr{k}{ \lambda_p }{n}}. \\
\end{eqnarray*}
Let
\begin{eqnarray*}
D \circ T_t & = & P \circ S_t \circ D, \hspace{10mm}  t \geq 0,
\end{eqnarray*}
be a dilation according to Theorem~\ref{thm:dilthm}, respectively according
to the remarks in the last part of section~\ref{sec:condil}. 
Since 
$ \{ T_t : t \geq 0 \} $
and
$ \{ S_t : t \in \mathbb{R} \} $
are strongly continuous (semi)groups:
\begin{eqnarray*}
 D \circ  \int\limits_{0}^{\infty}{k(t)  T_t} \, dt 
& \, = \, & \int\limits_{0}^{\infty}{k(t) D \circ T_t} \, dt \\
& \, = \, & \int\limits_{0}^{\infty}{k(t) D \circ S_t \circ P } \, 
dt \\
& \, = \, &  D \circ \int\limits_{0}^{\infty}{k(t) S_t \; dt \circ P } \\
& \,  = \, & D \circ S (k) \circ P, 
\hspace{10mm} k \in L^{1}(\mathbb{R}^+ , \lambda). 
\end{eqnarray*}
Hence, by Theorem~\ref{thm:cbrepn}, applied to the continuous representation 
$ S $ 
of 
$ \mathbb{R} $ 
\begin{eqnarray*}
\norm[(n)]{  \left(\int_{0}^{\infty}{k_{i,j}(t) \ T_{t} \, dt}\right)_{i,j=1}^{n}} 
& = & 
\norm[(n)]{  \left(\int_{0}^{\infty}{D \circ k_{i,j}(t) \ T_{t}\, dt}\right)_{i,j=1}^{n}} \\
& = & \norm[(n)]{ \left({ D \circ S (k_{i,j}) \circ P}\right)_{i,j=1} ^{n}}   \\
& \leq & \norm[(n)]{ \matr{k}{ S }{n}} \\
& \leq & \norm[(n)]{ \matr{k}{ \lambda_p }{n}}   \\
\end{eqnarray*}
which completes the proof of Corollary~\ref{cor:estimcor}.\qed
\section{Transference of Square Functions}
\label{sec:square}
In this section we shall apply our
knowledge about $p$-completely bounded maps 
to the transference of certain square functions.
\par
We are more interested to give some examples 
than to obtain sharp results. In fact, 
as far as representations of amenable groups
or representations of continuous one-parameter
semigroups on $\mathcal{L}^1$-spaces are concerned
our results are not optimal.
There are two reasons for this.
\begin{description}
\item[(i)~] For continuous  one-parameter
semigroups on $\mathcal{L}^1$-spaces
we did not prove the existence of a
dilation to a strongly continuous one-parameter group.
Our constructions heavily involved reflexivity and hence
our proof of the subsequently stated result of
Coifman and Weiss~\cite{CoifWeiss_77}, our Corollary~{\ref{semigrsquare}},
does not cover the case $p=1$.
\item[(ii)] Using $p$-completely bounded maps to derive
square function estimates, we pass by inclusions
of finite dimensional Hilbert spaces as norm
closed subspaces of $\mathcal{L}^p$-spaces. 
More even, we use a uniform norm bound
on the canonical projection onto these subspaces
(c.f.\ Lemma~{\ref{khinprojlem}}).
\end{description} 
\par
For the case of representations of amenable groups
results for $p$ in the whole range
$p \in [1,\,  \infty)$, with a better constant
than ours,
have been obtained by
\index{Asmar}%
Asmar, 
\index{Berkson}%
Berkson and 
\index{Gillespie}%
Gillespie in~\cite{AsmBerGil_91a}.
\par
Let
$ r_{1},\, r_{2},\, r_{3},\, \ldots$
denote the 
\index{Rademacher}%
\index{functions!Rademacher}%
Rademacher functions, defined by:
\begin{eqnarray*} 
r_{1}(t)  &  \, = \, & \left\{ \begin{array}{ll}
                       1  \quad & \mbox{ if $0 \, \leq \, t \, \leq \,
                       \frac{1}{2}$} \\
                       -1 \quad & \mbox{ if $\frac{1}{2} \, < \, t \, 
                       \leq \, 1$},
                          \end{array}
                 \right. \\
\mbox{and for $ n \, \geq \, 1 \, $} :&{}&{}\\
r_{n+1}(t) & \, = \, &  \left\{ \begin{array}{ll}
                       r_{n}(2t)  \quad & \mbox{ if $0 \, \leq \, t \, \leq \,
                       \frac{1}{2}$} \\
                       r_{n}(2t-1) \quad & \mbox{ if $\frac{1}{2} \, < \, t \, 
                       \leq \, 1$}
                          \end{array}
                 \right..
\end{eqnarray*}
The
\index{Khinchin}%
\index{inequalities!Khinchin}%
 Khinchin inequalities then assert that 
for $1\leq p < \infty$
there exist constants
$ c_{p}$
and
$C_{p}$
such that for all sequences 
$(\alpha_{i})_{i \in \mathbb{N}} \, \subset \, \mathbb{C}:$
\begin{displaymath}
c_{p} \left( \int_{0}^{1}{\abs[\sum_{i=1}^{\infty} \alpha_{i} r_{i}(t)]^{p}}\, dt
\right)^{\frac{1}{p}}
\, \leq \, \left(\sum_{i=1}^{\infty} \abs[ \alpha_{i}]^{2}
\right)^{\frac{1}{2}}
\, \leq \, C_{p} \left( \int_{0}^{1}{\abs[\sum_{i=1}^{\infty} \alpha_{i} 
r_{i}(t)]^{p}}\, dt
\right)^{\frac{1}{p}}.
\end{displaymath}
\par
\begin{rem}
\label{rem:Khinineq} \rm
Trivially 
$c_p=1$ for $1\leq p \leq 2$ and $ C_p = 1 $ for $ 2 \leq p < \infty$.
The best possible constants have been computed, for $p=1$ by
Szarek~\cite{Szarek_76} and for the other cases by Haagerup~\cite{Haagerup_78}.\end{rem}
Let
$E = \Lp{\Omega}{\mu} $
be an $\mathcal{L}^{p}$-space $(1<p<\infty)$
and for
$k=(k_{1},\ldots,k_{n}) \in B(E)^n$
let
\index{$m^{k}~\in~M_{2^{n}}~\otimes~B(E)~$}%
$m^{k} \in M_{2^{n}} \otimes B(E)$
be defined by:
\begin{eqnarray}
\label{rowmatr}
m^{k}_{ij }& = & \left\{ \begin{array}{lll}
                  { \displaystyle \left(\frac{1}{2^{n}}\right) ^{\frac{1}{p}}
                               \sum_{l=1}^{n}r_{l}(\frac{i}{2^{n}})k_{l}} 
                               & \mbox{ if $ j=1 $} &  \\ 
                             0 
                               & \mbox{ if $ j \not= 1 $} & 
                                 \mbox{ for  $i,j \in \{1, \ldots , 2^{n}\} $.}
                       \end{array} \right.
\end{eqnarray}     
So that 
$m$
is a 
\index{operators!row-matrix}%
row-matrix of operators in
$M_{2^{n}} \otimes B(E)$
which space is normed according to~(\ref{eq:p-cb}).            
\begin{lem}
\label{lem:rowkhin}
With the above notations the following inequalities hold true:
\begin{displaymath}
c_{p}\norm[(2^{n})]{\, m^{k} \,} 
\, \leq \, 
\sup \left\{ \norm[E]{
  \left( \sum_{j=1}^{2^{n}} \abs[k_{j}(g)]^{2} \right )^{\frac{1}{2}}} \, : \,
  \norm[E]{g}\leq 1 \right\}
\, \leq \, C_{p}\norm[(2^{n})]{m^{k}}.
\end{displaymath}
\end{lem}
{\bf Proof:} We take 
$k = (k_1, \ldots ,k_n) \in B(E)^n$,
form $m \in M_{2^{n}} \otimes B(E)$ accordingly
(we suppress the super index $k$ here and in the following proof),
and for
$g_{1},\ldots, g_{2^{n}} \in E$
we compute:
\begin{eqnarray*}
c_{p}^p \sum_{i=1}^{2^{n}}\norm[E]{\sum_{j=1}^{2^{n}} m_{ij}(g_{j})}^p
& \, = \, &
c_{p}^p \sum_{i=1}^{2^{n}} \int_{\Omega} 
\abs[m_{i1}(g_{1})(\omega)]^{p}\, d\mu(\omega)\\
& \, = \, &
c_{p}^p \sum_{i=1}^{2^{n}}\int_{\Omega} \frac{1}{2^n} 
\abs[\sum_{l=1}^{n}r_{l}(\frac{i}{2^{n}})k_{l}(g_{1})(\omega)]^{p}\,
d\mu(\omega) \\
& \, = \, &
c_{p}^p \int_{\Omega} \int_{0}^{1}
\abs[\sum_{l=1}^{n}r_{l}(t)k_{l}(g_{1})(\omega)]^{p}\, dt \,
d\mu(\omega) \\
& \, \leq \, &
\int_{\Omega}
\left( \sum_{j=1}^{n} \abs[k_{j}(g_{1})(\omega)]^{2} \right )^{\frac{p}{2}}
\, d\mu(\omega) \\
& \, \leq \, &
C_{p}^p \int_{\Omega} \int_{0}^{1}
\abs[\sum_{l=1}^{n}r_{l}(t)k_{l}(g_{1})(\omega)]^{p}\, dt \,
d\mu(\omega) \\
& \, = \, & 
C_{p}^p \sum_{i=1}^{2^{n}} \norm[E]{\sum_{j=1}^{2^{n}} m_{ij}(g_{j})}^{p}
\end{eqnarray*}
Hence on one hand
\begin{eqnarray*}
c_p \norm[(2^n)]{m} & = & c_p \sup \left \{ \left( \sum_{i=1}^{2^{n}} 
\norm[E]{\sum_{j=1}^{2^{n}} m_{ij}(g_{j})}^{p} \right)^{\frac{1}{p}} \, : \,
\sum_{j=1}^{2^{n}} \norm[E]{g_{j}}^{p}\, \leq \, 1 \, \right \} \\
& \leq &{\sup \left\{ \norm[E]{
\left( \sum_{j=1}^{n} \abs[k_{j}(g_1)]^{2} \right )^{\frac{1}{2}}} \, : \,
\sum_{j=1}^{2^{n}} \norm[E]{g_{j}}^{p}\, \leq \, 1 \, \right \} }\\
& \leq & \sup \left\{ \norm[E]{
\left( \sum_{j=1}^{n} \abs[k_{j}(g_1)]^{2} \right )^{\frac{1}{2}}} \, : \,
\norm[E]{g_1}\leq 1 \right\}.
\end{eqnarray*}
And on the other hand, choosing 
$ g_1=g, \quad g_2 = \ldots =g_{2^n} = 0$
\begin{eqnarray*}
\norm[E]{
\left( \sum_{j=1}^{n} \abs[k_{j}(g)]^{2} \right )^{\frac{1}{2}}}
& \leq & C_p \norm[(2^n)]{m} \left( \sum_{j=1}^{2^n} \norm[E]{g_j}^p
  \right)^{\frac{1}{p}}\\
& \leq & C_p \norm[(2^n)]{m} \norm[E]{g}. \qquad \qquad \qquad \qquad \qquad \qed
\end{eqnarray*}

%
\begin{thm}
\label{thm:rowtrans}
\index{theorem!transference}%
\index{theorem!transference!p-completely bounded}%
If 
$E,F $
are $\mathcal{L}^{p}$-spaces,
$1\leq p < \infty$,
$S\subset B(E)$  is
a closed subspace and
$u\, : \, S \rightarrow \, B(F)$
a p-completely bounded linear map,
then for any sequence
$k_{1},k_{2},\ldots \, \in S :$
\begin{eqnarray*}
\lefteqn{\sup \left\{ \norm[F]{\left( \sum_{j=1}^{\infty} \abs[u(k_{j})(g)]^{2}
    \right )^{\frac{1}{2}}}  :  \norm[F]{g}\leq 1 \right\}} & & \\
& \leq &
\frac{C_{p}}{c_{p}}\cbnorm[p]{u}\sup \left\{ \norm[E]{\left( \sum_{j=1}^{\infty} \abs[k_{j}(f)]^{2}
    \right )^{\frac{1}{2}}}  :  \norm[E]{f}\leq1 \right\}.
\end{eqnarray*}
\end{thm}
{\bf Proof:}
It suffices to prove this for a finite sequence
$k_{1},\ldots ,k_{n}\, \in S$.
For, if $k_{1},\ldots \in S$
is infinite, then for any $g\in F$ the sequence of functions 
$(h_n)_{n \in \mathbb{N}}$,
\[ 
h_n = \left( \sum_{j=1}^{n} \abs[u(k_{j})(g)]^{2} %
\right )^{\frac{p}{2}} \in {L}^1,
\]
is non-decreasing and the monotone convergence theorem can be applied.
But then, by  Lemma~\ref{lem:rowkhin}, the above left hand side is dominated by
\begin{eqnarray*}
C_{p} \norm[(2^{n})]{\matr{n}{}{2^{n}}} & \leq & \cbnorm[p]{u} C_{p}
\norm[(2^{n})]{\matr{m}{}{2^{n}}} \\
& \leq &
\frac{C_{p}}{c_{p}}\cbnorm[p]{u}\sup \left\{ \norm[E]{\left( \sum_{j=1}^{2^{n}}
      \abs[k_{j}(f)]^{2} \right)^{\frac{1}{2}}} \, : \, \norm[E]{f} \leq 1
\right\},
\end{eqnarray*}
where
$m \, = \, \matr{m}{}{2^{n}}$
is constructed from the sequence
$k_{1},\ldots ,k_{n}\, \in S$
and
$n  = \matr{n}{}{2^{n}}  =  \matr{m}{u}{2^{n}}  = \id \otimes u \left(\matr{m}{}{2^{n}} \right)$
from 
$u(k_{1}),\ldots ,u(k_{n})  \in B(F)$
just as in~(\ref{rowmatr}) before the lemma.
\qed
As mentioned in the introduction to this section 
we are not able to derive from the above theorem 
and from Corollary~\ref{cor:estimcor} the case 
$p=1$ of a theorem of Coifman, Weiss 
(their Corollary 4.17. in~\cite{CoifWeiss_77}).
By our means we only can prove:
\begin{cor}
\label{semigrsquare}
\index{transference!square functions!semigroup}%
\index{theorem!transference!square functions}%
Let 
$(T_{t})_{t\geq0}$
be a strongly continuous semigroup of positive or sub-positive
contractions on 
$E=\Lp{\Omega}{\mu}$, where $1<p<\infty$.
If
$k_{1},k_{2}, \ldots \in L^{1}(\mathbb{R},\lambda)$
with
$\mbox{\rm supp}\ k_{i} \subset [0,\infty),\ i \in \mathbb{N}$
are such that for some 
$M>0$
for all
$f\in \Lp{\mathbb{R}}{\lambda}$
\[
\int_{\mathbb{R}} \left(\sum_{i=1}^{\infty}
 \abs[\lambda_{p}(k_{i})f(x)]^{2}\right)^{\frac{p}{2}} \, dx 
\, \leq \, 
M^{p}\int_{\mathbb{R}}\abs[f(x)]^{p}\, dx, 
\]
then for all
$f \in \Lp{\Omega}{\mu}$ 
we have
\[
\norm[E]{\left( \sum_{i=1}^{\infty}\abs[\int_{0}^{\infty}
k_{i}(t)T_{t}f \ dt]^{2} \right)^{\frac{1}{2}}}
\, \leq \,
M \frac{C_{p}}{c_{p}}\norm[E]{f}.
\]
\end{cor}
{\bf Proof:}
The Corollary~\ref{cor:estimcor} of the last section just states that the map
\[
u \, : \, \lambda_{p}(k)  \, \mapsto \,  \int_{0}^{\infty}k(t)T_{t} \, dt,
\]
defined on 
$S \ = \ \left\{ \lambda_{p}(k) \ : \ \mbox{supp}(k)\in [0,\infty), \
 k\in L^{1}(\mathbb{R},\lambda)  \right\} \ \subset  B(\Lp{\mathbb{R}}{\lambda})$,
is p-completely bounded from
$S$
to
$B(E)$.
Hence, by the above theorem, the assertion of the corollary  is immediate.
\qed
\begin{rem} \hspace*{\fill}~\\[-2em] \rm
\begin{description}
\item[(i)~]
The corollary has an obvious counterpart for uniformly bounded
strong\-ly continuous representations of amenable
locally compact groups
\linebreak[4] 
$\pi \, : \, G \rightarrow B(E)$
on $\mathcal{L}^{p}$-spaces $E$.
Here the p-complete boundedness of 
$\pi \, : \, \lambda_{p}(L^{1}(G))) \rightarrow B(E) $
follows directly from Theorem~\ref{thm:cbrepn}. We leave it to the reader to 
formulate and prove this (the case 
$p=1$ can be included).
 For a different, more direct approach to this we refer the reader
to the paper~\cite{AsmBerGil_91a} of 
\index{Asmar}%
Asmar, 
\index{Berkson}%
Berkson and 
\index{Gillespie}%
Gillespie.
\item[(ii)] 
In the proof of Lemma~\ref{lem:rowkhin} we used just the embeddings of
$l^{2}_{n}$
in $L^p([0,1])$ given by the Rademacher functions.
If, instead of Rademacher functions, appropriate 
\index{functions!Rademacher}%
\index{random variables!Gaussian}%
Gaussian random variables are used, then
one can avoid the appearance of the constant $\frac{C_{p}}{c_{p}}$
in the final result.
(This requires 
an additional approximation procedure, 
since a finite
subset of the considered random variables should be realizable on a finite
set). 
\end{description}
\end{rem}
Similar to the former we may consider the problem of transfering,
by means of a p-completely bounded
$u\, : \, S \rightarrow B(\Lp{\Omega'}{\mu'})$, 
\index{inequalities!square-function}%
square-function inequalities  of the following type:\\
Given a sequence
$k_{1},k_{2}, \ldots \in S \subset B(\Lp{\Omega}{\mu})$ 
such that for some $M$
\[
\norm[p]{\left(\sum_{i=1}^{\infty}\abs[k_{i}(f_{i})]^{2}\right)^{\frac{1}{2}}}
\, \leq \,
M \norm[p]{\left(\sum_{i=1}^{\infty}\abs[f_{i}]^{2}\right)^{\frac{1}{2}}}
\]
for all sequences
$f_{1},f_{2},\ldots \in \Lp{\Omega}\mu$,
we may ask whether a similar inequality holds true for the sequence 
$u(k_{1}),u(k_{2}), \ldots \in B(\Lp{\Omega'}{\mu'})$, may be with some
additional constant depending on $p$ only.
\par
Whereas in our last problem we discussed, for varying
$n \in \mathbb{N}$ and $k=(k_1, \ldots , k_n )$,
operators
\index{$K_{k}$}%
\begin{eqnarray*}
K_{k} \, : \, \Lp{\Omega}{\mu}
& \rightarrow &
\Lp{[0,1]}{\lambda;\Lp{\Omega}{\mu}},\\
K_{k} \, : \, f 
& \mapsto & 
\sum_{i=1}^{n}r_{i}(.)k_{i}(f), 
\end{eqnarray*}
we shall now consider
\index{$K_{k}'$}%
\begin{eqnarray*}
K_{k}' \, : \, \Lp{[0,1]}{\lambda; \Lp{\Omega}{\mu}}
& \rightarrow &
\Lp{[0,1]}{\lambda;\Lp{\Omega}{\mu}},\\
K_{k}' \, : \, F 
& \mapsto & 
\int_{0}^{1}\sum_{i=1}^{n}r_{i}(.)k_{i}(F(t))r_{i}(t)\,dt. 
\end{eqnarray*}
Thus,
for a sequence
$k \, = \, (k_{i})_{i \in \mathbb{N}} \subset B(\Lp{\Omega}{\mu})$ 
denote 
\index{$\norm[{[2]}]{k}$}%
\index{norm!$\norm[{[2]}]{k}$}%
$\norm[{[2]}]{k}$ 
the least constant 
$M$ such that for all $ f_{1},f_{2}, \ldots \in \Lp{\Omega}{\mu} : $
\[
\norm[p]{\left( \sum_{j=1}^{\infty}
      \abs[k_{j}(f_{j})]^{2} \right)^{\frac{1}{2}}}
\, \leq \, M 
\norm[p]{\left( \sum_{j=1}^{\infty}
      \abs[f_{j}]^{2} \right)^{\frac{1}{2}}} \,. 
\]
That is,
$\norm[{[2]}]{k}$ 
is the norm of the sequence
$k$
as an operator on 
$\Lp{\Omega}{\mu ; l^{2}}$.
\begin{prop}
\label{prop:khin}
There exist constants
$d_{p}$ and $D_{p}$ depending only on $p$, 
\linebreak[4]
$ 1<p<\infty $,
such that
\[
d_{p} \norm{K_{k}'} \, \leq \, \norm[{[2]}]{k} \, \leq \, D_{p} \norm{K_{k}'}
\]
for all sequences 
$k \, = \, (k_{i})_{i \in \mathbb{N}} \subset B(\Lp{\Omega}{\mu})$. 
\end{prop}
To prove the proposition we note first that for the finite sequence
\linebreak[4]
$l = (\id,\ldots ,\id)$
the operator
$K_{l}'\, = \, P_{n} \otimes \id$ 
is just the tensor product of the projection
$P_{n}\, : \, \Lp{[0,1]}{\lambda} \rightarrow  \Lp{[0,1]}{\lambda}$
onto the span of the first n 
Rademacher functions and of the identity
on $\Lp{\Omega}{\mu}$. 
Thus, for 
$F \in \Lp{[0,1]}{\lambda; \Lp{\Omega}{\mu}}$,
\[
\label{khinproj}
K_{l}'(F)(s)\, = \,P_{n} \otimes \id (F)(s)\, = \, 
\sum_{i=1}^{n}\int_{0}^{1}F(t)r_{i}(t)\, dt\, r_{i}(s),\quad s \in [0,1],
\]
and the 
\index{inequalities!Khinchin}%
Khinchin inequalities imply: 
\begin{lem}
\label{khinprojlem}
If $p \in (1,\infty)$, then for all 
$F \in \Lp{[0,1]}{\lambda; \Lp{\Omega}{\mu}}:$
\[
\norm{P_{n} \otimes \id (F)}\, \leq \,\max \{c_{p}^{-1}, c_{q}^{-1}\}\norm{F},
\]
where $q\in (1,\infty)$ is such that $\frac{1}{p}+\frac{1}{q}=1$.
\end{lem}
{\bf Proof:}
If $2 \leq p$, then
for $\mu$ almost all $\omega \in \Omega $ 
\begin{eqnarray*}
c_{p}^{p}\int_{0}^{1}\abs[P_{n} \otimes \id (F)(t,\omega)]^{p}\,dt
&\, \leq \,&
\left(\int_{0}^{1}\abs[P_{n} \otimes \id (F)(t,\omega)]^{2}\,dt 
\right)^{\frac{p}{2}}\\
&\, \leq \,& 
\left(\int_{0}^{1}\abs[F(t,\omega)]^{2}\,dt \right  )^{\frac{p}{2}}\\
&\, \leq \,&
\int_{0}^{1}\abs[F(t,\omega)]^{p}\,dt,
\end{eqnarray*}       
since  
$P_{n}$
is an orthogonal projection on the Hilbert space
$L^{2}([0,1],\lambda)$.
Integrating this with respect to $\mu$ proves the lemma 
in this case. 
\par
For $p$ in the range $1<p <2$ the assertion is then established by duality:\\
Let $q \in [2,\infty)$ be such that $\frac{1}{p} + \frac{1}{q} = 1$
then for $G \in %
\Lq{[0,1]}{\lambda; \Lq{\Omega}{\mu}}$:
\begin{eqnarray*}
\abs[\int_{\Omega}\int_{0}^{1} \!\! P_{n} \otimes \id (F)(t,\omega)G(t,\omega)\,
dt\, d\omega ] \!\!%
& = & \!\! \abs[\int_{\Omega}\int_{0}^{1} \!\! F(t,\omega) %
P_{n} \otimes \id (G)(t,\omega)dt\, d\omega ]\\
& \leq & \norm[p]{F} %
\norm[q]{P_{n} \otimes \id (G)} \\
& \leq & c_{q}^{-1} \norm[p]{F} \norm[q]{G}.
\end{eqnarray*}
By the converse to H\"{o}lder's inequality we obtain the assertion.
\hspace{\fill}
\qed
{\bf Proof of Proposition~\ref{prop:khin}:}
As in the proof of Theorem~\ref{thm:rowtrans}  we may assume that
$k\, = \, (k_{1}, \ldots , k_{n}) \subset B(E)$
is finite.
Let 
$f_{1}, \ldots , f_{n}\, \in E=\Lp{\Omega}{\mu}$
be given.
Then, by the orthogonality relations of the Rademacher functions,
\begin{eqnarray*}
K_{k}' \left(\sum_{j=1}^{n}r_{j}f_{j} \right)(s,\omega)
& \, = \, &
\int_{0}^{1}\sum_{i=1}^{n}\sum_{j=1}^{n}
r_{i}(s)k_{i}(f_{j})(\omega)r_{i}(t)r_{j}(t) \, dt \\
& \, = \, &
\sum_{i=1}^{n}r_{i}(s)k_{i}(f_{i}(\omega)).
\end{eqnarray*}
Hence we obtain the right hand side inequality of the proposition, with 
$D_{p}= \frac{C_p}{c_{p}}$,
from
\begin{eqnarray*}
\int_{\Omega}\left( 
\sum_{i=1}^{n}\abs[k_{i}(f_{i})(\omega)]^{2}\right)^{\frac{p}{2}} \, d\mu(\omega)
& \, \leq \, &
C_{p}^{p}\int_{\Omega} \int_{0}^{1} 
\abs[\sum_{i=1}^{n}r_{i}(t)k_{i}(f_{i})(\omega)]^{p} \, dt d\mu(\omega)\\
& \, \leq \, &
C_{p}^{p}\norm[\Lp{[0,1]}{\lambda;E}]
{K_{k}'\left(\sum_{j=1}^{n}r_{j}f_{j}\right)}^{p}\\
& \, \leq \, &
C_{p}^{p} \norm{K_{k}'}^{p} \norm[\Lp{[0,1]}{\lambda;E}]
{\left(\sum_{j=1}^{n}r_{j}f_{j}\right)}^{p}\\
& \, \leq \, & \left( \frac{C_p}{c_{p}} \right) ^p %
 \norm{K_{k}'}^{p} \norm[E]{\sum_{i=1}^{n} %
\left(\abs[f_{i}]^{2}\right)^{\frac{1}{2}}}^p.
\end{eqnarray*}
\par
To deduce the other estimate we apply the Lemma~\ref{khinprojlem}.\\ 
Given
$F \in \Lp{[0,1]}{\lambda;E}$
define 
$f_{i} \, = \, \int_{0}^{1}F(t)r_{i}(t) \, dt \in \Lp{\Omega}{\mu}$,
so that  
\[
P_{n} \otimes \id (F)(s,\omega) \, = \, \sum_{i=1}^n r_i(s) \otimes f_i(\omega)\,. 
\]
Then,
\begin{eqnarray*}
\norm[\Lp{[0,1]}{\lambda;E}]{K_{k}' F}^{p}
& \, = \, &
\norm[\Lp{[0,1]}{\lambda;E}]{K_{k}' P_{n} \otimes \id (F)}^p\\
& \, = \, &
\int_{\Omega} \int_{0}^{1}\abs[ 
\sum_{i=1}^{n}r_{i}(s)k_{i}(f_{i})(\omega)]^{p} \, ds d\mu(\omega)\\
& \, \leq \, &
c_{p}^{-p}\int_{\Omega}
           \left( \sum_{i=1}^{n}
             \abs[k_{i}(f_{i})(\omega)]^{2} \right)^{\frac{p}{2}} 
                       d\mu(\omega)\\
& \, \leq \, &
c_{p}^{-p} \norm[{[2]}]{k}^{p}\int_{\Omega}
           \left( \sum_{i=1}^{n}
             \abs[f_{i}(\omega)]^{2} \right)^{\frac{p}{2}} 
                       d\mu(\omega)\\
& \, \leq \, &
c_{p}^{-p} \norm[{[2]}]{k}^{p}\int_{\Omega}C_{p}^{p}\int_{0}^{1}
           \left( \abs[\sum_{i=1}^{n}
            r_{i}(t) f_{i}(t,\omega)]^{p} \right) \, 
                       dt\,d\mu(\omega)\\
& \, = \, &
\left(\frac{C_{p}}{c_{p}} \right)^{p} \norm[{[2]}]{k}^{p} %
\norm[\Lp{[0,1]}{\lambda;E}]{P_{n} \otimes \id (F)}^{p}\\
& \, \leq \, &
\left(\frac{C_{p}}{c_{p}}\right)^{p}\max \{ c_{p}^{-1}, c_{q}^{-1}\}^{p}%
\norm[{[2]}]{k}^{p}\norm[\Lp{[0,1]}{\lambda;E}]{F}^{p}.
\end{eqnarray*}

Where now $d_{p}=\left(\frac{c_{p}}{C_{p}}\right) \min \{ c_{p}, c_{q}\}$.
\qed
The analogue of Theorem~\ref{thm:rowtrans} is:
\begin{thm}
\label{thm:squaretrans}
\index{theorem!transference!p-completely bounded maps}%
\index{theorem!transference!square functions}%
\index{$S~$}%
For 
$1<p<\infty$
there exists a constant
$C_{p}^{\ast}$
such that, whenever
$E,F \in \mathcal{L}^{p}$
are $\mathcal{L}^{p}$-spaces,
$S\subset B(E)$
is a closed subspace, and
$u\, : \,S \rightarrow \, B(F)$
is a p-completely bounded linear map, then 
for all sequences
$g_{1},g_{2}, \ldots \in F: $
\[
\norm[F]{\left( \sum_{j=1}^{\infty} \abs[u(k_{j})(g_{j})]^{2}
    \right )^{\frac{1}{2}}}
\, \leq \, C_{p}^{\ast}
M \norm[F]{\left( \sum_{j=1}^{\infty} \abs[g_{j}]^{2}
    \right )^{\frac{1}{2}}},
\]
whenever 
$k_{1},k_{2},\ldots \, \in S$
is a sequence for which there exists $M$ such that for all 
$f_{1},f_{2}, \ldots \in E: $
\[
\norm[E]{\left( \sum_{j=1}^{\infty} \abs[k_{j}(f_{j})]^{2}
    \right )^{\frac{1}{2}}}
\, \leq \, 
M \norm[E]{\left( \sum_{j=1}^{\infty} \abs[f_{j}]^{2}
    \right )^{\frac{1}{2}}}.
\]
\end{thm}
{\bf Proof:}
It suffices again to consider finite sequences 
\[
k\, = \,  (k_{1},k_{2}, \ldots , k_n) \subset S^n.
\]
Then 
$K_{k}'$
may be viewed as an element of 
\index{$M_{2^{n}}~\otimes~S$}%
$ M_{2^{n}} \otimes S $,
and denoting 
$ u(k) \, = \, (u(k_{1}), \ldots ,u(k_{n})) $
we have
\[
\id_{2^{n}}\otimes u(K_{k}') \,= \, K_{u(k)}'.
\]
In fact, the sigma algebra $\Sigma_n$ generated by the first $n$ 
Rademacher functions contains
just $2^n$ atoms of Lebesgue measure $2^{-n}$.
If we denote $\mathcal{E}_n$ the conditional
expectation 
$ \mathcal{E}_n : \Lp{[0,1],\mathcal{B}}{\lambda} \rightarrow
\Lp{[0,1],\Sigma_n}{\lambda}$,
then
\[
(\mathcal{E}_n \otimes \id) \circ K_{k}' \circ (\mathcal{E}_n \otimes \id) %
 \, =   \,   K_{k}' \mbox{ for all } k=(k_1, \ldots ,k_n) \in B(E)^n.
\]
Denote $\delta_1, \ldots, \delta_{2^n}$ the standard unit vector basis
of 
$l^{p}_{2^n}$
and
\linebreak[4]
$I_n : l^{p}_{2^n} \rightarrow  \Lp{[0,1],\Sigma_n}{\lambda}$
the canonical 
identification 
\[
I_n : \delta_i \mapsto 2^{\frac{n}{p}}
\chi_{(\frac{i-1}{2^{-n}},\frac{i}{2^{-n}}]},\;  i=1, \ldots, 2^n.
\]
Then $(I_{n} \otimes \id )^{-1} \circ  K_{k}' \circ (I_n \otimes \id )$ is given by a matrix of operators
\[
 m^k = (m^{k}_{i,j})_{i,j=1}^{2^n} \in M_{2^n} \otimes B(E),
\]
 where
\begin{eqnarray}
\label{def:matr}
m^{k}_{ij }& = &  \begin{array}{ll}
                               \sum_{l=1}^{n}r_{l}(\frac{i}{2^{n}})r_{l}(\frac{j}{2^{n}})k_{l} 
                           & 
                                 \mbox{ for  $i,j \in \{1, \ldots , 2^{n}\} $}
                       \end{array}.
\end{eqnarray}     
Thus, by 
the p-complete boundedness of $u$:
\[
\norm{m^{u(k)}} \, \leq \, \cbnorm[p]{u}\norm{m^{k}}.
\]
Since
\[
\norm{K_{k}'} = \norm{(\mathcal{E}_n \otimes \id ) \circ K_{k}' \circ
 (\mathcal{E}_n \otimes \id )} = \norm{ m^k},
\]
an application of 
Proposition~\ref{prop:khin}  yields:
\begin{eqnarray*}
 \norm[{[2]}]{u(k)}\ \leq &  \displaystyle{D_{p} \norm{K_{u(k)}'}} &  \leq \
D_{p}\cbnorm[p]{u}\norm{K_{k}'}\\
\leq  & \displaystyle{\frac{D_{p}}{d_{p}} \norm[{[2]}]{k} } & \leq \ \frac{D_{p}}{d_{p}}M.
\end{eqnarray*}
This proves the theorem with the constant
$C_{p}^{*}\, = \; \frac{D_{p}}{d_{p}}$.
\qed
\begin{rem}  \rm
For the convenience of the reader we formulate in this
remark a result due to Asmar, Berkson and Gillespie~\cite{AsmBerGil_91a}
(item (i) below).
That theorem holds true for $p \in [1,\infty)$. But for reasons discussed in
the introduction to this section, it appears as a corollary, with an additional
constant, to
Theorem~{\ref{thm:squaretrans} } only for $1<p<\infty$.
\hspace*{\fill}~\\[-2em]
\begin{enumerate}
\item[~(i)]
Let $ G $
be a locally compact, amenable group and 
$\pi : G \rightarrow B(E)$
a strongly continuous uniformly bounded representation of $G$ on
a space $E \in \mathcal{L}^{p}$.\\
If, for a sequence $ k_{1},k_{2},\ldots \in L^{1}(G,\lambda)$, there 
exists $M$ such that for all $ f_{1}, f_{2}, \ldots \in \Lp{G}{\lambda}:$ 
\[
\int_{G}\left(
\sum_{i=1}^{\infty}\abs[\lambda_p(k_{i})f_{i}(x)]^{2}
\right)^{\frac{p}{2}}\, d\lambda(x)
\, \leq \, M
\int_{G}\left(
\sum_{i=1}^{\infty}\abs[f_{i}(x)]^{2}
\right)^{\frac{p}{2}}\, d\lambda(x),
\]
then for all $g_{1},g_{2},\ldots \in E:$
\[
\norm[E]{\left(\sum_{i=1}^{\infty}\abs[\pi(k_{i})g_{i}]^{2}\right)^{\frac{1}{2}}}
\, \leq \, M \sup_{g \in G} \norm{\pi (g)}^2
\norm[E]{\left(\sum_{i=1}^{\infty}\abs[g_{i}]^{2}\right)^{\frac{1}{2}}}.
\]
\item[(ii)]
This time we leave it to the reader to formulate the version of
the last theorem for one-parameter semigroups. 
\end{enumerate}
\end{rem}
\section[Transference of Maximal Functions]{Transference of Maximal
  Functions}
\label{sec:maxtrans}
The story of studying dilation theorems for positive contractions on 
$\mathcal{L}^{p}$-spaces had been initiated by 
\index{Ak\c{c}oglu}%
Ak\c{c}oglu~\cite{Akcoglu_75}
as a means to prove a 
\index{theorem!maximal ergodic}%
maximal ergodic theorem:
\begin{thm}[Ak\c{c}oglu]
\label{thm:maxergodic}
Let 
$ 1 \leq p < \infty $
(the case $ p = \infty$ is trivial),
$ E = \Lp{\Omega}{\mu } \in \mathcal{L}^{p}$
and let
$ T \, : \, E \rightarrow E$
be a positive contraction.
\par
Then, for some constant
$c_{p}$
depending only on $p$,
there holds true for all $f \in \Lp{\Omega}{\mu}$:
\[
\left[ \int_{\Omega} \left( \sup_{n \in {\mathbb{N}}} \, \frac{1}{n} \, 
\sum_{k=1}^{n} \abs[T^{k}f\,(\omega)] \right)^{p}\,
d\mu(\omega)\right]^{\frac{1}{p}}\, \leq \, c_{p}\, \left( \int_{\Omega}
\abs[f(\omega)]^{p}\,d\mu(\omega) \right) ^{\frac{1}{p}}.
\]
\end{thm}
\begin{rem} \rm \hspace*{\fill}~\\[-2em]
\begin{description}
\item[{(i)~~}]
If 
$T$
is a contraction for $p=1$ and for $ p=\infty$, 
then the assertion of this theorem essentially is contained in the
Hopf-Dunford-Schwartz~maximal ergodic theorem. The above theorem
then follows simply by interpolation. Its value lies in the fact that it
deals with a single $p$.
\par
\item[{(ii)~}]
If
$T$
is an invertible isometry, then the conclusion of the theorem had been shown
to hold true by
\index{Ionescu-Tulcea}%
A.~Ionescu-Tulcea. Furthermore, Ak\c{c}oglu reduces, involving
a dilation, to this case.
\item[{(iii)}]
\index{de~la~Torre}%
A.~de~la~Torre in~\cite{Torre_76} has a nice proof of the ergodic theorem
relying on ideas of transference. This is further developed in the 
monograph~\cite{CoRoWe_77} of Coifman, Rochberg and Weiss.
\item[{(iv)~}]
\index{Asmar}%
Asmar, 
\index{Berkson}%
Berkson and 
\index{Gillespie}%
Gillespie in their paper~\cite{AsmBerGil_91}  generalise
from the group of integers to uniformly bounded representations of amenable
groups by 
\index{operator!separation preserving}%
separation preserving operators.
\end{description}
\end{rem}
We shall consider in this section a locally compact group
$G$ and a strongly continuous, uniformly bounded representation
\begin{eqnarray} 
\label{eq:unibd}
\pi \, : \, G \rightarrow B(E)
\end{eqnarray}
of it on $E=\Lp{\Omega}{\mu}$, ($1\leq p < \infty$), 
by separation preserving operators.
Such a representation we shall call a 
\index{representation!separation preserving}%
separation preserving representation.
\par
Alternatively, we could deal in the remainder of this chapter with a 
representation 
$\pi$ of $G$ on $E$ such that 
\[
\sup_{x \in G} \norm{ \abs[\pi(x)] } \le \infty.
\] 
We shall not do so, but we shall shortly discuss the definition and
properties of the map
$x \mapsto  \abs[\pi(x)]$,
when $\pi$ is a separation preserving representation.
\begin{defi}
\label{def:dominate}We shall say that a positive operator $P$ %
\index{dominating}%
\index{operator!dominating}%
dominates another operator $R$ %
(resp. $R$ is dominated by $P$),
if for all $f \in \Lp{\Omega}{\mu}$
\[
\abs[Rf] \leq P \abs[f].
\]
\end{defi}
\begin{prop}
\label{prop:absval}
If
$ R \, : \, E \rightarrow E$
is separation preserving, then there exists a positive operator
$P$ 
dominating $R$,
with $\norm{P} \, = \, \norm{R}$. 
Furthermore, 
$\abs[Rf] = P\abs[f]\,  \, \mbox{ for all } f \in E $,
and this justifies to denote this operator
$\abs[R]$.
\end{prop}
{\bf Proof:}
The assertion we know already if 
$R$
is a contraction, see Lemma~\ref{lem:disjoint} \\
and Remark~\ref{rem:seppreserv}.
This we may apply to 
$\frac{1}{\norm{R}} \cdot R$.
\qed
\begin{rem} \rm
\label{rem:measurable}%
Now the map $ \abs[\pi(.)] : G \rightarrow B(E)$ is defined properly and we
should note that it is in fact strongly continuous.
\end{rem}
Actually, for a non-negative $f \in \Lp{\Omega}{\mu}$ and 
$x,y \in G$ we may estimate $\mu \mbox{-almost everywhere}$:
\begin{eqnarray*}
\abs[ | \pi (x) |  f - | \pi (y) | f] &=& %
\abs[ | \pi (x) f| - | \pi (y) f|]\\
& \leq & \abs[ \pi (x) f - \pi (y) f].
\end{eqnarray*}
Hence we obtain the continuity of
$x \mapsto \abs[\pi(x)]f$,
for non-negative $f$, from the strong continuity of the 
representation in question.
Since any element of $\Lp{\Omega}{\mu}$ is a linear combination of
at most four non-negative ones, we are done. 
\par
For a finite sequence 
$k_{1}, \ldots , k_{l} \in \ L^{1}(G,\lambda) $
consider the maximal functions
\index{$m(f)$}%
\index{$M(h)$}%
\begin{eqnarray}
\label{eq:maxfkt1}
m(f)(x) & \, = \, & \sup_{1 \leq i \leq l} \abs[\lambda_{p}(k_{i})f](x), 
\quad f \in \Lp{G}{\lambda},\\
\label{eq:maxfkt2}
M(h)(\omega) & \, = \, & \sup_{1 \leq i \leq l} \abs[\pi_{p}(k_{i})h](x), 
\quad h \in \Lp{\Omega}{\pi}.
\end{eqnarray}
Our aim is to give an estimate of $M$ in terms of $m$.
\begin{lem}
Given the above conditions and non-negative
$\alpha \in \Lp{G}{\lambda}$,
$ \beta \in L^{q}(G,\lambda)$,
there holds for all 
$f \in E$
and
$i=1, \ldots,l$, 
$\mu\mbox{-almost everywhere}$:
\[
\abs[\pi((\beta \star \alpha^{\vee}) \cdot k_i) \, f] %
\leq \! \int_{G} \!\! \beta(y) %
\abs[\pi(y)]  \sup_{1 \leq i \leq l} \left\{ %
\abs[\int_{G}\!\! \alpha^{\vee}(x) k_{i}(yx)\pi(x)f  d\lambda(x)] \right\} %
d\lambda(y). 
\]
Here $q$ is such that 
$\frac{1}{p}+\frac{1}{q}=1$,
and if $q=\infty$, then we additionally assume 
that $\beta \in L^{\infty}(G,\lambda)$
has a support of finite Haar measure.
\end{lem}
{\bf Proof:}
\begin{eqnarray*}
\pi((\beta \star \alpha^{\vee}) \cdot k_i)f &  =  & \int_{G}\int_{G} %
\beta(y)\alpha^{\vee}(y^{-1}x)k_{i}(x)\pi(x)f \  d\lambda(y)d\lambda(x) \\
& \, = \, &
 \int_{G} \beta(y)\int_{G}
  \alpha^{\vee}(x)k_{i}(yx)\pi(y)\pi(x)f \,  d\lambda(x)d\lambda(y) \\
& \, = \, &
 \int_{G} \beta(y)\pi(y)\int_{G}
  k_{i}(yx)\alpha^{\vee}(x)\pi(x)f \,  d\lambda(x)d\lambda(y). 
\end{eqnarray*}
Hence, it follows that
\begin{eqnarray*}
\abs[\pi((\beta \star \alpha^{\vee})\cdot k_i)f]&  \leq  &
 \int_{G} \beta(y)\abs[\pi(y)] \abs[\int_{G}
  k_{i}(yx)\alpha^{\vee}(x)\pi(x)f \,  d\lambda(x)]d\lambda(y)\\
&  \leq  &
 \int_{G}\! \beta(y)\abs[\pi(y)] \sup_{1 \leq i \leq l}\abs[\int_{G}\!
  k_{i}(yx)\alpha^{\vee}(x)\pi(x)f   d\lambda(x)]d\lambda(y),
\end{eqnarray*}
where the last inequality is true since 
$\abs[\pi(y)]$
is a positive operator.
\qed
\begin{thm}
\label{thm:maxthm}
\index{theorem!maximal!group}%
Assume that
$G$
is amenable and that
$\pi \, :\, G \rightarrow B(E)$
is a representation of $G$ on $E \in \mathcal{L}^{p}$ as in~{\rm (\ref{eq:unibd})}.
Let
$k_{1}, \ldots  \in L^{1}(G,\lambda)$
be a sequence and define the maximal operators corresponding to this sequence
as in~{\rm (\ref{eq:maxfkt1})} and~{\rm (\ref{eq:maxfkt2})}.
If 
${\rm m}_k$
is a constant 
such that for all 
$ f \in \Lp{G}{\lambda}$:
\begin{eqnarray*}
\norm[\Lp{G}{\lambda}]{m(f)} & \, \leq \, & {\rm m}_k \, \norm[\Lp{G}{\lambda}]{f},
\end{eqnarray*}
then for all 
$h \in E$:
\begin{eqnarray*}
\norm[\Lp{\Omega}{\mu}]{M(h)} & \, \leq \, & {\rm m}_k \, 
\sup_{y \in G} \norm{\pi(y)}^{2} \norm[\Lp{\Omega}{\mu}]{h}.
\end{eqnarray*}
\end{thm}
{\bf Proof:}
Approximating 
$Mh  = \lim_{l \rightarrow \infty}M_{l}h$
from below by functions
\[M_{l}h(\omega) =  \sup_{1 \leq  i \leq l}\abs[\pi(k_{i})h(\omega)] 
  \qquad l=1,2,\ldots \; ,\]
we see, that it is sufficient to prove the theorem in the case that the sequence
$k=(k_{1}, \ldots, k_{l})$
is finite.
\par
Now, since 
$G$
is amenable, there exist  sequences 
$ (\alpha_{n})_{n\in {\mathbb N}}, \, (\beta_{n})_{n \in {\mathbb N}}$, 
with
$0 \leq \alpha_{n} \in \Lp{G}{\lambda}, $ 
$0 \leq \beta_{n} \in \Lq{G}{\lambda}$,
$\norm[p]{\alpha_{n}}= \norm[q]{\beta_{n}}=1$
for all 
$n \in {\mathbb N}$,
such that as
$ n \rightarrow \infty $:
\[
\norm[1]{(\beta_{n} \star \alpha_{n}^{\vee}) \cdot k_{i} \, - \, k_{i}} %
\rightarrow 0 \quad \mbox{ for } i \in \{1, \ldots , l \}.
\]
(See the proof of Theorem~\ref{thm:cbrepn} for this.) 
\par
The strong continuity and the uniform boundedness of the representation imply 
the convergence of 
$\pi \left( (\beta_{n} \star \alpha_{n}^{\vee} ) \cdot k_{i} \right) h$
to
$\pi \left( k_{i} \right) h$
in $L^{p}$-norm for each
$i \in \{1, \ldots , l \}$.
Since the lattice operations are continuous, we infer the $L^{p}$-norm
convergence
of
$\sup_{1 \leq i \leq l} %
\abs[\pi \left( (\beta_{n} \star \alpha_{n}^{\vee} ) \cdot k_{i} \right) h]$
to
$\sup_{1 \leq i \leq l} 
\abs[\pi \left( k_{i} \right) h]$.
\par
Since
$\norm[p]{M_{l}h} \, = \, \sup \left\{\int_{G}\, g(\omega)M_{l}h(\omega) \, 
d\mu(\omega) \, : \, g\geq0, \norm[q]{g}=1
\right\}
$,
it suffices to estimate (with obvious modifications if $q=\infty$):
\begin{eqnarray*}
{\int_{\Omega}\, g(\omega)\sup_{1 \leq i \leq l}
\abs[\pi((\beta_{n} \star \alpha_{n}^{\vee}) \cdot k_{i}) h]
(\omega)\, d\mu(\omega)} 
& \, \leq \, & 
{\rm m}_k \,  \norm[q]{g} \, \norm[p]{h} \, \sup_{x \in G}\norm{\pi(x)}^2.
\end{eqnarray*}
Now, by the lemma:
\begin{eqnarray*}
\lefteqn{\int_{\Omega}\, g(\omega)\sup_{1 \leq i \leq l}
\abs[\pi((\beta_{n} \star \alpha_{n}^{\vee}) \cdot k_{i})h]
(\omega)\, d\mu(\omega)
 }\\
&\leq &
\int_{\Omega} \!g(\omega)\int_{G} \!\beta_{n}(y)\abs[\pi(y)] \sup_{1 \leq i \leq
  l}\abs[\int_{G} %
\! k_{i}(yx)\alpha_{n}^{\vee}(x)\pi(x)h(\omega) \,  d\lambda(x)]d\lambda(y) d\mu(\omega)\\
&  =  &
\int_{G} \!\!\beta_{n}(y)\int_{\Omega}\!\!\abs[\pi(y)]^{\ast}g(\omega) %
\sup_{1 \leq i \leq l}\abs[\int_{G} %
\!\! k_{i}(yx)\alpha_{n}^{\vee}(x)\pi(x)h(\omega) \,  d\lambda(x)]d\mu(\omega)d\lambda(y) \\
&  \leq  &
\int_{G} \beta_{n}(y)\norm[q]{\abs[\pi(y)]^{\ast}g}\norm[p]{ \sup_{1 \leq i \leq l}\abs[\int_{G}
  k_{i}(yx)\alpha_{n}^{\vee}(x)\pi(x)h(\omega) \,  d\lambda(x)]}d\lambda(y)\\
&  \leq  &
\left(\int_{G} \abs[\beta_{n}(y)]^{q}
\sup_{z\in G} \norm{\abs[\pi(z)]^{\ast}}^{q}
\norm[q]{g}^{q}d\lambda(y)\right) ^{\frac{1}{q}}\\
& & \quad
\left(\int_{G} \int_{\Omega}{ \sup_{1 \leq i \leq l}\abs[\int_{G}
  k_{i}(yx)\alpha_{n}^{\vee}(x)\pi(x)h(\omega) \,  d\lambda(x)]^{p}}
d\mu(\omega)d\lambda(y) \right)^{\frac{1}{p}}\\
&  \leq  &
\sup_{z\in G} \norm{\abs[\pi(z)]^{\ast}}\norm[q]{g}
\left(\int_{\Omega} \norm[p]{m(\alpha_{n}(.)\pi(.)h(\omega))}^{p} \, d
  \mu(\omega) \right)^{\frac{1}{p}}\\
&  \leq  &
\sup_{z\in G} \norm{\abs[\pi(z)]^{\ast}}\norm[q]{g}
\left(\int_{\Omega} {\rm m }_{k}^{p}  
\int_{G}\abs[\alpha_{n}(x)\pi(x)h(\omega)]^{p} \, d\lambda(x)
d\mu(\omega) \right)^{\frac{1}{p}}\\
&  \leq  & {\rm m}_k \,
\sup_{z \in G}  \norm{\abs[\pi(z)]^{\ast}} \, \sup_{z \in G}
\norm{\abs[\pi(z)]}\, \norm[q]{g}\norm[p]{h} ,
\end{eqnarray*}
and Proposition~\ref{prop:absval} implies the statement.
\qed
As a consequence we can prove a version of the Hopf-Dunford-Schwartz 
maximal ergodic theorem. We shall deal only with the case that our
$\mathcal{L}^p$-spaces are reflexive, i.e.\ $ p>1$. For $p=1$, 
a probability measure space $ (\Omega,\mu)$, and a single operator
T with $T({\rm 1}) = {\rm 1}$, the theorem is due to 
\index{Chacon}%
Chacon and 
\index{Ornstein}%
\index{theorem!Chacon and Ornstein}%
Ornstein~\cite{ChacOrns_60}.
But before stating and proving it we remind the reader of the following 
simple fact which we shall need in our proof of the theorem.
\begin{lem}
\label{lem:posiso}
If 
$ E = \Lp{\Omega}{\mu} ,\  E'=\Lp{\Omega'}{\mu'}$
are
$\mathcal{L}^p$-spaces,
$1<p<\infty$, and 
$ D : E \rightarrow E'$
is a 
\index{operator!positive}%
\index{isometry!positive}%
positive isometry, then for all real
valued
$f,g \in E:$
\[
D(f \vee g) = D(f) \vee D(g).
\]
\end{lem}
{\bf Proof:}
Since $D$ is positive we obtain
the inequality
\begin{eqnarray}
\label{eq:posop}
D(f \vee g) & \geq & D(f)  \vee D(g)
\end{eqnarray}
from 
$D(f \vee g) -D (f) =  D(f \vee g -f) \geq 0$ and
$D(f \vee g) -D (g) =  D(f \vee g -g) \geq 0$.
Now, let 
\begin{eqnarray*}
A&=& \{\, \omega \, : \, f(\omega)=g(\omega) \, \},\\
B&=& \{\, \omega \, : \, f(\omega) \vee g(\omega) = f(\omega)\, \}
 \setminus A ,\\
C&=& \{\, \omega \, : \, f(\omega) \vee g(\omega) = g(\omega)\, \}
 \setminus A. 
\end{eqnarray*}
Clearly,
\begin{eqnarray*}
D(f \vee g) &=& D(f \cdot  \chi_A ) +D(f\cdot \chi_B ) +D(g\cdot \chi_C )\\
             &=& D(g \cdot  \chi_A ) +D(f\cdot \chi_B ) +D(g\cdot \chi_C ).
\end{eqnarray*}
Since $D$ is positive and isometric we know from Remark~\ref{rem:isometry}~that the supports
\begin{eqnarray*}
A_1&=& {\rm supp}\ D(f\cdot \chi_A ) \, = \, {\rm supp}\ D(g\cdot \chi_A ) ,\\
B_1&=& {\rm supp}\ D(f\cdot \chi_B ),\\
C_1&=&  {\rm supp}\ D(g\cdot \chi_C ) 
\end{eqnarray*}
are mutually disjoint. Furthermore,
\begin{eqnarray*}
A_1  \cap {\rm supp}\ D(f\cdot \chi_{\Omega \setminus A}) & = & \emptyset,\\
B_1  \cap {\rm supp}\ D(f\cdot \chi_{\Omega \setminus B}) & = & \emptyset,\\
C_1  \cap {\rm supp}\ D(g\cdot \chi_{\Omega \setminus C}) & = & \emptyset.
\end{eqnarray*}
Hence on  $A_1,B_1,C_1$, respectively, we have:\\
\begin{eqnarray*}
D(f \vee g) & = & D(f\cdot \chi_A ) \, = \, D(f)  \, \leq \, D(f) \vee D(g),\\
D(f \vee g) & = & D(f\cdot \chi_B ) \, = \, D(f)\, \leq \, D(f) \vee D(g),\\
D(f \vee g) & = & D(g\cdot \chi_C ) \, = \, D(g)\, \leq \, D(f) \vee D(g).
\end{eqnarray*} 
\qed
\vspace{-2em}
\begin{thm}
\label{thm:dunschw}
\index{theorem!maximal!semigroup}%
Let
$(T_t)_{t \geq 0}$ be a strongly continuous one-parameter semigroup
of positive contractions acting on some $\mathcal{L}^p$-space
$E=\Lp{\Omega}{\mu}$. For $p>1$ there exists a constant $C_p$,
depending only on p, such that for all $f \in E$:
\[
\norm[E]{ \sup_{t>0} \ \abs[ \frac{1}{t} \int_{0}^{t} \ T_s f \ ds]}
\, \leq \, C_p \norm[E]{f}.
\] 
\end{thm}
{\bf Proof:}
First we note that we are taking a supremum, pointwise, over an uncountable
set of functions
$\omega \mapsto \frac{1}{t} \int_{0}^{t} \ T_s f(\omega)\ ds\ $. 
But 
$ t \mapsto  T_t f  $
is continuous from
$(0, \infty)$ to $ \Lp{\Omega}{\mu}$ and $ {\mathbb{Q}}_{+} $
is a countable dense subset in $ (0,\infty)$,
hence, as Dunford and Schwartz show on page 686 of their work~\cite{DunfSchw_58}: 
\begin{eqnarray}
\label{eq:dense}
\sup_{t>0} \ \abs[ \frac{1}{t} \int_{0}^{t} \ T_s f \ ds] & = &
\sup_{t \in  {\mathbb{Q}}_{+}} \ \abs[ \frac{1}{t} \int_{0}^{t} \ T_s f \ ds] \, ,
\end{eqnarray}
except on a possibly  $ f $-dependent set of $ \mu $-measure $0$.
We note further, that it suffices to prove  that for each finite subset
$ \{t_1, \ldots , t_n\} =:N \subset {\mathbb{Q}} $:
\[
\norm[E]{ \sup_{t \in N} \abs[\frac{1}{t} \int_{0}^{t}
T_s f \ ds ]} \leq C_p \norm[E]{f}.
\]
Because then the monotone convergence theorem will provide the 
estimate for
the right hand side of the above equality~(\ref{eq:dense}).
\par
Now, for $t>0$, we let 
$k_t= \frac{1}{t} \chi_{[0,t]}$.
The corresponding maximal operator on 
$\Lp{{\mathbb{R}}}{\lambda}$:
\[
m(h)(x) = \sup_{t>0} \abs[ \lambda_p(k_t)h \ (x)] =
\sup_{t>0} \ \abs[ \frac{1}{t} \int_{0}^{t} \ h(x-s) \ ds], 
\]
is just the left sided 
\index{operator!maximal!Hardy-Littlewood}%
\index{maximal function!Hardy-Littlewood}%
Hardy-Littlewood maximal operator, for 
whose $L^p$ boundedness we refer the reader to \cite{HardLitt_30}
and the theorems 384 and 398 of \cite{HaLiPa_52}.
\par
Now we apply the dilation theorem, Theorem~\ref{thm:dilthm} to the semigroup to write
\[
D \circ T_t = P\circ S_t \circ D , \qquad t\geq 0.
\]
Our preparing Lemma~\ref{lem:posiso} and Theorem~\ref{thm:maxthm}
yield:
\pagebreak[4]
\begin{eqnarray*}
\norm[E]{\sup_{t \in N} \ \abs[ \frac{1}{t} \int_{0}^{t} \ T_s f \ ds] } 
& = & \norm[\tilde{E}]{ D \left(\sup_{t \in N} \ \abs[ \frac{1}{t} \int_{0}^{t} \ T_s f \ ds] \right)} \\
& = & \norm[\tilde{E}]{ \sup_{t \in N} \ \abs[ \frac{1}{t} \int_{0}^{t} \ D T_s f \ ds] } \\
& = & \norm[\tilde{E}]{ \sup_{t \in N} \ \abs[ \frac{1}{t} \int_{0}^{t} \ P S_sD f \ ds] } \\
& \leq & \norm[\tilde{E}]{ \sup_{t \in N} \ \abs[ \frac{1}{t} \int_{0}^{t} \ S_sD f \ ds] } \\
& = & \norm[\tilde{E}]{ \sup_{t \in N}\  M(D f) }  \, \leq \, C_p \norm[\tilde{E}]{Df} \\
& = &  C_p \norm[E]{f},
\end{eqnarray*}
where $ M $ is the maximal operator corresponding to the representation
$t \mapsto S_t$ and the finite sequence
$ k= (k_t)_{t \in N}$.
\qed 
\chapter[Submarkovian Semigroups]{Submarkovian Semigroups}
\section[Examples of Sub-Markovian Semigroups]{Some Examples of Sub-Markovian
  Semigroups}
\label{sec:examples}
For a $\sigma$-finite measure space $(\Omega,\mu)$ we consider a semigroup $(T_t)_{t\ge 0}$
acting ``simultaneously'' on all $L^p(\Omega,\mu)$ spaces, 
$1 \le p \le \infty$,
such that for all
\index{semigroup}
$ f \in L^2(\Omega,\mu) \cap L^p(\Omega,\mu) $ and all $ t \geq 0:$
\begin{eqnarray}
\label{eq:conditions}
\begin{array}{lrcl}
\mbox{\it i)~~} \quad \qquad &  \norm[p]{T_t f} & \le & \norm[p]{ f } \\
\mbox{\it ii)~}  \quad  \qquad &   T_t^\ast & = & T_t \; \mbox{\rm on}\; L^2(\Omega,\mu) \\
\mbox{\it iii)~} \qquad   \quad &    T_0 & = & \id\\
\mbox{\it iv)~~}\quad  \qquad &   T_t f &  \ge & 0 \; \mbox{\rm if}\; f \ge 0.
\end{array} & & 
\end{eqnarray}
Furthermore, for $ p < \infty$ we assume strong continuity of the map
\begin{eqnarray*}
 T & : & t  \mapsto  T_t, \\
 T & : & [0, \infty )   \rightarrow   L^p(\Omega,\mu).
\end{eqnarray*}
In this chapter we require $\sigma$-finiteness of the measure space to
use the duality $L^1(\Omega,\mu)^{\ast} = L^{\infty}(\Omega,\mu)$.
\par
The above conditions are surely to a large extend redundant and to an even
larger extend not really needed for our development. To be more precise we
define:
\begin{defi}
\label{defi:submar}
A one-parameter semigroup
$(T_t)_{t\geq 0}$
acting strongly continuously on 
$L^2(\Omega,\mu)$
is called a 
\index{submarkovian}%
\index{semigroup!submarkovian}%
submarkovian semigroup if:
\[
\begin{array}{lrcl}
\mbox{\rm i)~~~} \quad \qquad &  \norm[2]{T_t f} & \le & \norm[2]{ f } \\
\mbox{\rm ii)~~} \quad  \qquad &   T_t^\ast & = & T_t \; \mbox{\rm on}\; L^2(\Omega,\mu) \\
\mbox{\rm iii)~} \qquad   \quad &    T_0 & = & \id\\
\mbox{\rm iv)~~} \quad  \qquad &   T_t f &  \ge & 0 \; \mbox{\rm if}\; f \ge 0\\
\mbox{\rm v)~~~}\qquad \quad &  \norm[\infty]{T_t f} & \le & \norm[\infty]{ f } 
\mbox{\rm if}\;f \in L^2(\Omega,\mu) \cap L^{\infty}(\Omega,\mu).
\end{array}
\]
\end{defi}
\begin{rem} \rm \hspace*{\fill}~\\[-2em]
\begin{description}
\item[(i)~]
A use of the selfadjointness, the norm bounds and an application of the 
\index{Marcinkiewicz}%
Marcinkiewicz
\index{theorem!Marcinkiewicz}%
\index{theorem!interpolation}%
interpolation theorem (see e.g.\ \cite{SteiWeis_71} Chap. 5 sect. 2)
implies that our original set of conditions is fulfilled by a submarkovian 
semigroup. Strictly speaking we should notationally distinguish
the different extensions 
$(T_{t}^{p})_{t \geq 0}$,
obtained this way on the different $L^p$-spaces, $1 \leq p < \infty $.
But in the sequel it is always understood
that on 
$L^p$ we are considering the extension by continuity of an operator
defined on
$L^p \cap L^2$.
\item[(ii)]
Thus far our presented results are all proved for a single
$ p \in (1, \infty)$,
but now we need a little more to define functions of the negative of the
infinitesimal generator
$-A_{p}$ on $L^p$. At the moment we see at least two ``philosophies'':
\begin{enumerate}
\item
Extend everything from $L^p \cap L^2$. 
\item
Assume, for one $p \in (1,\infty)$, some analyticity of the semigroup.
\end{enumerate}
In the first, more traditional case, to which we essentially shall stick on,
we use von~Neumann's spectral calculus as a start.
Then a good set of conditions is the above {\rm i)-iii)}, and instead of 
{\rm iv)} and {\rm v)}, that for
some 
$r_0 \in (1,2)$  
$(T_t)_{t\geq 0}$ extends, by continuity, to a semigroup of sub-positive
contractions on $L^{p}, \quad r_0 \leq p \leq r_{0}'$. We shall return to this
in Remark~(\ref{rem:angleanalytic}).
\par
In the second case, 
one can define functions of the generator, in fact more elementarily, by the
integral calculus given for the resolvent of the infinitesimal generator.
If one still considers only a single $p \in (1, \infty)$,
it would be necessary for us to have a semigroup of sub-positive
contractions which is strongly continuous on
$\Lp{\Omega}{\mu}$. Some of our following conclusions are improved
by interpolation, and such as the maximal theorem, Theorem~\ref{thm:aeconv},
essentially need selfadjointness on $L^2$.  
\end{description}
\end{rem}
We denote $A$ the
\index{generator}%
\index{generator!infinitesimal!on $L^2$}%
infinitesimal generator of $(T_t)_{t \ge 0}$ on $L^2$, such that
\begin{eqnarray*}
T_t & = & e^{tA},\\
Af & = & \lim_{t \searrow 0}\ \frac{T_t-1}{t}\ f, \\
D(A) & = & \{ f\in L^2(\Omega,\mu): \lim_{t \searrow 0}\ \frac{T_t -1}{t}\ f
\; \mbox{\rm exists in norm} \}\,.
     \end{eqnarray*}
     Then $A$ is a selfadjoint operator with $\sigma (A) \subset {\mathbb{R}}_{-}
     \cup \{ 0 \}, $ hence $-A$ is positive.
\par
Further, we know from von~Neumann's spectral theory that there exists a unique
spectral resolution of the identity $(P_\lambda)_{\lambda \in {\mathbb{R}}}$, such
that
\begin{eqnarray*}
Af & = & \int^0_{-\infty}\ \lambda\ dP_\lambda f \qquad f \in D(A),\\  
T_t f & = & \int^0_{-\infty}\ e^{t\lambda}dP_\lambda f \qquad f \in L^2\,.
\end{eqnarray*}
We shall suppose $P_0 = 0$, i.e.\ $A$ is $1-1$ on $D(A)$.
\par
Next we recall some examples of semigroups, contractive on some
$\mathcal{L}^{p}$-spaces.
\begin{exam} \rm \hspace*{\fill}
\label{examples}
\begin{description}
\item[(i)~~] On $L^p({\mathbb{R}},\lambda)\; \mbox{\rm let}\; T_s f(.) = f(.-s),\;
f \in L^p({\mathbb{R}},\lambda),\; s \ge 0 $, denote the semigroup of
\index{semigroup!translation-operators}%
translation-operators, which are not selfadjoint.  Note that we artificially
made a one-parameter semigroup out of an, on $L^2({\mathbb{R}},\lambda)$, unitary
group.
\par
The generator $A$ is an extension of $ - \frac{d}{dx} : f \mapsto - f'$ from
the space 
$\{ h : h'\, \mbox{\rm exists and}\, h' \in L^2 ({\mathbb{R}},\lambda) \}$. 
Using the notation 
$\hat{} : f \mapsto \hat{f}$ for the
\index{$\hat{f}$}%
\index{Fourier transform}%
\index{transform!Fourier }%
Fourier transform, defined on $ L^1 ({\mathbb{R}},\lambda) \cap
L^2({\mathbb{R}},\lambda) $ by:
\[
\hat{f}(\xi) = \int_{\mathbb{R}} e^{-ix \xi} f(x) \ dx , \qquad \xi \in {\mathbb{R}},
\]
we have
\begin{eqnarray*}
(T_t f)\hat{\ } (\xi) & = & e^{- i t \xi} f(\xi), \\
(Af)\hat{\ }  (\xi) & = & -i \xi \hat{f} (\xi).
\end{eqnarray*}
\item[(ii)~] Convolution with
\index{semigroup!Gaussian}%
Gaussian kernels on $L^p({\mathbb{R}}^n, \lambda)$,
\[
T_t f (x) = \frac{1}{(2 \pi t)^{\frac{n}{2}}}\ \int_{{\mathbb{R}}^n}\ 
e^{-\frac{(x-y)^2}{2t}} f(y)dy\,.
\]
Here the generator $A$ is one half of the
\index{Laplacian}%
Laplacian,
\[
Af(x) = \frac{1}{2}\triangle f(x) = \frac{1}{2}\sum^n_{i=1}\ \left(
  \frac{\partial}{\partial x_i} \right)^2 f(x).
\]
\par
On the side of the Fourier transform:
\begin{eqnarray*} 
(T_t f)\hat{\ } (\xi ) & = & e^{-\frac{t}{2}\abs[\xi ]^2} \hat{f}(\xi ), \\
(Af)\hat{\ }  (\xi)& = & - \frac{\abs[\xi]^2}{2} f(\xi) \,.
\end{eqnarray*}
\item[(iii)] \index{semigroup!Poisson} Poisson semigroup on 
$L^p({\mathbb{R}}^n,\lambda)$,
\begin{eqnarray*}
T_t f(x) & = &{ \left( \frac{1}{\pi} \right)^n %
\int_{{\mathbb{R}}^n}\ \frac{t}{(t^2 + \abs[y]^2)^{\frac{n+1}{2}}}\ f(x-y)dy }, \\
A & = & - (- \triangle)^{\frac{1}{2}}.
\end{eqnarray*}
For the Fourier transforms this means:
\begin{eqnarray*}
(T_t f)\hat{\ } (\xi) & = & e^{-t \abs[\xi ]} \hat{f}(\xi) ,\\
(Af)\hat{\ } (\xi) & = & - \abs[\xi ] \hat{f}(\xi).
\end{eqnarray*}
\item[(iv)~] An example not fulfilling all of the above conditions:\\
On $L^p({\mathbb{R}},\lambda)$ let for $t \geq 0:$
\[
S_t\ f(x) = e^{\frac{t}{p}} f(e^t x)\,.
\]
Then, only on the fixed $\mathcal{L}^p$-space, $ \norm{S_t } \le 1$.
On that space the semigroup is strongly continuous,
and additionally all the $S_t$ are positive operators.\\
For $f \in C^1:$
\begin{eqnarray*}
  Af(x) & = & \lim_{t\searrow 0} \frac{ e^{\frac{t}{p}} f(e^t x) - f(x)}{t}
 \\
 & = & \lim_{t \searrow 0} \frac{e^{\frac{t}{p}} (f(e^t x) - f(x)) }{t}+
\frac{e^{\frac{t}{p}}-1}{t}\ f(x)  \\
 & = &  f'(x) \cdot x\ + \frac{1}{p} f(x).
\end{eqnarray*}
\item[(v)~~]
Evolution semigroups to the Schr\"{o}dinger operator with certain
perturbations, for example point interactions
are not always bounded on the whole range $1 \leq p < \infty$. 
In three dimensions one has to deal with
a semigroup which is  uniformly
bounded on $\Lp{\mathbb{R}}{\lambda}$ only for $p \in (\frac{1}{3},3)$.
(see \cite{CaspClem_94} and 
\cite{AlbBrzDab_95} for more on this.)
\end{description}
\end{exam}
\section[Fourier $L^p$-Multipliers]{Fourier $L^p$-Multipliers}
In this section we shall discuss some properties of distributions which map,
by convolution, the space $\Lp{{\mathbb{R}}}{\lambda}$ into itself.
We shall be interested in the case $ 1<p<\infty$.
Here, and in the sequel, by a  
\index{distribution}%
\index{distribution!Schwartz}%
distribution we mean
a functional $\phi$ defined on the
\index{Schwartz space}%
Schwartz space 
$\mathcal{S}({\mathbb{R}})$. For information on this class, the class of tempered
distributions, we refer the reader to Chap.1 sect.3 of the book of Stein and Weiss~\cite{SteiWeis_71}. 
In a lot of cases the distributions, we have to deal with, 
will be given by integration against a function which we
shall denote $x \mapsto \phi(x), \, x\in{\mathbb{R}}$.
\begin{defi}
\label{defi:p-convolver}
A tempered distribution
$\phi$
is called an
\index{convolver}%
\index{convolver!$L^p$}%
$L^p$-convolver,
respectively its Fourier transform
$\hat{\phi}$ 
is called an 
\index{multiplier}%
\index{multiplier!$L^p$}%
$L^p$-multiplier, if for some 
$ C > 0 $ for all 
$ f \in \mathcal{S}({\mathbb{R}})$
\[
 \norm[p]{ (\hat{\phi} \cdot \hat{f})\ \check{} } \le C \norm[p]{ f }.
\]
Here, 
\index{$\check{\varphi}$}%
$\, \check{}: \varphi \mapsto \check{\varphi}$ 
denotes the inverse Fourier transform.
\end{defi}
\begin{rem} \rm
\label{rem:p-conv}%
Sometimes it is useful to note that
$\mathcal{L}^p$-convolvers are just those bounded operators
on $ L^p(\mathbb{R})$
which commute with all translations (c.f.\
\cite{SteiWeis_71} chap. I). 
\end{rem}
Some remarks on the space $Cv_p({\mathbb{R}})$ of all 
$L^p$-convolvers seem to be useful for our development.
\begin{rem} \rm \hspace*{\fill}~\\[-2em]
\label{rem:Cv_p}
\begin{description}
\item[(i)~~]
From our definition it is natural to consider 
\index{$Cv_p({\mathbb{R}})$}%
\index{algebra!$Cv_p({\mathbb{R}})$}%
$Cv_p({\mathbb{R}})$
as a sub-algebra of the bounded operators on 
$\Lp{{\mathbb{R}}}{\lambda}$.
In there, it is not only norm closed, but it can be shown to be closed in the
weak operator topology. For 
$p=2$ we just deal with the von~Neumann algebra of ${\mathbb{R}}$. 
\item[(ii)~]
If 
$(u_k)_{k=1}^{\infty}$
is a norm bounded approximate identity in 
$L^1({\mathbb{R}}, \lambda)$, and if
$\phi \in Cv_p({\mathbb{R}})$,
then
$\phi \star u_k$
tends to 
$\phi$
in the strong operator topology, as
$k \rightarrow \infty$.
This provides a method to approximate a convolver,
some awkward distribution, by quite regular ones, i.e.\ some given by
integration against 
$C^{\infty}$-functions.
\item[(iii)]
One should be a bit careful when defining the pointwise product of a
convolver and a function. That is, it is not always true that this product is
again a convolver, or even exists as a distribution, if the function is only 
assumed to be, say,
bounded and continuous.
\par
There is anyway a nice class of functions for which
the pointwise product can be defined, moreover, 
this class constitutes an algebra
\index{$B_p({\mathbb{R}})$}%
$B_p({\mathbb{R}})$, with respect to pointwise
multiplication. This is the so called algebra of
\index{algebra!Herz-Schur}%
Herz-Schur multipliers (see~\cite{Fendler_87}, \cite{Herz_73}, 
\cite{BozeFend_84} for more on this).
\item[(iv)~]
A purely operator theoretic definition of  $B_p({\mathbb{R}})$
can be given as follows:\\
 $B_p({\mathbb{R}})$ is just the space of weakly continuous, 
p-completely bounded, pointwise multipliers of 
$Cv_p({\mathbb{R}})$. Though this sounds a bit awkward, simply because of its
length, it has been on the roots of introducing the notions of complete and
p-complete boundedness. We mention here only that p-completely bounded
maps have nice functorial properties.
\par
An alternative characterisation runs as follows:
A function 
$\varphi$ belongs to  $B_p({\mathbb{R}})$, if there
exists a subspace of a quotient-space of an 
$\mathcal{L}^p$-space, call it E,
and a strongly continuous, uniformly bounded rep\-re\-sen\-ta\-tion 
$S:{\mathbb{R}} \rightarrow B(E)$ on it, such that
$\varphi(t) = (S_t \xi, \eta), \quad t \in {\mathbb{R}}$,
for some 
$\xi \in E, \, \eta \in E^{\ast}$~\cite{BozeFend_84}~\cite{Fendler_87}.
\item[(v)~~]
It is useful to include ``subspaces of quotients'' of 
$\mathcal{L}^p$-spaces in the above, since any Hilbert space is of this type
and hence any continuous positive definite function $\varphi$
is in $B_p({\mathbb{R}})$. Moreover, for some such $\varphi$ one has:
$\norm[p,p]{\varphi \cdot \phi} \leq \varphi(0)\norm[p,p]{\phi}$ for all $L^p$-convolvers 
$\phi$.  
\par.  
\par
This gives a tool for a further regularisation of an convolver. Namely it can
be boundedly weakly approximated by convolvers with compact support, at least
on the amenable group ${\mathbb{R}}$.
\par
We shall use the functions
$v_j(x) = \max \{ 1 - \frac{\abs[x]}{j},0\}, \quad x \in {\mathbb{R}}$,
$j \in {\mathbb{N}}$
for this.
\end{description}
\end{rem}
\par
There are some celebrated theorems providing non-trivial
$L^p$-convolvers. Because of its relation to analytic functions
the 
\index{theorem!H\"{o}rmander}%
\index{theorem!multiplier}%
multiplier theorem of
\index{H\"{o}rmander}%
H\"{o}r\-man\-der, see sect.\ 2 of~\cite{Hoermander_60}, will be of importance to us and we
want to introduce the reader to this theme next.
\begin{defi}
\label{defi:hoermcond}
We say that a distribution  $\phi$ 
fulfils the Mihlin-H\"{o}rmander 
\index{conditions!Mihlin-H\"{o}rmander}%
conditions with constant
$ C>0$, if
\begin{eqnarray*}    
       \abs[ \hat{\phi}(\nu)] & \le& C,\\
       \abs[ \nu \frac{\partial}{\partial \nu}\ \hat{\phi}(\nu)  ] & \le & C,\ 
 \qquad \forall \nu \in {\mathbb{R}}.
\end{eqnarray*}
\end{defi}
The H\"{o}rmander multiplier theorem then asserts that, for 
$1 < p < \infty$, 
\[
 \norm[p]{ \phi  \star  f } \le C \cdot C_1 [ p + (p-1)^{-1}] \norm[p]{ f } 
\quad \forall \  f \in L^p \,.
\]
\begin{rem} \rm
In the above we wrote on the right hand side the dependence on $p$ of a
constant. This dependence enters in the proof of  a result of Cowling,
which we shall state as Corollary~\ref{cor:genpowers}. 
\par
The theorem of
H\"{o}rmander, we in the above referred to, is stated as theorem 2.5 in 
\cite{Hoermander_60}. Regarding the assumptions
of that theorem it would be more correct to call the condition 
required in there ``H\"{o}rmander's condition'' and the one we stated in
the definition ``Mihlin's condition'' (c.f.\ \cite{Mihlin_56}, resp.\ 
Theorem~2~\cite{Mihlin_65}, Appendix).
In our one-dimensional setting the square-integrability condition,
H\"{o}rmander requires, on the derivative of $\phi$ is slightly
weaker than the uniform growth, respectively decrease, estimate which is
required above. It should be noted, however, that in the $n$-dimensional case
H\"{o}rmander reduced the differentiability requirement from the order $n$ to
the least integer not less than half of the dimension.  
\par
For a proof of the multiplier theorem, we used, we
refer to Theorem 3 and its corollary in chapter IV \S3 of 
the book of Stein~\cite{Stein_70}.
\par
We agree, if the interested reader finds it hard to chase for the constant
there. So we shall give some more arguments here. The multipliers we just
consider are, in fact bounded from $L^1$ to weak-$L^1$. Interpolating now with
the $L^2$--$L^2$ estimate by means of the Marcinkiewicz interpolation
theorem  exhibits the dependence of the constant on $p$. Alternatively one
can show that the multipliers are bounded from the real Hardy space
$H^1_{\mathbb{R}}$ to $L^1$ and interpolate again. 
\end{rem}
If, for 
$\theta \in (0,\pi), \quad \Gamma_\theta$ 
denotes the 
\index{cone}%
\index{cone!of angle $\theta$}%
cone
\[      \begin{array}{rclcl}
  \Gamma_\theta & = & \{ \ z \in {\mathbb{C}}\setminus\{0\} \ : \  \abs[\arg z ]< \theta \ \} 
& = & \{ \ \rho e^{i\psi} \ : \  \rho > 0, \abs[ \psi ] < \theta \} \, \\
         \end{array} ,
\]   
where we represented 
$z \not= 0$ 
\index{coordinates}%
\index{coordinates!polar}%
as $ z =  \rho e^{i \arg z},\mbox{ with } \arg z \in [- \pi, \pi),
\; \rho \ > 0 \  $,
then a bounded holomorphic function $m$ on $\Gamma_\theta$ has almost 
everywhere non-tangential 
\index{boundary values}%
boundary values and we denote $m$ again the boundary function
$m : \partial \Gamma_\theta \rightarrow {\mathbb{C}}$.
\par
The cone 
$\Gamma_{\psi}$
is strictly contained in
$\Gamma_\theta$ for 
$\psi \in (0,\theta)$,
and we define 
a function
$m_{\psi}$
on 
$\mathbb{R}$:
 \begin{eqnarray}
\label{eq:boundaryvalues}
m_{\psi} (x) & = & \left\{ \begin{array}{ll}
                      m(x e^{i \psi}) & \mbox{ if } x \geq 0\\
                      m(\abs[x] e^{- i \psi}) & \mbox{ if } x < 0.
                     \end{array}
              \right.
\end{eqnarray}
\begin{defi}
\label{defi:boundarydistribution}
Given
$m$ on $\Gamma_\theta$
as above, we define a distribution
$\phi$
by:
\[ 
\phi(h) = \lim_{\psi \nearrow \theta} \int m_{\psi}(x) \hat{h}(-x) \, dx,
\qquad h\in \mathcal{S}.
\]
We then shall say that the Fourier transform 
of $\phi$ coincides with the boundary values of $m$ on $\Gamma_\theta$.
\end{defi}
\begin{exam} \rm
\label{ex:complpow}
An interesting example is the distribution
$\phi$ 
with Fourier transform
\[
  \hat{\phi}(\nu) = (i \nu)^{i\gamma}.
\]
In this case 
\[
 \phi(x) = \left\{ \begin{array}{ll}
   \Gamma (-i \gamma)^{-1}\, x^{i \gamma -1} & x > 0\\
    0 & x \le 0. 
  \end{array} \right.
\]
$ \hat{\phi}$ 
is the boundary value of the, in 
${\mathbb{C}} \setminus [- \infty, 0]$
holomorphic, function
\[
m_{\gamma} (z) = z^{i\gamma} = \exp(i\gamma \log |z| - \gamma \arg z),
\qquad z \in {\mathbb{C}}.
\]
We note that $m$ fulfils the 
Mihlin-H\"ormander conditions~\ref{defi:hoermcond}:
\begin{eqnarray*}
   \abs[m_\gamma (i \nu)] & \le & e^{\abs[\gamma] \frac{\pi}{2}}, \\
   \abs[ \frac{\partial}{\partial \nu} m_\gamma (i \nu)] & = & \abs[ i \frac{\gamma}{\nu}
   e^{\pm \gamma \frac{\pi}{2}} ] \ \le \ \frac{1}{\abs[\nu ]} \abs[ \gamma]\ 
   e^{\abs[\gamma]\frac{\pi}{2}}   \quad \nu \in {\mathbb{R}}.
\end{eqnarray*}
\end{exam}
\par
\begin{figure}[b]
\begin{minipage}{8cm}
Cauchy's 
\index{formula!Cauchy}%
formulae:
\begin{eqnarray}
\label{eq:cauchy1}
m(\zeta) & = & \frac{1}{2 \pi i} \oint_{C_r} \frac{m(z)}{z-\zeta}\ dz \\
\label{eq:cauchy2}
\frac{d}{d\zeta}m  (\zeta) & = & \frac{1}{2 \pi i} \oint_{C_{r}} \frac{m(z)}{(z-\zeta)^2}\ dz
\end{eqnarray}
\end{minipage}
\hfill
\begin{minipage}{6cm}
\begin{picture}(0,0)%
\includegraphics{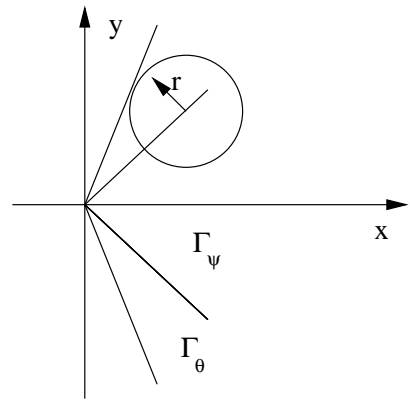}%
\end{picture}%
\setlength{\unitlength}{4144sp}%
\begingroup\makeatletter\ifx\SetFigFont\undefined%
\gdef\SetFigFont#1#2#3#4#5{%
  \reset@font\fontsize{#1}{#2pt}%
  \fontfamily{#3}\fontseries{#4}\fontshape{#5}%
  \selectfont}%
\fi\endgroup%
\begin{picture}(2499,2499)(439,-1738)
\end{picture}%
\caption{Two Cones}
\label{fig:gamma}
\end{minipage}
\end{figure}

If 
$\overline{\Gamma_\psi} \subset \Gamma_\theta$ 
are two cones with
$\psi < \theta$, like in Figure~\ref{fig:gamma},
then the ``boundary'' value on 
$\Gamma_\psi$ 
of a function which is bounded and  holomorphic
on the larger cone
$\Gamma_\theta $,
fulfils the conditions of H\"{o}rmander's Fourier
multiplier theorem. This is a consequence of the 
\index{Cauchy}%
\index{theorem!Cauchy}%
Cauchy integral theorem.
\par
In fact, if $\zeta= \abs[x] e^{i \psi} \in \partial \Gamma_{\psi}$, then
a circle $C_{r}$ with centre $\zeta$ and radius $r$ lies in $\Gamma_{\theta}$,
provided $r<\sin(\theta -\psi)\abs[x]$ (compare Figure~\ref {fig:gamma}).
It follows from~(\ref{eq:cauchy1}) and~(\ref{eq:cauchy2})  that
\begin{eqnarray*}
\abs[m_{\psi}(x)] & \leq & \frac{1}{2 \pi} 
\int_{C_r} \ {\abs[m(z)]} \, \frac{\abs[dz]}{r} \, 
\leq \, {\rm M}_{\theta},\\
\abs[\, x \ \frac{d}{dx}m_{\psi}(x)\, ] & \leq & \frac{\abs[x]}{2 \pi} 
\int_{C_r} \ {\abs[m(z)]}\, \frac{\abs[dz]}{r^2} \, 
\leq \, {\rm M}_{\theta} \frac{\abs[x]}{r},
\end{eqnarray*}
where 
${\rm M}_{\theta} = \sup_{z \in \Gamma_{\theta}} \abs[m(z)]$
denotes the least upper bound of 
$\abs[m]$ on $\Gamma_{\theta}$.
The conditions are thus fulfilled for  
$m_{\psi}$, with bound
$  {\rm M}_{\theta} \sin^{-1}(\theta-\psi)$.
\section[Application to Functional Calculus]{Application to Functional
  Calculus}
Now let 
$(T_t)_{t \geq 0}$
be a semigroup, related to a measure space
$ (\Omega,\mu)$,
as in section~(\ref{sec:examples}) which fulfils the conditions 
i)--iii) of~(\ref{eq:conditions}).
We denote $A $ its infinitesimal generator on
$L^2(\Omega,\mu)$.
If $m$ is a bounded Borel function on 
$[0,\infty)$, 
then 
$\int^{0}_{-\infty}\ m (- \lambda)dP_\lambda = m(-A)$ 
is defined on 
$L^2(\Omega,\mu)$.
For those semigroups
\index{Cowling}%
Cowling~\cite{Cowling_83} discusses the question whether 
$m(-A)$ 
extends to a
bounded operator on 
$L^p(\Omega,\mu)$ 
and considers the consequences of this functional calculus.
I shall sketch some of his results and proofs in this section.
\par
We first note that no assumption on positivity, as in~(\ref{eq:conditions}) iv)
is made. The considered class is thus more general than the class of
submarkovian semigroups.
But there is still a relation to positivity. Namely, a contraction on
$L^1$ is automatically a sub-positive contraction 
(see e.g.\ Theorem~1.1 in chapter~4~\S  \cite{Krengel_85}).
This enables one to make use of the dilation theorems and their implications.
\par 
If $P_0 \not= 0$ then to define $m(-A)$ it would be necessary to know $m(0)$
which might not be defined by the above argumentation. Not to worry about
this was the reason to assume $P_0 = 0$. Furthermore, if a uniformly 
bounded sequence 
$(m_k)_{k \in \mathbb{N}}$ of
 bounded Borel functions on $(0,\infty)$
converges to some function $m$ uniformly on compact sets in  $(0,\infty)$, 
then $m_k (-A)$ converges to $m(-A)$ in the weak operator
topology of $B(L^2(\Omega,\mu))$,
since the spectral
measure is concentrated on $(-\infty,0)$.
 
\par
Now as a corollary of our previous work we obtain a theorem of M.~Cowling,
Theorem~1 in~\cite{Cowling_83}, which is his starting point for discussing
the work~\cite{Stein_70} of 
E.~Stein.
Our proof is a variation of that one given by  Cowling.
\begin{thm}
\label{thm:multiplier}
\index{theorem!multiplier}
Let $m$ be bounded, holomorphic in 
$\Gamma_{\frac{\pi}{2}}$. Suppose that
the distribution $\phi$ 
whose Fourier transform coincides  almost everywhere with the
boundary values of 
$m$
fulfils, for some finite constant 
$C_p$ 
depending only on $p$,
$ 1 < p < \infty $,
\[
 \norm[p]{ \phi \star f } \le C_p \norm[p]{ f } \quad \forall\, f \in
 C^{\infty}_{cp}({\mathbb{R}}) \,.
\]
Then for all $ f \in L^p(\Omega,\mu) \cap L^2(\Omega,\mu) $
\[
   \norm[p]{ m (-A) f } \le C_p \norm[p]{ f }\,.
\]
\end{thm}
{\bf Proof:} The distribution
$\hat{\phi}$ can be extended to a function, still denoted 
$\hat{\phi}$,
which is holomorphic  and bounded in the lower half-plane,
since 
$m$ 
is such in 
$\Gamma_{\frac{\pi}{2}}$. 
By the Paley-Wiener theorem
the support of 
$\phi$ 
is contained in 
$[0,\infty)$.\\
Next one regularises
$\phi:$
\begin{description}
\item[(i)~]
If, for $k \in {\mathbb{N}}$, $ u_k \in  C^\infty_{cp} $ is such that
\[
 \mbox{\rm supp}\ u_k \in [ \frac{1}{k},\ \frac{2}{k} ]
\mbox{ \rm and } \int_{\mathbb{R}}\ u_k (x) \ dx  = 1\ ,
\]
then 
\[
\mbox{\rm supp}\ \phi \star u_k \subset [ \frac{1}{k},\infty),\, \phi  \star  u_k
\in L^p \cap C^\infty.
\]
Furthermore, for 
$f \in \Lp{{\mathbb{R}}}{\lambda}$ we have convergence in $ L^p$-norm:\\
$ u_k \star f \rightarrow f \mbox{ \rm as } k \rightarrow \infty$.
\item[(ii)]
If,  for $j \in {\mathbb{N}}$,
the function 
$v_j$ 
is defined by
\[
 v_j(x) = \max \{ 1 - \frac{\abs[x]}{j},0\}, \quad x \in {\mathbb{R}},
\]
then each 
$v_j$ 
is a positive definite function, and 
$v_j \nearrow 1$ 
uniformly on compacts as $ j \rightarrow \infty$.
\end{description}
It follows, see Remark~\ref{rem:Cv_p}, that for all $ j,k \in {\mathbb{N}}$
\begin{eqnarray}
\label{eq:approx}
  \norm[p,p]{v_j(\phi  \star u_k )} &  \le & \norm[p,p]{\phi  \star  u_k \ }
 \ \le \ \norm[p,p]{\phi} \ \le \   C_p.
\end{eqnarray}
By the conditions imposed in (i) the supports of $u_k$ and 
$\phi \star u_k$ are contained in the interval
$[\frac{1}{k}, \infty) \subset [0,\infty)$.
Thus, a bounded holomorphic extension to the lower half-plane 
\index{${\mathbb{C}}_{-}$}%
${\mathbb{C}}_{-}$
exists for each of them.
We note that, for those extensions and  for $ z \in {\mathbb{C}}_{-}$,  
\begin{eqnarray*}
(\phi \star u_k)\ \hat{}\ (z) 
& = & \int_{0}^{\infty}\ \phi \star u_k (x) e^{-izx} \ dx \\
& = & \hat{\phi}(z) \hat{u}_{k}(z).
\end{eqnarray*}
Hence,
\begin{eqnarray*}
\int_{-\infty}^{\infty}(\phi \star u_k)(t) e^{-t \lambda} \ dt 
& = & \int_{-\infty}^{\infty}(\phi \star u_k)(t) e^{-it\, (-i\lambda)} \ dt \\
& = &(\phi \star u_k)\ \hat{}\ (-i\lambda) \, = \,  \hat{\phi}(-i\lambda)%
 \hat{u}_{k}(-i\lambda)\\
& = & m(\lambda) \hat{u}_{k}(-i\lambda).
\end{eqnarray*}
Using again the conditions imposed on 
$u_k$
we see from
\[
\hat{u}_{k}(-i\lambda) =  \int_{-\infty}^{\infty} u_k(t) e^{-it\, (-i\lambda)} \ dt
= \int_{-\infty}^{\infty} u_k(t) e^{-t\lambda} \ dt,
\]
that
$\hat{u}_{k}(-i\lambda) \rightarrow 1 $
boundedly and uniformly on finite intervals in ${\mathbb{R}}_+$.
Now it is not hard to see that
\[
  \mathcal{S}_{j,k}f := \int^\infty_0\ (v_j \cdot (\phi  \star  u_k))(t)T_t f\ dt
\]
converges to 
$m(-A)f$ 
in 
$L^2$ 
when first
$j \rightarrow \infty$ 
and then 
$k \rightarrow \infty$.
In fact,
\begin{eqnarray*}
\lim_{k \rightarrow \infty}\ \lim_{j \rightarrow \infty} 
\mathcal{S}_{j,k}f & = & \lim_{k \rightarrow \infty}\ \lim_{j \rightarrow \infty} \int^\infty_0\ v_j(t) (\phi  \star  u_k)(t)T_t \
dt \\
& = & \lim_{k \rightarrow \infty}\ \lim_{j \rightarrow \infty} 
\int^\infty_0\ v_j (\phi  \star  u_k)(t)
\int_{- \infty}^0 e^{t \lambda}\ dP_{\lambda} f\ dt\\
& = & \lim_{k \rightarrow \infty}\ \lim_{j \rightarrow \infty} 
\int_{- \infty}^0 e^{t \lambda}\ \int^\infty_0\ v_j (\phi  \star  u_k)(t)\ dt\
dP_{\lambda} f \\
& = & \lim_{k \rightarrow \infty}\ \int_{- \infty}^0 \int^\infty_0\
e^{t \lambda}\
(\phi  \star  u_k)(t)\ dt\ dP_{\lambda} f \\
& = & \int_{- \infty}^0 \int^\infty_0\ e^{t \lambda}\ \phi(t)\ dt
dP_{\lambda} f \\
& = & \int_{- \infty}^0 m(-\lambda) dP_{\lambda} f.
\end{eqnarray*}
For all $ j,k \in {\mathbb{N}}$ we know, using the transference theorem
of Coifman and Weiss, our Corollary~\ref{cor:maptransfer}, that
(\ref{eq:approx}) implies
\[
  \norm[p]{  \mathcal{S}_{j,k}f } \le C_p \norm[p]{ f }.
\]
Because of this boundedness we obtain the existence of
$ \lim_{k \rightarrow \infty}\ \lim_{j \rightarrow \infty}  \mathcal{S}_{j,k}$
in the strong operator topology on the closure of 
$ \Lp{\Omega}{\mu} \cap L^2(\Omega,\mu) $, that is on all of 
$ \Lp{\Omega}{\mu}$.
The uniqueness of the limit on $ \Lp{\Omega}{\mu} \cap L^2(\Omega,\mu) $
shows that the limit must be the extension to $ \Lp{\Omega}{\mu}$ of the
restriction to
$ \Lp{\Omega}{\mu} \cap L^2(\Omega,\mu) $
of the spectrally defined
$m(-A)$.
\qed
In the case that $m$ is given as in Example~\ref{ex:complpow} it follows that
\[
  (-A)^{i \gamma} =m_\gamma (-A)
\]
extends to a bounded operator on 
$L^q(\Omega,\mu),\ 1<q<\infty$, 
of norm at most
\begin{eqnarray}
\label{eq:pnorm}
  \norm[q,q]{(-A)^{i \gamma}} \le C_1 [q + (q-1)^{-1}]\ (1 + \gamma)e^{\gamma \frac{\pi}{2}}\,.
\end{eqnarray}
Since 
\begin{eqnarray}
\label{eq:2norm}
 \norm[2,2]{(-A)^{i \gamma}} & = & \sup_{x\in{\mathbb{R}}} \abs[m_{\gamma}(x)] \
 = \  1,
\end{eqnarray}
it is possible to improve this by interpolation.
\begin{cor}
\label{cor:genpowers}
For some constant $C(p)$, depending only on $ p, \, 1<p<\infty:$
\[ 
\norm[p,p]{(-A)^{i\gamma}}\le C(p)[1+\abs[\gamma ]^{12}
]^{| \frac{1}{p}\ - \frac{1}{2}|} e^{\pi | \frac{1}{p}\ - \frac{1}{2}| \abs[ 
  \gamma ] }\,.
\]
\end{cor}
{\bf Proof:} We may assume $p \neq 2$.
\par 
The function 
\[
f: y \mapsto \frac{y}{y - \log y} , \qquad y \in [e,\infty)
\]
decreases monotonically to $1$.
Thus there exists at most one $\gamma_0 > e $ such that
\[
f(\gamma_0) \abs[ \frac{2}{p}-1] = 1.
\]
If a solution to this equation exists, then we define $\gamma_p$ equal
to it. Otherwise, i.e.\, if $p$ is such that
$ f(y)\abs[ \frac{2}{p}-1] < 1$ 
for all $y\geq e$,
then we let $\gamma_p = e$.
\par
If $\abs[\gamma] \leq \gamma_p $, then we estimate
\begin{eqnarray*}
\norm[p,p]{(-A)^{i \gamma}}
& \leq & c_1(p) (1+ \abs[\gamma]) e^{\frac{\pi}{2} \abs[\gamma]}\\
& \leq & c_1(p) (1+ \abs[\gamma]) e^{\frac{\pi}{2}  \abs[ \frac{2}{p} \,-1]
  \abs[\gamma]} e^{\frac{\pi}{2} (1-\abs[\frac{2}{p} \, -1])\gamma_p}\\
& \leq & c_2(p) (1+ \abs[\gamma]) e^{\frac{\pi}{2}  \abs[ \frac{2}{p} \, -1]
  \abs[\gamma]}\\
& \leq & c_3(p) (1+ \abs[\gamma])^{ \abs[ \frac{2}{p} \, -1]} e^{\pi
  \abs[ \frac{1}{p} \, - \frac{1}{2} ] \abs[\gamma]},
\end{eqnarray*}
which proves the corollary in this case.
\par
Otherwise, if $ \abs[\gamma] > \gamma_p $,
then we define $ \theta \in (0,1) $ by
\[
\theta = f(\abs[\gamma]) \abs[ \frac{2}{p}-1]. 
\]
Further we can define $q$ uniquely by
\[
\frac{1}{p}= \frac{\theta}{q} +(1-\theta)\frac{1}{2}.
\]
Then we have 
\[
1= f(\abs[\gamma]) \abs[ \frac{2}{q}-1],
\]
and it follows from this that
\begin{eqnarray*}
\abs[\gamma] \geq \frac{\abs[\gamma]}{\log \abs[\gamma]}& = & \left\{
\begin{array}{ll}
\displaystyle{\frac{q}{2}} \qquad &\mbox{ if }2 < q \\
\displaystyle{\frac{q}{2} \frac{1}{q-1}} \qquad & \mbox{ if } 1<q<2
\end{array}
\right.\\
& \geq & \frac{1}{4} (q + \frac{1}{q-1} ).
\end{eqnarray*}
By the Riesz-Thorin interpolation theorem we obtain from (\ref{eq:pnorm})
and (\ref{eq:2norm}):
\begin{eqnarray*}
\norm[p,p ]{(-A)^{i \gamma}}& \leq & C_{1}^{\theta}[q + \frac{1}{q-1}]^{\theta}
(1+ \abs[\gamma)]^{\theta} e^{\frac{\pi}{2} \abs[\gamma] \theta}\\
& \leq & C_2 \abs[\gamma]^{\theta}(1+\abs[\gamma])^{\theta}
e^{\frac{\pi}{2}  \abs[ \frac{2}{p} \, -1]\abs[\gamma]}e^{\frac{\pi}{2}  \abs[
    \frac{2}{p} \, -1]\abs[\gamma](f(\abs[\gamma])-1)}\\
& = & C_2 \abs[\gamma]^{\theta}(1+\abs[\gamma])^{\theta}
e^{\frac{\pi}{2}  \abs[ \frac{2}{p} \, -1]\abs[\gamma]}e^{\frac{\pi}{2} 
\abs[\frac{2}{p} \, -1] \log \abs[\gamma]f(\abs[\gamma])  }\\
& \leq & C_2(1+\abs[\gamma])^{(1+1+\frac{\pi}{2})\theta}e^{\frac{\pi}{2}
  \abs[ \frac{2}{p} \, -1]\abs[\gamma]}\\
& \leq & C_3(1+\abs[\gamma])^{(4+\pi) f(\abs[\gamma])\abs[ \frac{1}{p}-\frac{1}{2}]}
e^{\frac{\pi}{2} \abs[ \frac{2}{p} \, -1]\abs[\gamma]}.
\end{eqnarray*}
Here we used that by the definition of $f$:
\[
y\ (f(y)-1) = y\  \frac{\log y}{y - \log y} = \log y\  f(y).
\]
Since $f(\abs[\gamma]) \leq \frac{e}{e-1}$, we have proved the assertion of
the corollary.
\qed
\par
In his above cited paper~\cite{Cowling_83} Cowling has a better estimate, 
but I could not figure out his computation. Furthermore he deduces, 
and we refer 
the reader to this deduction, the following theorem: 
\begin{thm}[Cowling]
\label{thm:univmultiplier}
\index{theorem!multiplier!universal}%
Let 
$m$ 
be bounded, holomorphic in 
$\Gamma_{\psi}$, 
where 
$ 0 < \psi \le \frac{\pi}{2}$.
If 
$p$ 
is such that 
$\abs[ \frac{1}{p} - \frac{1}{2} ] < \frac{\psi}{\pi}$,
then for some constant $C$,  not depending on $ p  $,
\[
 \norm[p]{ m(-A)f} \le C 
\left[ \frac{\psi}{\pi} - \abs[ \frac{1}{p} - \frac{1}{2} ]
 \right]^{-\frac{5}{2}} \norm[\infty]{m }\,.
\]
\end{thm}
 \par
Now we consider the family of functions
\[
  m_\theta: \lambda \mapsto \exp (e^{\theta}\lambda)\ \quad \theta \in \Gamma.
\]
Then 
$m_0 (tA) = e^{tA} = T_t$, 
and proving $L^p$-boundedness of the, on $L^2$ spectrally defined operators
\[
 m_\theta (tA) = e^{te^{i\theta}A} =: T_{te^{i\theta}}, \qquad 
\abs[ \theta ] < \eta,
\]
will give an extension of the semigroup 
$ (T_t)_{t \geq 0} $
to the complex sector 
$\Gamma_\eta$. Here, $\eta$ is a possibly $p$-dependent constant,
ideally we should find the largest possible one. 
The boundedness and the weak analyticity on 
$L^2(\Omega,\mu)$
imply the  analyticity of this extension on this sector. 
\par
For this purpose the
operators 
$(-A)^{i\gamma}, \ \gamma \in {\mathbb{R}}$ 
are quite useful. By inverting 
\index{transform!Mellin}%
Mellin transforms we obtain for $x\leq 0, \, t
\geq 0:$
 \begin{eqnarray}
\label{eq:Melinv}
\exp(te^{i\theta}x) - \exp (tx) & =& \frac{1}{2\pi}\ \int_{\mathbb{R}}\ \left[ e^{-\gamma \theta}
  -1 \right] \Gamma (-i\gamma)t^{ i \gamma}  (-x)^{i\gamma}d \gamma\,.
     \end{eqnarray}
 Compute for this:
\begin{eqnarray*}
\int^\infty_0\ \exp (-e^{i\theta}\lambda) \lambda^{-i\gamma} \frac{d\lambda}{\lambda}
&  = &\lim_{\varepsilon \searrow 0}\ \int^\infty_0\ \exp(-e^{i\theta}\lambda)
  \lambda^{\varepsilon -i\gamma -1}d\lambda \\
& = & \lim_{\varepsilon \searrow 0}\ \int^\infty_0\ \exp (-\lambda)
  e^{-i\theta (\varepsilon - i \gamma)} \lambda^{\varepsilon - i \gamma - 1}d \lambda\\
& = & e^{-\gamma \theta} \int^\infty_0\
\lambda^{-i\gamma}e^{-\lambda}\frac{d\lambda}{\lambda}\\
& = & 
  e^{-\gamma \theta} \Gamma (-i\gamma)\,.
\end{eqnarray*}
One can not directly make use of this last formula, since 
$\abs[\Gamma(-i\gamma)] \sim 
\frac{1}{ \abs[\gamma ]}$ 
for 
$\gamma $ 
close to zero, and so for $x>0$ an Integral
\[
  \frac{1}{2\pi}\ \int_{\mathbb{R}}\ e^{-\gamma \theta} \Gamma (-i\gamma)t^{i\gamma}
  x^{i\gamma}d\gamma \simeq e^{-te^{i\theta}x}
\]
does not converge. But in equation~(\ref{eq:Melinv}) it can be used that for
some
$C>0$:
\begin{eqnarray}
\label{eq:Melest}
\abs[ e^{- \gamma \theta} - 1] \abs[ \Gamma(-i \gamma)] 
& \leq & e^{(\abs[\theta] - \frac{\pi}{2}) \abs[\gamma]}.
\end{eqnarray}
This is used by Cowling to prove the following corollary due to 
\index{Stein}%
Stein, c.f.\ \cite{Stein_70a} chap III, Theorem 1.
\index{theorem!Stein}%
\begin{cor}[Stein]
\label{cor:analyticextension}
\index{analytic!extension}%
If
$0 \leq \theta \leq \frac{\pi}{2}- \pi\abs[\frac{1}{p} \ - \frac{1}{2}]$,
then the semigroup
$T_t$ extends to an analytic semigroup on $\Gamma_{\theta}$.
\end{cor}
{\bf Sketch of a Proof:}
We may assume that $1<p<\infty$, since otherwise the
angle $\theta$ has to be zero.
\par  
The $L^2$-spectral theory shows that the semigroup has
a unique analytic extension on $L^2(\Omega,\mu)$, even to the cone
$\Gamma_{\frac{\pi}{2}}$.
We further obtain, for $z= t e^{i\theta}$:
\[
e^{t e^{i\theta} A} - e^{tA} =   \frac{1}{2\pi}\ \int_{\mathbb{R}}\ \left[ e^{-\gamma \theta}
  -1 \right] \Gamma (-i\gamma)t^{ i \gamma}  (-A)^{i\gamma}d \gamma\, , 
\]
the integral being convergent in the strong operator topology, still on $L^2$.
\par
Now, if $f \in \Lp{\Omega}{\mu}$ and $g \in \Lq{\Omega}{\mu}$
$(\frac{1}{p} +\frac{1}{q}=1)$
are simple functions, then they are in $L^2$.
From the above inequality~(\ref{eq:Melest}) and Corollary~\ref{cor:genpowers}
we infer that
for $\abs[\psi] < \theta$:  
\begin{eqnarray*}
\abs[<T_{te^{i\psi}}f,g>] & \leq & \frac{C(p)}{2\pi}\int_{\mathbb{R}} %
[1+\abs[\gamma ]^{12} ]^{| \frac{1}{p}\ - \frac{1}{2}|} %
e^{\pi | \frac{1}{p}\ - \frac{1}{2}| \abs[ \gamma ] } %
e^{(\abs[\psi] - \frac{\pi}{2}) \abs[\gamma]} d\gamma \norm[p]{f}\norm[q]{g}\\
& \leq &  \frac{C(p)}{2\pi}\int_{\mathbb{R}} %
[1+\abs[\gamma ]^{12} ]^{| \frac{1}{p}\ - \frac{1}{2}|} %
e^{-(\theta - \abs[\psi]) \abs[\gamma]} d\gamma \norm[p]{f}\norm[q]{g}\\
& \leq & C_{\psi} \norm[p]{f}\norm[q]{g}.
\end{eqnarray*}
\par
If general elements $f \in \Lp{\Omega}{\mu}$ and 
$g \in \Lq{\Omega}{\mu}$ are given, then they may be approximated
by sequences of simple functions, $(f_k)_k$ and $(g_k)_k$
in the respective norms. The above estimate then shows that the 
sequence of analytic functions
$z \mapsto <T_z f_k , g_k>$ converges locally uniformly inside $\Gamma_{\theta}$.
\qed
 \par
A formula similar to the above~(\ref{eq:Melinv}), together with the appropriate estimate:
\begin{eqnarray}
\label{eq:Melinv2}%
\!\!\!\!\!e^{t e^{i\theta} x} - \frac{1}{t}\int_{0}^{t} \!e^{sx}  ds 
\!& = & \! \frac{1}{2\pi} \int_{\mathbb{R}} \left[ e^{-\gamma \theta}
  -\frac{1}{1+i\gamma} \right] \Gamma (-i\gamma)t^{ i \gamma}
(-tx)^{i\gamma}d \gamma,
\end{eqnarray} 
\begin{eqnarray}
\label{eq:Melest2}%
\abs[e^{(-\theta \gamma)} - \frac{1}{1+i\gamma}]\abs[ \Gamma(-i\gamma)]
& \leq & C e^{(\abs[\theta] - \frac{\pi}{2}) \abs[\gamma]},
\end{eqnarray}
can be used prove a maximal theorem and an abstract
\index{theorem!non-tangential convergence}%
non-tangential convergence theorem. 
(Cowling, in his paper proves an
interesting extension of it):
\begin{thm}[Stein]
\label{thm:aeconv}
If 
$ 0 \leq \theta < \frac{\pi}{2}(1- \abs[\frac{2}{p} \ - 1])$,
then
for some constant $C_{p}>0$ for all $f \in \Lp{\Omega}{\mu}$:
\[
\norm[p]{\sup \left\{\abs[T_z f] \, : \, \abs[z]<1, z \in
    \Gamma_{\theta}\right\}} \leq C_{p} \norm[p]{f},
\]
and 
\[
T_z f \rightarrow f \quad \mu \mbox{-almost everywhere as } z \rightarrow 0 \mbox{ in }
\overline{\Gamma_{\theta}}.
\]
\end{thm}
{\bf Proof:}
We estimate pointwise
\begin{eqnarray*}
\lefteqn{\sup \left\{\abs[T_z f]\, :\, \abs[z]<1,z \in \Gamma_{\theta} \right\}}\\
& \leq & \sup_{\abs[\psi]<\theta, \ 0 \leq t <1}%
\left| \frac{1}{2\pi}\ \int_{\mathbb{R}}\ 
\left[ e^{-\gamma \psi}-\frac{1}{1+i\gamma} \right] \Gamma (-i\gamma)t^{ i
  \gamma}  (-A)^{i\gamma}f\, d \gamma\ \right| + \\%
& & \qquad \qquad \qquad \sup_{0\leq t <1} \abs[\frac{1}{t}\int_{0}^{t} T_sf ds].
\end{eqnarray*}
Thus the norm of the first term on the right hand side may be estimated using the above inequality~(\ref{eq:Melest2})and
Corollary~\ref{cor:genpowers}:
\begin{eqnarray*}
\lefteqn{\norm[p]{ \sup_{\abs[\psi]<\theta, \ 0 \leq t <1} \frac{1}{2\pi}\ \int_{\mathbb{R}}\
\abs[ e^{-\gamma \psi}-\frac{1}{1+i\gamma} ] \abs[\Gamma (-i\gamma)t^{ i
  \gamma}]  \abs[(-A)^{i\gamma}f]\ d \gamma }}\\
& \leq & C \norm[p]{\frac{1}{2\pi}\ \int_{\mathbb{R}}\
e^{(\abs[\theta] - \frac{\pi}{2})\abs[\gamma]}  \abs[(-A)^{i\gamma}f]\ d
\gamma } \\
& \leq & C \frac{1}{2\pi}\ \int_{\mathbb{R}}\
e^{(\abs[\theta] - \frac{\pi}{2})\abs[\gamma]}  \norm[p]{(-A)^{i\gamma}f}\ d
\gamma \\
& \leq & C \frac{1}{2\pi}\ \int_{\mathbb{R}}\
C(p)[1+\abs[\gamma]^3 \log^2 \abs[\gamma]]^{\abs[\frac{1}{p} -\frac{1}{2}]}e^{(\abs[\theta] -
  \frac{\pi}{2})\abs[\gamma] + \pi{\abs[\frac{1}{p}
    -\frac{1}{2}]\abs[\gamma]}}\ d\gamma \norm[p]{f}\\
& \leq & C^{1}_{p} \norm[p]{f}.
\end{eqnarray*}
Together with Theorem~\ref{thm:dunschw} this establishes the maximal inequality of the theorem.
\par
The assertion on the almost everywhere convergence relies on an application of
Banach's principle.
\begin{thm}[Banach's principle]
\label{thm:banachprinciple}
\index{Banach's principle}
Let $(T_n)_{n\in {\mathbb{N}}}$
be a sequence of operators on 
$\Lp{\Omega}{\mu}$,
$1\leq p < \infty$,
such that for some constant
$C>0$,
for all 
$f \in \Lp{\Omega}{\mu}:$
\[
\norm[p]{\sup_{n \in \mathbb{N}}\abs[T_n f]} \leq C \norm[p]{f}.
\]
If for all elements 
$f$ in a dense subspace $E$ of $\Lp{\Omega}{\mu}$\\
\[
\lim_{n \rightarrow \infty} T_n f 
\]
exists pointwise almost everywhere,
then for all $f \in \Lp{\Omega}{\mu}$
\[
\lim_{n \rightarrow \infty} T_n f
\]
exists pointwise almost everywhere.
\end{thm}
{\bf Proof:}
We are given
$f \in  \Lp{\Omega}{\mu}$.
For
$\epsilon >0$ we write 
$f=f_1+f_2$ with $f_1 \in E$ and $\norm[p]{f_2}<\epsilon$.\\
Then for $\mu$ almost all $\omega \in \Omega$
\begin{eqnarray*}
0 & \leq & {\lim_{n \rightarrow \infty} \sup}{\mbox{ \rm Re }T_n f(\omega)} 
- {\lim_{n \rightarrow \infty} \inf}{\mbox{ \rm Re }T_n f(\omega)}\\
 & \leq & {\lim_{n \rightarrow \infty} \sup}{\mbox{ \rm Re }T_n f_1(\omega)} 
- {\lim_{n \rightarrow \infty} \inf}{\mbox{ \rm Re }T_n f_1(\omega)}\\
& & \qquad \qquad \qquad + {\lim_{n \rightarrow \infty} \sup}{\mbox{ \rm Re }T_n f_2(\omega)} 
- {\lim_{n \rightarrow \infty} \inf}{\mbox{ \rm Re }T_n f_2(\omega)}\\
 & \leq & 0 + 2 \sup_{n \in {\mathbb{N}}} \abs[T_nf_2](\omega). 
\end{eqnarray*}
Hence
\begin{eqnarray*}
\norm[p]{ {\lim_{n \rightarrow \infty} \sup}{\mbox{ \rm Re }T_n f(.)} 
- {\lim_{n \rightarrow \infty} \inf}{\mbox{ \rm Re }T_n f(.)}}
& \leq &2C\epsilon.
\end{eqnarray*}
Since $\epsilon > 0$ is arbitrary this proves that 
\[
{\lim_{n \rightarrow \infty} \sup}{\mbox{ \rm Re }T_n f(\omega)} 
- {\lim_{n \rightarrow \infty} \inf}{\mbox{ \rm Re }T_n f(\omega)}
\]
vanishes almost everywhere.
An analogous argumentation for the imaginary parts then proves the almost
everywhere convergence of the sequence
$(T_n f)_{n \in \mathbb{N}}$.
\qed
The subspace $E$ is exhibited using the analyticity of the semigroup in
question. For this denote
$\Gamma_{p}$ the cone to which the semigroup extends analytically:
\begin{lem}
Given the conditions of Theorem~\ref{thm:aeconv}. For 
$f \in \Lp{\Omega}{\mu}$ and 
$z \in \Gamma_{p}$ it is then  possible to redefine
$T_zf$ on a set of measure zero, such that 
for all $\omega \in \Omega$ 
\[
 z \mapsto T_zf(\omega) \mbox{ is analytic on } \Gamma_p.
\]
\end{lem}
{\bf Proof:}
The weakly analytic map 
$ \Phi : z \mapsto T_z\ f$, from 
$\Gamma_p$ to $\Lp{\Omega}{\mu}$, 
is strongly analytic (See e.g.\ chap.5 sect.3 of~\cite{Yosida_78}). 
Given a closed ball $B(z,r)$ of radius $r$ with centre $z$ contained in $
\Gamma_p$, $\Phi$ has an, in the interior  $B(z,r)^{\circ}$ of the ball,
convergent 
expansion as a power series
\begin{eqnarray}
\label{eq:powerser1}
T_{\zeta} \ f & = & \sum_{n=0}^{\infty} h_n \, (\zeta -z)^{n} \mbox{ with } \\
\label{eq:powerser2}
\sum_{n=0}^{\infty} \norm[p]{h_n} \rho^n & < & \infty \mbox{ for all }
\rho < r.
\end{eqnarray}
For $n \in \mathbb{N}$ let $ H_n : \Omega \rightarrow \mathbb{C}$
be a function representing $ h_n $ in  $\Lp{\Omega}{\mu}$. From the estimate
in
(\ref{eq:powerser2}) it follows that for $\rho < r$
$ \sum_{n=0}^{\infty} \abs[{H_n}(\omega)] \rho^n  <  \infty $ almost
everywhere on $\Omega$.
In fact, the triangle inequality ensures that for all $k \in \mathbb{N}$
\[
\norm[p]{\sum_{n=0}^{k} \abs[H_n] \rho^n} \leq \sum_{n=0}^{\infty} %
\norm[p]{H_n}\rho^n < \infty.
\]
The monotone convergence theorem then says that
$\left(\sum_{n=0}^{k} \abs[H_n] \rho^n \right)^p$
converges almost everywhere to an integrable function.
Where this function is finite there 
$\sum_{n=0}^{k} \abs[H_n] \rho^n$
must converge to a finite value too.
\par
Choosing a sequence $\rho_k \nearrow r$ we find a set $A_z$
of measure zero such that for $\omega \in \Omega \setminus A_z$
the series 
$\sum_{n=0}^{\infty} {H_n}(\omega)(\zeta -z)^{n}$
converges absolutely for all $\zeta \in B(z,r)^{\circ}$.
On this ball let
\[
\begin{array}{rcl}
F_z(\zeta,\omega) & = & \left\{  \begin{array}{lcl} %
\sum_{n=0}^{\infty} {H_n}(\omega)(\zeta -z)^{n} & \quad & \omega \notin A_z \\
0 & \quad & \omega \in A_z.
\end{array}  \right.
\end{array}
\]
Now the (open) cone 
$\Gamma_p $
can be covered by a countable union of balls
$B(z_k,\frac{r_k}{2})^{\circ}$ with $B(z_k,r_k) \subset \Gamma_p$.
If two such open balls have non-void intersection, say
$B(z,\frac{r}{2})^{\circ} \cap B(z',\frac{r'}{2})^{\circ} \not= \emptyset$,
then, for almost all 
$\omega \in \Omega$,
the corresponding analytic functions
$F_z(\zeta,\omega)=\sum_{n=0}^{\infty}{H_n}(\omega)(\zeta -z)^{n} $
and 
$F_{z'}(\zeta,\omega)=\sum_{n=0}^{\infty} {H'_n}(\omega)(\zeta -z')^{n} $
coincide on this intersection.
To see this we may assume $r' \leq r$, then $z' \in B(z,r)^{\circ}$,
and the series representing $F_z(\zeta,.)$ converges outside $A_z$
absolutely and uniformly in a neighbourhood of $z'$.
Moreover, for $\omega \notin A_z \cup A_{z'}$:
\[
H'_0 (\omega) = \sum_{n=0}^{\infty}{H_n}(\omega)(z' -z)^{n}. 
\]
We can even differentiate (k times) to the result, that
\[
\begin{array}{rcccl}
\displaystyle{H'_k (\omega)} &=&\displaystyle{\frac{1}{k!} (\frac{d}{d\zeta})^k %
\sum_{n=0}^{\infty}{H'_n}(\omega)(\zeta -z')^{n} |_{\zeta =z'}}&&\\%
&=&\displaystyle{\frac{1}{k!} (\frac{d}{d\zeta})^k %
\sum_{n=0}^{\infty}{H_n}(\omega)(\zeta -z)^{n}|_{\zeta =z'}}
&=&\displaystyle{\sum_{n=k}^{\infty}{H_n}(\omega)%
\left( {}^{\textstyle n}_{\textstyle k} \right)(z' -z)^{n-k}.}
\end{array}
\]
The last sum here is still absolutely convergent, and
\begin{eqnarray*}
\sum_{k=0}^{\infty} {H'_k}(\omega)(\zeta -z')^{k} %
&=&\sum_{k=0}^{\infty}\sum_{n=k}^{\infty}{H_n}(\omega)%
\left( {}^{\textstyle n}_{\textstyle k} \right)(z' -z)^{n-k}(\zeta -z')^{k}\\
&=&\sum_{n=0}^{\infty}{H_n}(\omega) \sum_{k=0}^{n} %
\left( {}^{\textstyle n}_{\textstyle k} \right)(z' -z)^{n-k}(\zeta -z')^{k}\\
&=&\sum_{n=0}^{\infty}{H_n}(\omega) (\zeta -z)^n.
\end{eqnarray*}
Now we define $T_{\zeta}f(\omega) = 0 $ for $\omega \in \bigcup_{k} A_{z_k}$
and
\[ T_{\zeta}f(\omega) = F_{z_k}(\zeta, \omega) \quad \mbox{for} \quad %
\omega \in \Omega \setminus \bigcup_{k} A_{z_k},
\]
where $z_k$ is any element of our selection such that only
$\zeta \in B(z_{k},\frac{r_k}{2})^{\circ}$.
\qed


We may now continue with the {\bf proof of Theorem~\ref{thm:aeconv}:}\\
From the strong continuity of the semigroup 
$(T_t)_{t \geq 0}$
at zero we see that 
$E=\left\{ \, T_s \ f \, : \, f\in \Lp{\Omega}{\mu}, \, s > 0 \, \right\}$
is a dense linear subspace and for a sequence $z_n \in \Gamma_{\theta},
\, z_n \rightarrow 0$ we have for one of its elements
$T_sf$ outside the set where $ z \mapsto T_zf$ is not analytic:
\begin{eqnarray*}
\lim_{n \rightarrow \infty} T_{z_n}T_sf(\omega) 
& = & \lim_{n \rightarrow \infty} T_{z_n +s}f(\omega) \\
& = & T_s f(\omega).
\end{eqnarray*}
\qed
\begin{rem} \rm \hspace*{\fill}~\\[-2em]
\label{rem:angleanalytic}
\begin{description}
\item[(i)~~]
In the submarkovian case the angle of the cone of 
\index{angle!analyticity}%
analyticity
given here is not optimal. 
\index{Liskevich}%
Liskevich and 
\index{Perelmuter}%
Perelmuter \cite{LiskPere_95}
obtain 
$\frac{\pi}{2}- \arctan \frac{\abs[p-2]}{2\sqrt{p-1}}$.
\item[(ii)~]
If one no longer assumes that the semigroup is contractive on the whole scale
of 
$L^p$-spaces, then the above presented methods still apply.
\par
If we assume that our strongly continuous semigroup 
fulfils, for some $1<r<2$ the following set of conditions:
\begin{enumerate}
\item[i')~] 
For those $p$ such that 
$r \leq p \leq r'$ each operator $T_t$ is a sub-positive contraction on 
$\Lp{\Omega}{\mu}$ 
($ \frac{1}{r} + \frac{1}{r'}=1 $),
\item[ii)~]
Selfadjointness on $L^2$:
$  T_t^\ast  =  T_t \; \mbox{\rm on}\; L^2(\Omega,\mu) $,
\item[iii)] 
$T_0  =  \id$,
\end{enumerate}
then, for 
$p<2$,
the angle of analyticity would be at least
$\frac{\pi}{2} \frac{r}{p} \frac{2-p}{2-r}$.
Of course we obtain the maximal theorem for the cone of this angle and the
other consequences too.
\end{description}
\end{rem}
\addcontentsline{toc}{chapter}{Bibliography}
\bibliographystyle{plain}

\addcontentsline{toc}{chapter}{List of Figures}
\listoffigures

\addcontentsline{toc}{chapter}{Index}
\printindex

{\bf Author's address}\\
Gero Fendler\\
Finstertal 16\\
69514 Laudenbach\\
Germany\\
email: gero.fendler<at>univie.ac.at
\end{document}